\documentclass[twoside,10pt]{article}

\usepackage[english]{babel}
\usepackage[latin1]{inputenc}
\usepackage{amsmath,amsfonts,amssymb}
\usepackage{graphicx}
\usepackage{color}
\usepackage{psfrag}

\def\R{\mathbb{R}}
\def\N{\mathbb{N}}
\def\u{\textit {\textbf u}}

\newcommand\dps{\displaystyle }
\newtheorem{remark}{Remark}

\def\qed{\relax
     \ifmmode
       ~\hfill\Box
     \else
        \unskip\nobreak ~\hfill$\Box$
      \fi \par}

\newtheorem{algorithm}{Algorithm}

\begin{document}
\title{Symmetric parareal algorithms for Hamiltonian systems}
\author{X. Dai$^{1,2}$, C. Le Bris$^{3}$, F. Legoll$^{3}$ and Y. Maday$^{1,4}$ \\
{\footnotesize $^1$ UPMC Univ. Paris 06, UMR 7598, Laboratoire
  J.-L. Lions,
} \\ {\footnotesize 
Bo\^{\i}te courrier 187, 75252 Paris Cedex 05, France}
\\
{\footnotesize $^2$ LSEC, Institute of Computational Mathematics and
  Scientific/Engineering Computing,
} \\ {\footnotesize 
 Academy of Mathematics and Systems
  Science, Chinese Academy of Sciences, Beijing 100190, China}
\\
{\footnotesize $^3$ \'Ecole Nationale des Ponts et
Chauss\'ees, 6 et 8 avenue Blaise Pascal, 
} \\ {\footnotesize 
77455 Marne-La-Vall\'ee Cedex 2, France}\\
{\footnotesize and
} 
\\
{\footnotesize  INRIA Rocquencourt, MICMAC team-project, Domaine de
  Voluceau,
} \\ {\footnotesize 
B.P. 105, 78153 Le Chesnay Cedex, France}
\\ {\footnotesize   $^4$ Division of Applied Mathematics, Brown University,
  Providence, RI, USA} 
}
\date{\today}

\maketitle

\begin{abstract}
The parareal in time algorithm allows to efficiently use parallel computing
for the simulation of time-dependent problems. It is based on a
decomposition of the time interval into subintervals, and on
a predictor-corrector strategy, where the propagations over
each subinterval for the corrector stage are concurrently performed  on
the processors.

In this article, we are concerned with the long time integration of
Hamiltonian systems. Geometric, structure-preserving integrators are
preferably employed for such systems because they show interesting
numerical properties, in particular excellent preservation of the total
energy of the system. Using a symmetrization procedure and/or a
(possibly also symmetric) projection step, we introduce here several
variants of the original plain parareal in time
algorithm~\cite{lions-maday-turinici-01,baffico-bernard-maday-02,bal-maday-02}
that are better adapted to the Hamiltonian context. 
These variants are compatible with the geometric structure of the exact dynamics, and are easy to implement.
 
Numerical tests on several model systems illustrate the remarkable
properties of the proposed parareal integrators over long integration times. Some formal
elements of understanding are also provided. 

\medskip


\medskip

AMS: 65L05, 65P10, 65Y05
\end{abstract}

\section{Introduction}

Increasingly intensive computations now become possible thanks to the
improvement of both the efficiency and the
clock rate of processors, the interprocessor's connections and the
access to the different levels of memory. In addition, parallel
computing platforms, which allow many 
processors to work concurrently, also become available. This second
feature can only be useful if the problem to
be solved can be decomposed into a series of independent tasks, each
of them being assigned to one of the processors. 

The design of efficient algorithms for parallel architectures is the
subject of intense current research. In the case 
of models governed by partial differential equations, most of --- if
not all --- the contributions of the last three decades perform
domain decomposition. We refer to~\cite{Quarteroni,Widlund}
for a review on recent advances, and also to the proceedings of the Domain
Decomposition Method meetings (see {\tt www.ddm.org}) for various achievements. When the problem is
time-dependent, or when the problem is solely governed by a system of
ordinary differential equations, relatively few contributions are
available. We refer e.g. to the book of K. Burrage~\cite{Bur:book} for
a synthetic approach on the subject (see also~\cite{burrage_review}). In
this book, the various 
techniques of parallel in time algorithms are classified into three
categories: (i) parallelism {\em across the system}, (ii) parallelism
{\em across the method} and (iii) parallelism {\em across time}. 

The parareal-in-time method, our focus here,
  was first proposed in the work of
J.-L. Lions, Y. Maday and G. Turinici in
2001~\cite{lions-maday-turinici-01}. It belongs to the third category
where parallelism is achieved by 
  breaking up the integration interval into sub-intervals and
  concurrently solving over each sub-interval. The obvious difficulty is
  to provide the correct initial value for the integration over each of
  these sub-intervals. Most of the techniques in the third category are
  {\it multishooting techniques}, starting from the precursory work by
  A. Bellen and M. Zennaro~\cite{BZ89}. This next led to the {\it waveform
  relaxation methods} introduced by E. Lelarasmee, A.E. Ruehli and
  A.L. Sangiovanni-Vincentelli~\cite{LRSV82}, and to the {\it multigrid
  approaches} introduced by W. Hackbush~\cite{Hack}. See
  also~\cite{chartier93}. 

As it has been
explained by M. Gander and S. Vandewalle in~\cite{gander-vandewalle-07},
the parareal in time algorithm can be interpreted both as a type of multishooting technique and
also as a type of multigrid approach, even though the bottom line of the approach is closer to that of spatial domain decomposition with coarse
grid preconditioners.  

It is intuitively clear why so few contributions propose
parallel-in-time algorithms:  
time-dependent problems are intrinsically sequential. On the other hand, the development of parallel
computing provides computational opportunities that exceed the needs of 
parallelization in space. This motivates the development of efficient parallel-in-time 
approaches.

The parareal in time algorithm is based on the decomposition of the
temporal domain into several temporal subdomains and the  combination of a
coarse and a fine solution procedure. The name ``parareal'' has been adopted to
indicate that the algorithm has, in principle, the potential to so
efficiently speed up the simulation process that real time approximation
of the solution of the problem becomes plausible. Since the original
work~\cite{lions-maday-turinici-01}, where the convergence is proven to
be linear, the parareal algorithm has received
some attention and new developments have been
proposed. In~\cite{bal-maday-02}, G. Bal and Y. Maday
have provided a new interpretation of the scheme
as a predictor-corrector algorithm (see also L. Baffico et
al.~\cite{baffico-bernard-maday-02}). The scheme  
involves a prediction step based on a coarse approximation for a
global propagation of the system and a correction step computed in
parallel and based on a fine approximation. G. Bal (in~\cite{bal-03})
and E.M. Ronquist and G. Staff (in~\cite{staff-ronquist-03}) next
provided some analysis on the convergence and 
the stability of the scheme. In~\cite{gander-vandewalle-05},
M.J. Gander and  S. Vandewalle proved a superlinear convergence of the
algorithm 
when used on bounded time intervals, and a linear convergence on
unbounded time intervals. 

In the past few years, the parareal algorithm has been
successfully applied to various types of problems (see~\cite{maday-08} for a review):
nonlinear PDEs~\cite{bal-maday-02}, control
problems~\cite{maday-turinici-02}, quantum control
problems~\cite{maday-turinici-03}, Navier Stokes
equations~\cite{ficher-hecht-maday-05}, structural
dynamics~\cite{farhat06}, reservoir simulation~\cite{garido}.  

Although the plain parareal algorithm has proved to be efficient
for many  time-dependent problems, it also has some drawbacks in some
specific cases. This is for instance the case for molecular dynamics
simulations (as pointed out in~\cite{baffico-bernard-maday-02}), or,
more generally, for 
Hamiltonian problems. The exact flow of the system then enjoys many specific {\em geometrical} properties
(symplecticity, possibly time-reversibility, \ldots). Some
quantities are preserved along the trajectory: the Hamiltonian (e.g. the
total energy of the system), and, in some cases, other quantities such
as the angular 
momentum, \ldots See the textbooks~\cite{hlw,ss-calvo,leim_reich} for a systematic
introduction to numerical integration techniques for Hamiltonian systems. 
Symplectic and symmetric integrators are known to be suitable 
integrators for Hamiltonian systems. It turns out that, even in the case when
the parareal algorithm is based on coarse and fine integrators that
enjoy adequate geometrical properties (such as symplectic or
symmetric integrators), the global parareal algorithm itself does not enjoy any of  these
properties. Consequently, the long time properties of the numerical flow
(including {\em e.g.} energy preservation)
are not as good as  expected. In this article, our aim is to 
design a parareal scheme that preserves some geometrical structure of
the exact dynamics. 

Our article is organized as follows. 
In Section~\ref{sec:plain_para}, we briefly recall the parareal
algorithm, as presented
in~\cite{baffico-bernard-maday-02,bal-maday-02}. We next discuss, in
Section~\ref{sec:para4ham_gene}, general, commonly used 
numerical schemes for Hamiltonian problems. We show the deficiencies of
the parareal algorithm on such problems. 
In Section~\ref{sec:sym}, we develop a {\em symmetric} version of the
parareal algorithm, in a sense made precise at the end of
Section~\ref{sec:design_parasym}. As an
alternative to 
symmetrization, we explore in Section~\ref{sec:proj} the idea of
projecting the trajectory onto the constant energy manifold.  Next, in
Section~\ref{sec:sym_proj}, we couple a projection step with the
symmetric algorithm developed in Section~\ref{sec:sym}, while keeping
the overall scheme symmetric. Combining these two ideas, symmetry and
projection, yields the most 
efficient algorithms. The performances of all the
schemes proposed in these sections are illustrated by numerical
simulations on two low-dimensional systems, namely the harmonic
oscillator and the two-dimensional Kepler problem. 
The computational complexities of the different parallel integrators are
analyzed in details in each of the respective sections.  
In Section~\ref{sec:outer}, we consider a test case in a higher
dimensional phase space, namely the simulation of the solar
system (see also~\cite{tremaine} for the derivation of algorithms
specific to this test case, involving some parallel computations). We
demonstrate that the efficiency and the qualitative properties observed
in the previous test cases  
carry over to  this more challenging example. Using the symmetric
parareal scheme with symmetric 
projection developed in Section~\ref{sec:sym_proj} (see
Algorithm~\ref{algorithm-sps-symm-proj}), we obtain a speed-up 
of more than 60 with respect to a fully sequential
computation (provided a sufficient number of processors are
available), for an equal accuracy. Section~\ref{sec:conclusion}
summarizes our conclusions.  

\section{The  plain parareal algorithm}
\label{sec:plain_para}

In this section, we review the plain parareal method for a
time-dependent problem in a general setting. Consider $u$ solution to the Cauchy problem
\begin{equation}
\label{eq:problem-1order}
\frac{\partial u }{\partial t} + f(u) = 0, \ t \in [0, T],
\quad \hbox{supplied with the initial condition }
u(0) = u_0.
\end{equation}
We assume standard appropriate conditions on $f$ ensuring existence, uniqueness
and stability with respect to perturbations, of the solution $u$ to this
problem. Let $\mathcal{E}$ 
be the exact propagator, defined by $\mathcal{E}_\tau(u_0) =
u(\tau)$, where $u(\tau)$ denotes the solution at time $\tau$ of
problem~\eqref{eq:problem-1order}. 

In the sequel, for reference, we approximate the exact propagator $\mathcal{E}$
by using an accurate numerical scheme. We can use any of the classical
one-step schemes (explicit or implicit Euler schemes for simple
problems, velocity Verlet scheme in the case of Hamiltonian dynamics)
with a sufficiently small time step $\delta t$. 
The associated discrete
propagator is denoted as $\mathcal{F}$, with no reference to the time
step $\delta t$, in order not to make the notation too heavy.
$\mathcal{F}_{\tau}(u_0)$ is thus an approximation of 
$\mathcal{E}_{\tau}(u_0)$. Assuming that the approximation $u_n$ of
$u(T_n)$ is known at some time $T_n$, the approximation of $u(T_{n+1})$
at a latter time $T_{n+1}$ is computed by  
performing $(T_{n+1} - T_{n})/\delta t$ steps of the fine scheme, that
is denoted, with the previous notations, 
$$
u_{n+1} = \mathcal{F}_{T_{n+1} - T_{n}}(u_n).
$$
In the sequel, we will
consider the dynamics
$$
\dot{q} = M^{-1} p, 
\quad
\dot{p} = -\nabla V(q),
\quad
(q,p) \in \R^d \times \R^d,
$$
where $M$ is a diagonal matrix, and $V$ is a smooth scalar function depending
on~$q$. We will integrate this dynamics using the velocity
Verlet scheme, which is explicit and of order 2, and reads
\begin{equation}
\label{eq:vv}
\begin{array}{rcl}
q_{n+1} &=& \dps
q_{n} + M^{-1} \left( \delta t \ p_n - \frac{\delta t^2}{2}
  \nabla V(q_n) \right), 
\\
p_{n+1} &=& \dps
p_n  - \frac{\delta t}{2} \left( \nabla V(q_n) + \nabla
  V(q_{n+1}) \right).
\end{array}
\end{equation}

\begin{figure}[htbp]
\begin{center}
\includegraphics[width = 6cm,angle =0]{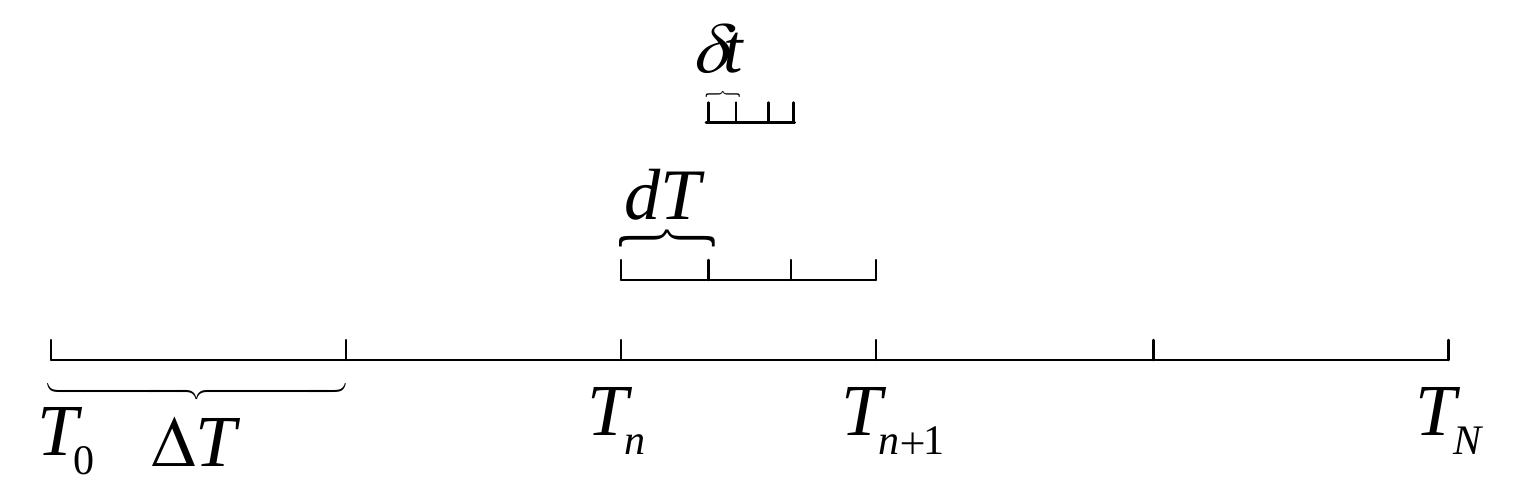}
\end{center}
\caption{Decomposition of the interval $[0,T]$ in sub-intervals
  $[T_n,T_{n+1}]$, on which the solution can be integrated using a
    coarse scheme of time step $d T$, or using a fine scheme of
    time step $\delta t$.}
\label{fig:parareel}
\end{figure}

The plain parareal algorithm  builds a sequence of $N$-tuples $\u^k \equiv
\left\{ u_n^k \right\}_{1 \leq n \leq N}$ that converges, when $k \rightarrow
\infty$, to the solution given by the fine scheme $\mathcal{F}$:
$\lim_{k \to \infty} u_n^k = u_n$.
In the sequel, we will consider Hamiltonian problems, for which
designing schemes with a non-uniform discretization is not
straightforward. We thus restrict ourselves to the case of a regular
discretization, namely $T_{n} = n \Delta T$ with $\Delta T = T/N$ for
some $N \in \N^\star$, where $[0,T]$ is the time interval of
interest. We introduce another approximation $\mathcal{G}$ of the exact
flow 
$\mathcal{E}$, which is not as accurate as $\mathcal{F}$, but is
much less expensive to use than $\mathcal{F}$. For example, we can choose the same
discretization scheme as for $\mathcal{F}$, but with a larger time
step $d T \gg \delta t$ (see Figure~\ref{fig:parareel}). Hence,
computing $\mathcal{G}_\tau$ amounts to performing $\tau/dT$ time steps of
length $dT$. Here again, the
dependency of $\mathcal{G}$ with respect to $dT$ is not explicitly written.
Another
possibility is to define the solver $\mathcal{G}$ from a simpler
problem, which does not contain as much information as the original
problem, and thus is easier to solve. We will use this opportunity in
Section~\ref{sec:outer}. In all the other sections, the only difference
between ${\cal G}$ and ${\cal F}$ lies in the choice of the time step.

Assume that we know an approximation 
$\left\{ u_n^k \right\}_{0 \leq n \leq N}$ of the solution
to~\eqref{eq:problem-1order}, at the end of some iteration $k$. Then the
parareal scheme defines the 
next iteration $\left\{ u_n^{k+1} \right\}_{0 \leq n \leq N}$ by 
\begin{equation}
\label{eq:scheme-parareal}
u_{n+1}^{k+1} = \mathcal{G}_{\Delta T}(u_{n}^{k+1}) +  
\mathcal{F}_{\Delta T}(u_{n}^{k})
- \mathcal{G}_{\Delta T}(u_{n}^{k}),
\end{equation}
with the initial condition $u_0^{k+1} = u_0$. 
At the beginning of iteration $k+1$, we first compute 
$\mathcal{F}_{\Delta T}(u_{n}^{k}) - \mathcal{G}_{\Delta
  T}(u_{n}^{k})$ in {\em parallel} over each interval $I_n$. Once this 
is completed, we only need to compute $\mathcal{G}_{\Delta T}(u_{n}^{k+1})$
and add to it the stored correction
$\mathcal{F}_{\Delta T}(u_{n}^{k}) - \mathcal{G}_{\Delta
T}(u_{n}^{k})$. This is a sequential process, the complexity of which
is negligible compared to the computation of $\mathcal{F}_{\Delta
T}(u_n^k)$. Note that an improved implementation in parallel has
recently been proposed in~\cite{batchelor} and explained hereafter. The
analysis of the complexities of the different parareal schemes proposed
in the sequel will use this implementation.

From the construction of the scheme, convergence in at most $N$
iterations can be demonstrated~\cite{maday-08}. 
It can also be proven that, on a fixed time interval $[0,T]$ and 
under some regularity conditions, the
scheme (\ref{eq:scheme-parareal}) yields a numerical 
solution $u_n^k$ after $k$ iterations that approximates $u(T_n)$ with an
error which is upper 
bounded (up to a multiplicative constant independent of the time steps
$\delta t$, $dT$ and $\Delta T$) by 
$\hbox{err } \mathcal{F} + [\hbox{err } \mathcal{G}]^k$, 
where $\hbox{err } \mathcal{F}$ is
the error on the solution at time $T$ for the fine solver and $\hbox{err }
\mathcal{G}$ is the error on the solution at time $T$ for the coarse
solver (see~\cite{lions-maday-turinici-01,bal-03,maday-08}).

\section{Parareal integration of Hamiltonian systems}
\label{sec:para4ham_gene}

As announced in the introduction, the purpose of this article is to design parallel-in-time numerical
schemes, derived from the parareal in time
algorithm~\eqref{eq:scheme-parareal}, for the integration of
Hamiltonian dynamical systems. We first review the specificities of such
dynamics before discussing their numerical integration with the parareal
algorithm~\eqref{eq:scheme-parareal}. 

\subsection{Numerical integration of Hamiltonian systems}
\label{sec:num_integ}

In this article, we consider finite dimensional Hamiltonian systems,
namely dynamical systems that read 
\begin{equation}
\label{eq:hamilton-eq}
\dot{q} = \frac{\partial H}{\partial p},
\quad
\dot{p} = -\frac{\partial H}{\partial q}, 
\end{equation}
where the Hamiltonian $H(q,p)$ is a smooth scalar function 
depending on the variables $q=(q_1,\ldots,q_d) \in \R^d$ and
$p=(p_1,\ldots,p_d) \in \R^d$. 

The evolution of many physical systems can be written as a
Hamiltonian dynamics. Examples include systems in molecular
dynamics (where $q$ and $p$ respectively represent the positions and momenta of the
atoms composing the system, see~\cite{frenkel}), celestial
mechanics (where 
$q$ and $p$ represent the positions and momenta of planets or
satellites~\cite{celestial1,celestial2,celestial3}), solid mechanics
(after space 
discretization, the wave equation modelling elastodynamics yields a
Hamiltonian dynamics of the type~\eqref{eq:hamilton-eq},
see~\cite{joly1,joly2}). In all 
these cases, $H(q,p)$ is, physically, the total energy of the system. 

Hamiltonian dynamics have very specific properties that we now briefly review
(see~\cite{hlw,ss-calvo} for more comprehensive expositions). 
First, the energy $H(q,p)$ is preserved
along the trajectory of~\eqref{eq:hamilton-eq}. Second, the flow of a
Hamiltonian system is {\em symplectic}. 
It is well-known that symplectic schemes, such as the velocity Verlet
scheme~\eqref{eq:vv}, are appropriate schemes for long time numerical
integration of Hamiltonian dynamics. They indeed define a 
symplectic numerical flow, and thus preserve, at the discrete level, a
fundamental geometrical property of the exact flow. In addition, the
numerical flow given by a symplectic scheme applied on 
the dynamics~\eqref{eq:hamilton-eq} can be shown to be almost equal to
the exact flow of a Hamiltonian dynamics, for a so-called {\em modified
  Hamiltonian}, which is close to the original Hamiltonian $H$
of~\eqref{eq:hamilton-eq}. This is one of the main results of the
celebrated backward error analysis for Hamiltonian dynamics
(see~\cite{benettin,hairer97,reich99}, and the comprehensive survey
in~\cite[Chap. IX]{hlw}). As a consequence, preservation of the energy
can be shown, in the sense 
that, for all $n$ such that $n \delta t \leq C \exp (c/\delta t)$,
$$
\left| H(q_n,p_n) - H(q_0,p_0) \right| \leq C (\delta t)^r,
$$
where $\delta t$ denotes the time step, $r$ the order of the scheme, and the constants $c$
and $C$ are independent of $\delta t$. The numerical error on the energy
hence does not grow with the simulation interval (at least for time
intervals exponentially long in the inverse of the time step). In some
cases (namely, the case of completely integrable systems and of almost
completely integrable systems, such as the solar system), the difference
between the numerical and the exact 
trajectories can be shown to grow linearly in time, rather than
exponentially, as would be the case for a generic integrator used on a generic dynamics $\dot{u}= f(u)$.

\medskip

In this article, we consider Hamiltonian systems for which the
energy reads
\begin{equation}
\label{eq:separable}
H(q,p) = \frac{1}{2} p^T M^{-1} p + V(q),
\end{equation}
where, as above, $M$ is a diagonal matrix, and $V$ is a smooth scalar function depending
on the positions $q$. For such an energy, the
dynamics~\eqref{eq:hamilton-eq} reads
\begin{equation}
\label{eq:hamilton-eq-sep}
\dot{q} = M^{-1} p, 
\quad
\dot{p} = -\nabla V(q). 
\end{equation}
Setting $u = (q,p)$, we recast this system as $\dot{u} = f(u)$. We then
observe that the system is reversible with 
respect to $\rho$ defined by $\rho(q,p) = (q,-p)$, that is $f$ satisfies
$f \circ \rho = - \rho \circ f$. As a consequence, for any $t$,
the exact flow $\Psi = {\mathcal E}_t$ of~\eqref{eq:hamilton-eq-sep} is
{\em $\rho$-reversible}, that is
\begin{equation}
\label{eq:reversibility}
\rho \circ \Psi = \Psi^{-1} \circ \rho.
\end{equation}

Another class of interesting schemes for {\em reversible} Hamiltonian
systems such as~\eqref{eq:hamilton-eq-sep} is the family of $\rho$-reversible
schemes: a one-step scheme defined by $u_{n+1} =
\Psi_{\delta t}(u_n)$ is $\rho$-reversible if $\Psi_{\delta t}$
satisfies~\eqref{eq:reversibility}. 
Properties similar to the ones obtained 
with a symplectic scheme (good preservation of the energy, \ldots) can be
shown when a reversible Hamiltonian dynamics is integrated using a
$\rho$-reversible scheme. Conditions under which such properties are
proven are however more restrictive than those of the backward error
analysis for symplectic schemes: the Hamiltonian system should be {\em reversible integrable} (this
implies that, if $q$ and $p$ are in $\R^d$, then the dynamics preserves
at least $d$ invariants) and a {\em non resonant} condition should be
satisfied~\cite[Chap. XI]{hlw}. 

The construction of $\rho$-reversible schemes is often done by first
designing {\em symmetric} schemes. A one-step scheme defined by $u_{n+1} =
\Psi_{\delta t}(u_n)$ is {\em symmetric} (see~\cite[Def. V.1.4]{hlw})
if, for any $\delta t$,
\begin{equation}
\label{eq:symmetry}
\Psi_{\delta t} \circ \Psi_{-\delta t} = \mbox{Id}.
\end{equation}
If a symmetric scheme satisfies the additional property
\begin{equation}
\label{eq:relation}
\rho \circ \Psi_{\delta t} = \Psi_{-\delta t} \circ \rho,
\end{equation}
then the scheme is $\rho$-reversible in the sense
of~\eqref{eq:reversibility} (see~\cite[Theo. V.1.5]{hlw}). In practice,
condition~\eqref{eq:relation} is much less restrictive than the symmetry
condition~\eqref{eq:symmetry}, and is satisfied by many schemes
(including e.g. the forward Euler scheme). In the sequel, we will design
symmetric schemes, satisfying~\eqref{eq:symmetry}, and check {\em a
  posteriori} that they indeed satisfy~\eqref{eq:relation}. These
schemes are hence $\rho$-reversible, namely satisfy~\eqref{eq:reversibility}.  

Despite the restrictions on the type of Hamiltonian systems for which
$\rho$-reversible schemes allow for good long-time properties, symmetric
schemes represent a very interesting alternative to symplectic schemes,
because they are often easier to design.

We conclude this section by recalling that the velocity Verlet
scheme~\eqref{eq:vv}, which is symplectic, as pointed out above, is also
symmetric and $\rho$-reversible.

\subsection{Parareal integration of Hamiltonian dynamics}
\label{sec:para4ham}

We observe that,
even if the coarse and fine propagators employed in the definition of
the parareal algorithm  (\ref{eq:scheme-parareal}) are symplectic or 
respectively symmetric, then, for a given parareal iteration $k \geq 1$ and
at a given time step $n$, with $n > k$, the application $u_0 =(q_0,p_0)
\mapsto u_n^k = (q_n^k,p_n^k)$ is neither symplectic nor respectively
symmetric. 

This lack of structure has immediate practical consequences when
employing the plain parareal algorithm~\eqref{eq:scheme-parareal} on
Hamiltonian systems. Consider as a first example the one-dimensional
harmonic oscillator, with Hamiltonian
\begin{equation}
\label{eq:H-ho}
H(q,p) = \frac{1}{2} p^2 + \frac{1}{2} q^2, 
\quad p \in \R, \ q \in \R,
\end{equation}
that we integrate up to time $T=10^{4}$. Set $\Delta T=0.2$, and 
consider the velocity Verlet scheme (which, we have recalled, is both
symplectic and 
symmetric) for the fine and the coarse propagators, with time
steps $\delta t = 10^{-3}$ and $dT = 0.1$, respectively. 

\begin{figure}[htbp]
\begin{center}
\includegraphics[width=9cm,angle=0,origin=0]{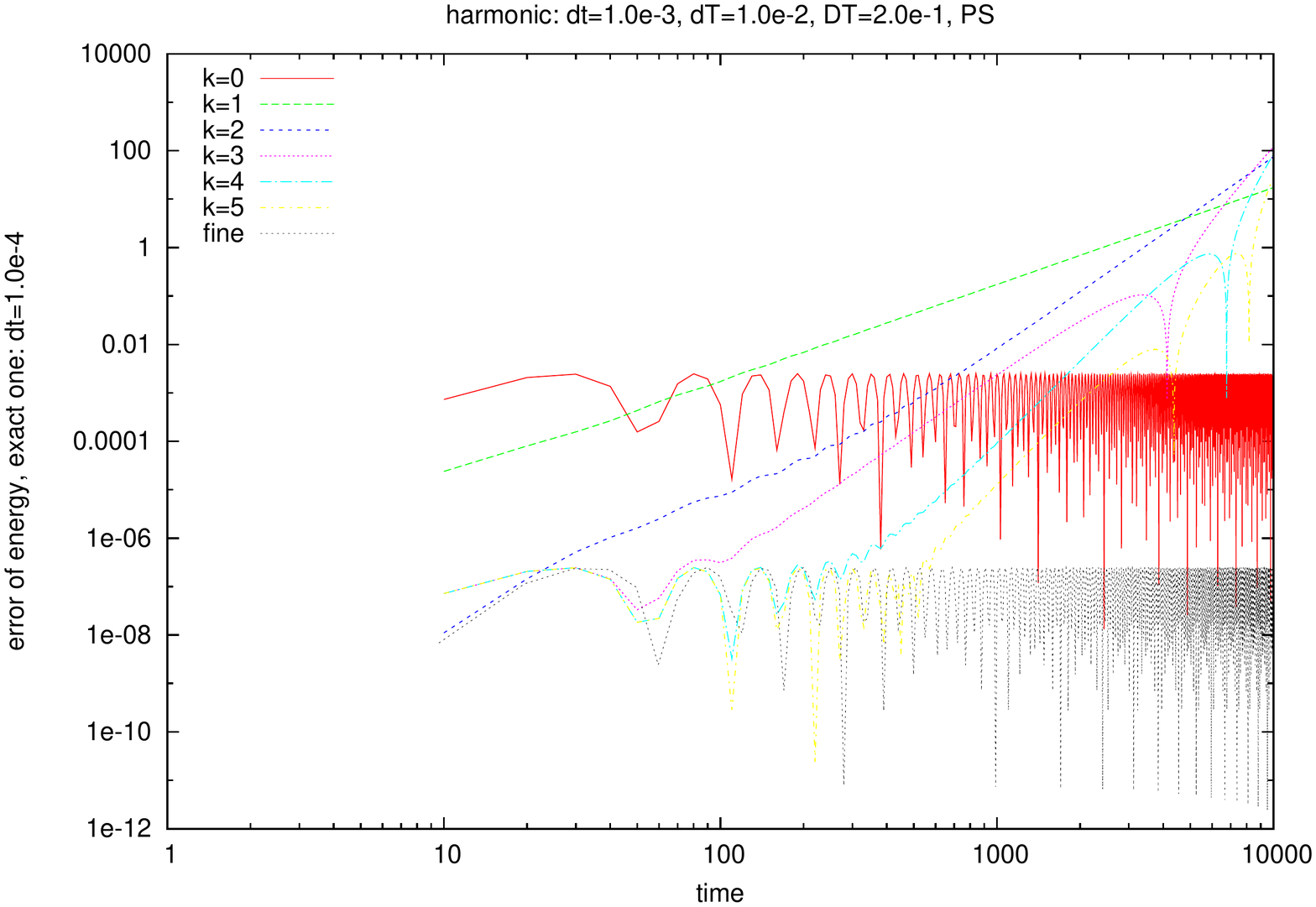}
\includegraphics[width=9cm,angle=0,origin=0]{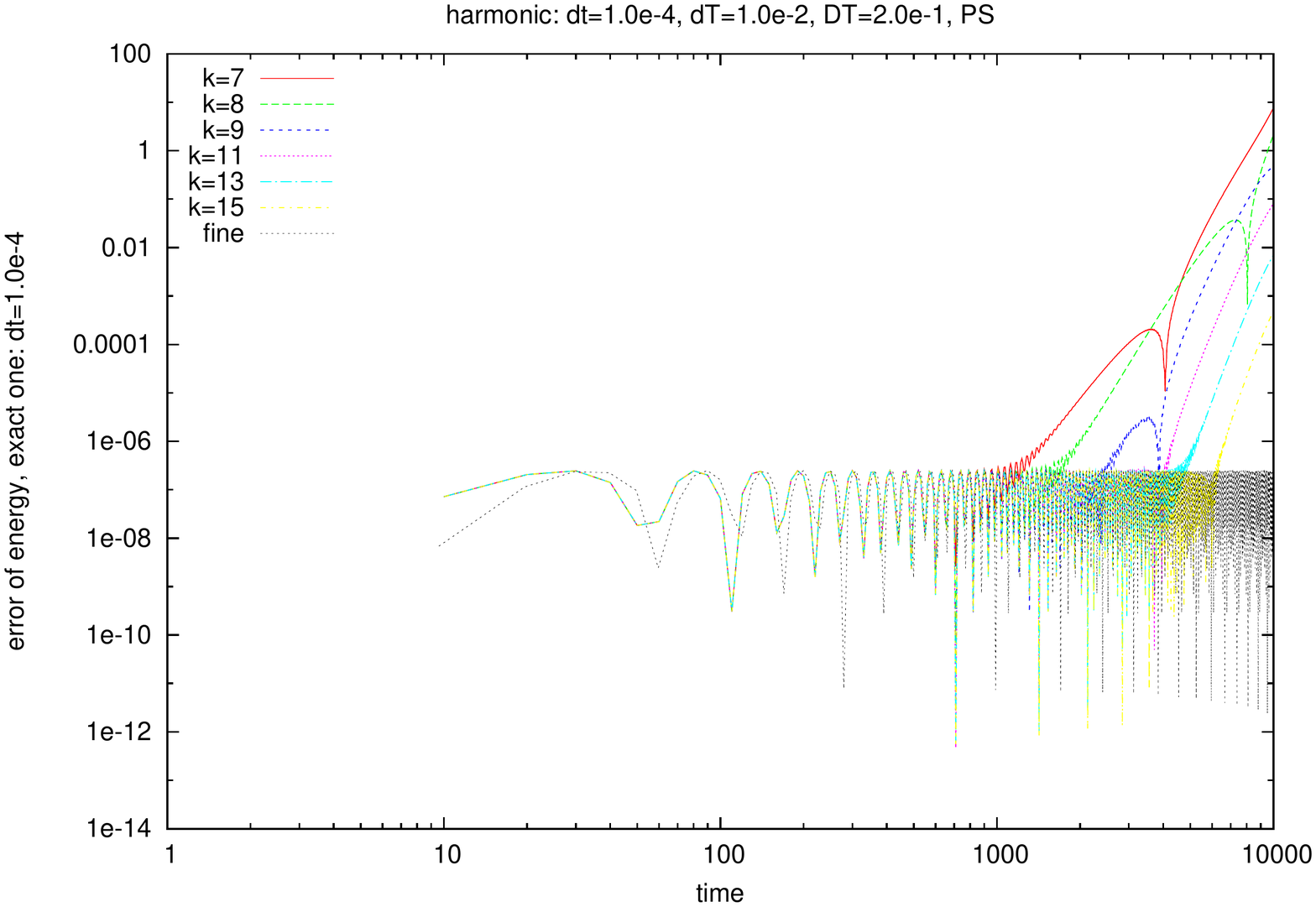}
\caption{\label{fig:harm-ener}
Error~\eqref{eq:err_H} on the energy for the harmonic oscillator problem,
obtained by the parareal method~\eqref{eq:scheme-parareal} with
$\delta t = 10^{-3}$, $dT = 0.1$, $\Delta T = 0.2$.} 
\end{center}
\end{figure}

As a confirmation of the lack of geometric structure, we observe the
lack of energy preservation. On Fig.~\ref{fig:harm-ener}, we plot the
relative error on the energy preservation, defined, at time $n \Delta T$
and iteration $k$, as 
\begin{equation}
\label{eq:err_H}
{\rm err}_n^k = \frac{\left| H(q_n^k,p_n^k) - H(q_0,p_0) \right|}
{\left| H(q_0,p_0) \right|}.
\end{equation}
We note that, for short times
(say up to $T=10^3$ for the iteration $k=5$), energy preservation is
equally good for the parareal algorithm and for the sequential fine
scheme. However, it deteriorates for larger times. Not unexpectedly, we
do not observe the traditional behavior of geometric schemes (either
symplectic or symmetric), namely a good preservation of energy, even
when the numerical trajectory is very different from the exact
trajectory.

\begin{figure}[htbp]
\begin{center}
\includegraphics[width=9cm,angle=0,origin=0]{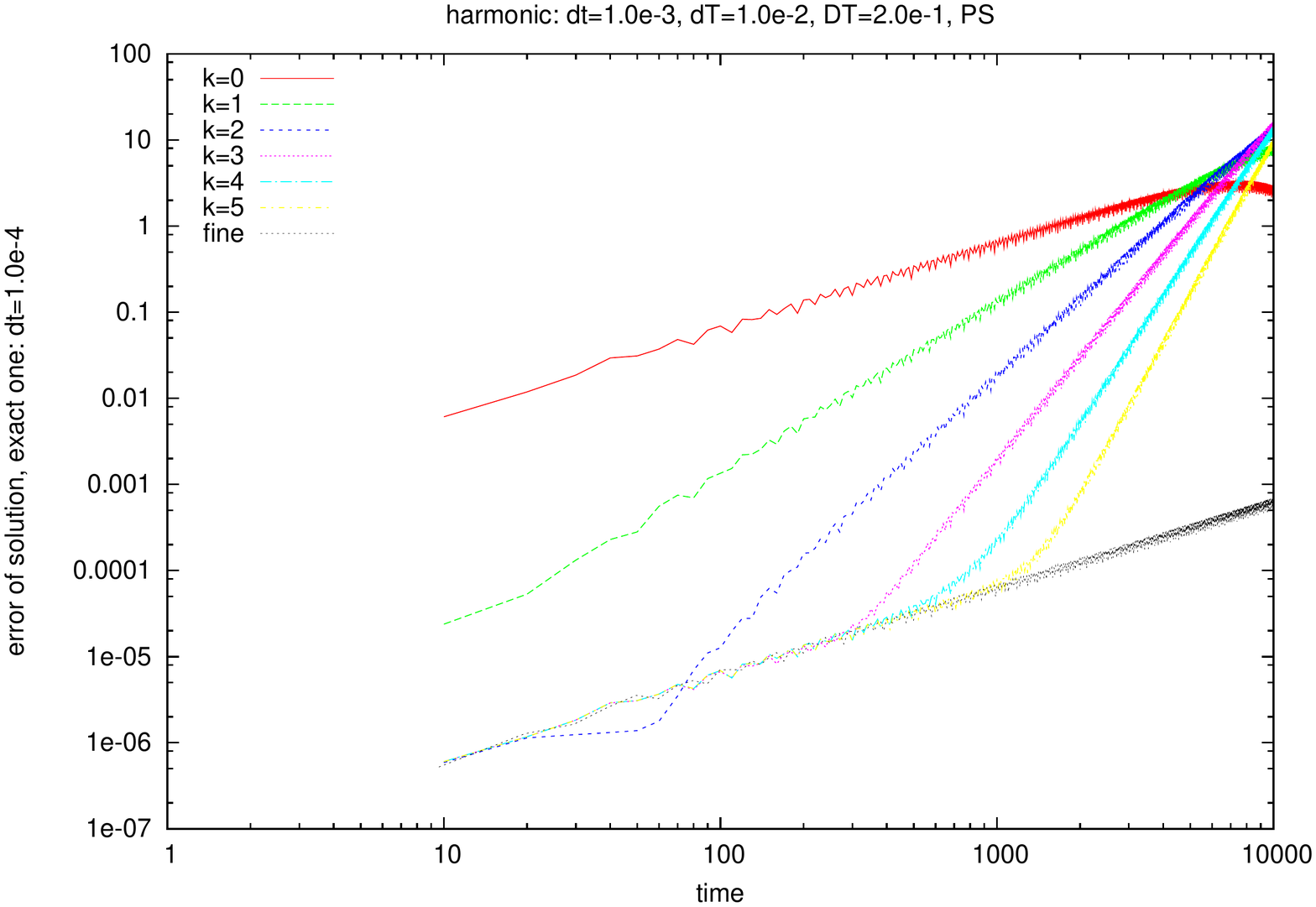}
\includegraphics[width=9cm,angle=0,origin=0]{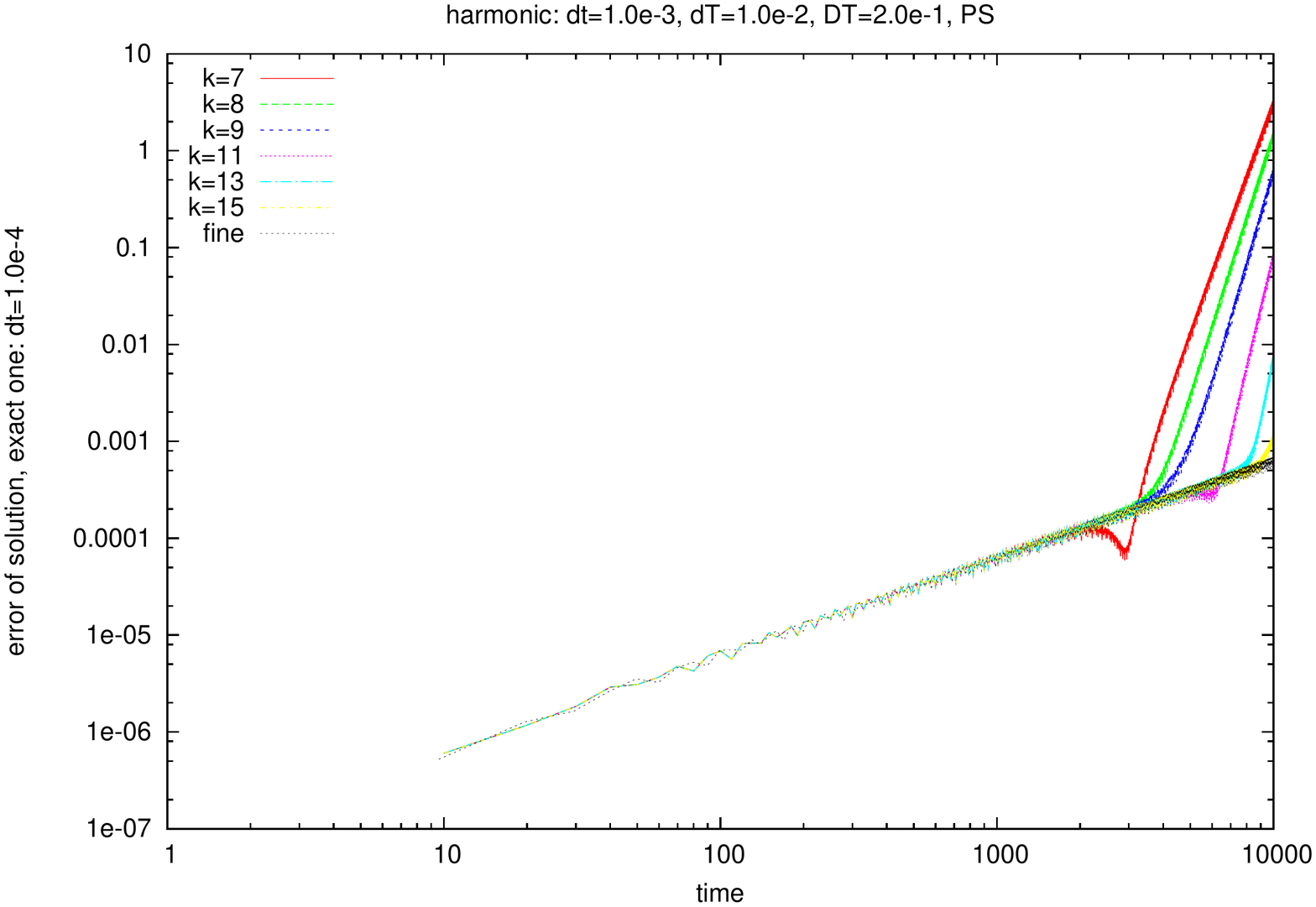}
\caption{\label{fig:harm-pos}
Error~\eqref{eq:err_traj} on the trajectory for the harmonic oscillator
problem, obtained by the parareal method~\eqref{eq:scheme-parareal} with
$\delta t = 10^{-3}$, $dT = 0.1$, $\Delta T = 0.2$.}
\end{center}
\end{figure}

On Fig.~\ref{fig:harm-pos}, we plot, for several iterations $k$,  the error 
\begin{equation}
\label{eq:err_traj}
{\rm err}_n^k = \| q_n^k - q_n^{\rm ref} \| + \| p_n^k - p_n^{\rm ref} \|
\end{equation}
on the trajectory $(q,p)$, with
respect to a reference trajectory $(q_{\rm ref},p_{\rm ref})$ computed with the
velocity Verlet algorithm --- used sequentially --- with the time step $\delta
t/10$, where $\delta t$ is the time step of the fine propagator
${\mathcal F}$. 
For $k=5$, on the time interval $[0,10^3]$, we note that results are very good: 
the parareal trajectory is very close to the fine scheme
trajectory, for a much smaller computational effort. Indeed, at $k=5$,
assuming that the complexity of the coarse propagations is negligible, each
processor has only computed 5 fine scale trajectories 
of length $\Delta T$, in contrast to the situation when a sequential algorithm
is used, and the processor has to compute $T/\Delta T = 5000$ fine scale
trajectories of length $\Delta T$ to reach the time $T=10^3$. However,
we also see that the error at iteration $k=5$ increases for $t \geq
10^3$, and becomes 
unacceptable for times $t$ of the order of $10^4$. With more iterations (say $k=15$),
convergence of the trajectories up to the time $T=10^4$ is obtained.  

Similar observations hold for the two-dimensional Kepler problem, where 
\begin{equation}
\label{eq:H-kepler}
H(q,p) = \frac{1}{2} p^T p - \frac{1}{\|q\|}, 
\quad p \in \R^2, \ q \in \R^2.
\end{equation}
We do not include them here for the sake of brevity.

\subsection{Evaluation of the complexity of the plain parareal
  algorithm} 
\label{sec:para4ham_complexity} 

In this section, for future reference, we assess the computational
complexity of the parareal algorithm~\eqref{eq:scheme-parareal} for the
integration of the Hamiltonian dynamics~\eqref{eq:hamilton-eq-sep}.  
We perform this evaluation under the
assumptions that (i) the coarse and the fine propagators integrate the same
dynamics, using an algorithm that requires the same fixed number of
calls (set here to one) to the right hand side
of~\eqref{eq:hamilton-eq-sep} per time step (we consider that the
complexities of evaluating $\nabla H$ or $\nabla V$ are equal, and we
denote it by $C_{\nabla}$), and
that (ii) we have all the necessary processors to perform all the fine
scale computations in parallel.

We use here the implementation quoted to us by~\cite{batchelor}. It
consists in starting the computation of $\mathcal{F}_{\Delta
  T}(u_{n}^{k+1})$ immediately after $u_{n}^{k+1}$ has been computed
in~(\ref{eq:scheme-parareal}), and not waiting that all
$(u_{n}^{k+1})_{1 \leq n \leq N}$ are available. This allows to start
the parareal iteration $k+2$ much sooner. The complexity of the first
coarse propagation scales as $C_{\nabla} \ T/dT$. We next
distinguish two extreme cases, whether the complexity of the fine propagator on
a time interval of size $\Delta T$ is smaller (respectively larger) than
the complexity of the coarse solver on the whole interval $[0,T]$. The
first case corresponds to the situation when 
$\dps \frac{\Delta T}{\delta t} <\frac{T}{dT}$, the second case when 
$\dps \frac{\Delta T}{\delta t} >\frac{T}{dT}$.

In the first case, the complexity of each iteration is dominated by the
complexity of the coarse solver on the whole interval $[0,T]$. The
complexity of such a parareal iteration is thus again of the order of
$C_{\nabla} \ T/dT$. In the second case, the complexity of each
parareal iteration is of the order of 
$C_{\nabla} \ \Delta T/\delta t$. 
In both cases, a final coarse iteration needs to be performed.

Denoting by $K_{\rm P}$ the number of parareal iterations, the
approximate complexity of the scheme~(\ref{eq:scheme-parareal}) is of
the order of 
$$
C_{\rm P} =  (K_{\rm P} +1) \, C_{\nabla} \, \frac{T}{dT}
$$ 
in the first case, and
$$
C_{\rm P} =  K_{\rm P} \, C_{\nabla} \, \frac{\Delta T}{\delta t}
$$ 
in the second one. 

In comparison, the complexity of the fully sequential algorithm, using
the small time step $\delta t$, is
$$
C_{\rm seq} = C_{\nabla} \, \frac{T}{\delta t}.
$$ 

Note that the second case above corresponds more to the paradigm of the
parareal scheme (remind also that, if possible, the coarse solver should
be based on a less expensive dynamics than the fine solver, as in
Section~\ref{sec:outer} below), especially for the integration of large
Hamiltonian systems as the one considered at the end of the article. In
that case, the speed-up is of the order $\dps \frac{T}{K_{\rm P} \Delta
  T}$, which is the number of processors divided by the number of
iterations.  

In the sequel, for the complexity analysis, we shall assume that we are
in the second case above, i.e. we assume that 
\begin{equation}
\label{assumption}
\frac{\Delta T}{\delta t} >\frac{T}{dT}.
\end{equation}

\section{A symmetric variant of the parareal in time algorithm}
\label{sec:sym}

In~\cite{bal}, a {\em symplectic} variant 
of the parareal algorithm has been introduced. The strategy there is
based on the reconstruction of the generating function $S : (q,p) \in
\R^{d} \times \R^d \mapsto \R$ associated with the
symplectic map 
$\mathcal{F}_{\Delta T} \circ \mathcal{G}^{-1}_{\Delta T}$. 
An interpolation procedure is used to obtain an approximation of that
generating function. The optimal way to perform the interpolation is
still an open question, especially when working with high dimensional
problems ($d \gg 1$). In this article, we limit ourselves to the 
design of a {\em symmetric} scheme, using the fact that it is often 
simpler to symmetrize a given scheme than to make it symplectic. We hope
to return to symplectic variants in future publications.  

\subsection{Derivation of the scheme}
\label{sec:design_parasym}
 
Our idea is based on the following well-known observation. Consider a general one-step
scheme $U_{n+1} = \Psi_{\Delta T}(U_n)$. Then the scheme
\begin{equation}
\label{eq:sym-gene}
U_{n+1} = \Psi_{\Delta T/2} \circ \left( \Psi_{-\Delta T/2} \right)^{-1}(U_{n})
\end{equation}
is symmetric (see \cite[Chap. V]{hlw}). For future use, we introduce the
intermediate variables $U_{n+1/2} = \left( \Psi_{-\Delta T/2}
\right)^{-1}(U_{n})$, and write the above scheme as
\begin{equation}
\label{eq:sym-gene_bis}
U_{n} = \Psi_{-\Delta T/2}(U_{n+1/2}),
\quad
U_{n+1} = \Psi_{\Delta T/2}(U_{n+1/2}).
\end{equation}

Let us now write the parareal algorithm~(\ref{eq:scheme-parareal}) in a form
more appropriate for our specific purpose. We first choose a number $K$ of iterations for the
parareal algorithm, and set
$$
U_n := (u_n^0, u_n^1, \cdots, u_n^K).
$$
Then the parareal scheme~\eqref{eq:scheme-parareal} can be written as
$$
U_{n+1} = \Psi_{\Delta T}(U_n),
$$
where the map $\Psi_{\Delta T}$ is defined by
$$
\Psi_{\Delta T}(U_n) = \left(
\begin{array}{l}
\mathcal{G}_{\Delta T}(u_n^0)
\\
\mathcal{G}_{\Delta T}(u_n^1) + \mathcal{F}_{\Delta T}(u_n^0)
- \mathcal{G}_{\Delta T}(u_n^0) 
\\ \cdots \\
\mathcal{G}_{\Delta T}(u_n^K) + \mathcal{F}_{\Delta T}(u_n^{K-1})
  - \mathcal{G}_{\Delta T}(u_n^{K-1})
\end{array}
\right).
$$

We now apply the symmetrization procedure~\eqref{eq:sym-gene_bis}
to the map $\Psi_{\Delta T}$, and consider the scheme 
\begin{equation}
\label{eq:sym-para}
\begin{array}{ccc}
U_{n} &=& \Psi_{-\Delta T/2}(U_{n+1/2}),
\\
U_{n+1} &=& \Psi_{\Delta T/2}(U_{n+1/2}).
\end{array}
\end{equation}
Let us write this scheme in a more detailed manner. The first line 
of~\eqref{eq:sym-para} reads
\begin{eqnarray*}
u_n^0 &=& \mathcal{G}_{-\Delta T / 2}(u_{n+1/2}^0), 
\\
u_n^1 &=& \mathcal{G}_{-\Delta T / 2}(u_{n+1/2}^1) +
\mathcal{F}_{-\Delta T / 2}(u_{n+1/2}^0) - \mathcal{G}_{-\Delta T / 2}(u_{n+1/2}^0),
\\
\ldots & = & \ldots 
\\
u_n^{K-1} &=& \mathcal{G}_{-\Delta T / 2}(u_{n+1/2}^{K-1}) +
 \mathcal{F}_{-\Delta T / 2}(u_{n+1/2}^{K-2}) -
 \mathcal{G}_{-\Delta T / 2}(u_{n+1/2}^{K-2}),
\\
u_n^{K} &=& \mathcal{G}_{-\Delta T / 2}(u_{n+1/2}^{K}) +
 \mathcal{F}_{-\Delta T / 2}(u_{n+1/2}^{K-1}) -
 \mathcal{G}_{-\Delta T / 2}(u_{n+1/2}^{K-1}).
\end{eqnarray*}
We hence obtain
\begin{eqnarray*}
u_{n+1/2}^0  &=& \mathcal{G}_{-\Delta T / 2}^{-1}(u_n^0), 
\\
u_{n+1/2}^1  &=& \mathcal{G}_{-\Delta T / 2}^{-1} \left[u_n^1 - 
\mathcal{F}_{-\Delta T / 2}(u_{n+1/2}^0) + 
\mathcal{G}_{-\Delta T / 2}(u_{n+1/2}^0) \right],
\\
\ldots &=& \ldots
\\
u_{n+1/2}^{K-1} &=& \mathcal{G}_{-\Delta T / 2}^{-1} \left[u_n^{K-1} -
\mathcal{F}_{-\Delta T / 2}(u_{n+1/2}^{K-2}) +
\mathcal{G}_{-\Delta T / 2}(u_{n+1/2}^{K-2}) \right],
\\
u_{n+1/2}^{K}  &=& \mathcal{G}_{-\Delta T / 2}^{-1} \left[u_n^{K} -
\mathcal{F}_{-\Delta T / 2}(u_{n+1/2}^{K-1}) +
\mathcal{G}_{-\Delta T / 2}(u_{n+1/2}^{K-1}) \right].
\end{eqnarray*}
We now collect the above set of relations with the second line
of~\eqref{eq:sym-para} and obtain the following formulation:  
\begin{equation}
\label{eq:sym-para-0}
u_{n+1/2}^{0} = \mathcal{G}_{-\Delta T/2}^{-1}(u_n^0), \quad
u_{n+1}^0 = \mathcal{G}_{\Delta T/2}(u_{n+1/2}^0),
\end{equation}
at iteration $k=0$, and, for $k \geq 1$, 
\begin{equation}
\label{eq:sym-para-k}
\left\{
\begin{array}{rcl} 
u_{n+1/2}^{k+1} &=&
\mathcal{G}_{-\Delta T/2}^{-1}\left[u_{n}^{k+1} -
\mathcal{F}_{-\Delta T/2}(u_{n+1/2}^{k}) + 
\mathcal{G}_{-\Delta T/2}(u_{n+1/2}^{k})\right], 
\\ \noalign{\vskip 4pt}
u_{n+1}^{k+1} &=& \mathcal{G}_{\Delta T/2}(u_{n+1/2}^{k+1})
+ \mathcal{F}_{\Delta T/2}(u_{n+1/2}^{k}) - 
\mathcal{G}_{\Delta T/2}(u_{n+1/2}^{k}).
\end{array}
\right.
\end{equation}
In the sequel, the scheme (\ref{eq:sym-para-0})-(\ref{eq:sym-para-k}) is
called the {\em symmetric parareal scheme}.

\begin{remark}
We have recalled in Section~\ref{sec:num_integ} that $\rho$-reversible
schemes have good geometrical properties. By construction, the scheme 
(\ref{eq:sym-para-0})-(\ref{eq:sym-para-k}) is symmetric. 
We have checked that it also satisfies condition~\eqref{eq:relation}. As
a consequence (see again Section~\ref{sec:num_integ}), the scheme
(\ref{eq:sym-para-0})-(\ref{eq:sym-para-k}) is $\rho$-reversible. 
\end{remark}

We note that, apart from the computation of $\mathcal{G}_{-\Delta T /
  2}^{-1}$, the scheme (\ref{eq:sym-para-0})-(\ref{eq:sym-para-k}) is
completely explicit, as soon as the propagators $\mathcal{G}$ and
$\mathcal{F}$  themselves are explicit. In particular, we point out that
the inverse of the fine propagator ${\mathcal F}$ is not needed. 

Let us now show that, in this new algorithm, all the expensive
computations (involving $\mathcal{F}$) 
can actually be performed in parallel. Assume that, at a given iteration $k$, 
we know $\left\{ u_n^{k} \right\}_{0 \leq n \leq N}$ and 
$\left\{ u_{n+1/2}^{k} \right\}_{0 \leq n \leq N-1}$. Then
the correction terms
$$
\mathcal{F}_{\Delta T/2}(u_{n+1/2}^{k}) -
\mathcal{G}_{\Delta T/2}(u_{n+1/2}^{k})
$$
and
$$
\mathcal{F}_{-\Delta T/2}(u_{n+1/2}^{k}) -
\mathcal{G}_{-\Delta T/2}(u_{n+1/2}^{k})
$$
can be computed independently, over each processor. The algorithm can
next proceed 
with sequential, but inexpensive, computations. Note again that, as
pointed out in Section~\ref{sec:para4ham_complexity}, we can start the
computation of $\mathcal{F}_{\Delta T/2}(u_{n+1/2}^{k})$ and 
$\mathcal{F}_{-\Delta T/2}(u_{n+1/2}^{k})$ as soon as $u_{n+1/2}^{k}$ is
available (we do not have to wait for all the computations of iteration
$k$ to be completed). 

The symmetric parareal algorithm is summarized as follows.
\begin{algorithm}
\label{algorithm-sps}
Let $u_0$ be the initial condition. 
\begin{enumerate}
\item Initialization: set $u_{0}^{0} = u_{0}$, and sequentially compute 
$\left\{ u_{n+1/2}^0 \right\}_{0 \leq n \leq N-1}$ and
$\left\{ u_{n+1}^0 \right\}_{0 \leq n \leq N-1}$ as
$$
u_{n+1/2}^{0} = \mathcal{G}_{-\Delta T/2}^{-1}(u_{n}^0), 
\quad
u_{n+1}^{0} = \mathcal{G}_{\Delta T/2}(u_{n+1/2}^0).
$$
\item Assume that, for some $k \geq 0$, the sequences 
$\left\{ u_{n}^k \right\}_{0 \leq n \leq N}$ and 
$\left\{ u_{n+1/2}^k \right\}_{0 \leq n \leq N-1}$ are known. Compute
the sequences $\left\{ u_{n}^{k+1} \right\}_{0 \leq n \leq N}$ and 
$\left\{ u_{n+1/2}^{k+1} \right\}_{0 \leq n \leq N-1}$ at iteration
$k+1$ by the following steps:
\begin{enumerate}
\item For all $0 \leq n \leq N-1$, compute in parallel 
$\mathcal{F}_{-\Delta T/2}(u_{n+1/2}^{k})$, 
$\mathcal{F}_{\Delta T/2}(u_{n+1/2}^{k})$,
$\mathcal{G}_{-\Delta T/2}(u_{n+1/2}^{k})$ and 
$\mathcal{G}_{\Delta T/2}(u_{n+1/2}^{k})$;
\item Set $u_{0}^{k+1} = u_{0}$;
\item Compute, for $0 \le n \le N-1$,
$$
\hspace{-15mm}
\left\{
\hspace{-1mm}
\begin{array}{rcl} 
u_{n+1/2}^{k+1} &=&
\mathcal{G}_{-\Delta T/2}^{-1} \left[ u_{n}^{k+1} -
\mathcal{F}_{-\Delta T/2}(u_{n+1/2}^{k}) + 
\mathcal{G}_{-\Delta T/2}(u_{n+1/2}^{k})\right], 
\\ \noalign{\vskip 5pt}
u_{n+1}^{k+1} &=& \mathcal{G}_{\Delta T/2}(u_{n+1/2}^{k+1})
+ \mathcal{F}_{\Delta T/2}(u_{n+1/2}^{k}) - 
\mathcal{G}_{\Delta T/2}(u_{n+1/2}^{k}).
\end{array}
\right.
$$
\end{enumerate}
\end{enumerate}
\end{algorithm}

Some comments are in order. The parareal
scheme~\eqref{eq:scheme-parareal} is not a one-step scheme 
defining 
$u_{n+1}^{k+1}$ from $u_n^{k+1}$, and hence the symmetric form
(\ref{eq:sym-para-0})-(\ref{eq:sym-para-k}) cannot be considered,
strictly speaking, as a symmetric integrator in the classical
sense. However, several reasons lead us to believe that this
algorithm is the appropriate generalization of symmetric integrators
when dealing with parareal-type algorithms.

First, the scheme at iteration $k=0$, defined by~(\ref{eq:sym-para-0}),
is exactly the symmetrization~\eqref{eq:sym-gene} of the coarse
propagator $\mathcal{G}_{\Delta T}$. 

Second, the flow (\ref{eq:sym-para-0})-(\ref{eq:sym-para-k}) is
symmetric in the sense that, if for some $n$ and some~$k$,
$$
(u_{n+1}^{0}, \cdots, u_{n+1}^{k}, u_{n+1}^{k+1})
$$
is obtained from 
$$
(u_{n}^{0}, \cdots, u_{n}^{k}, u_{n}^{k+1})
$$
by the flow implicitly defined in
(\ref{eq:sym-para-0})-(\ref{eq:sym-para-k}), then 
$$
(u_{n}^{0}, \cdots, u_{n}^{k}, u_{n}^{k+1})
$$  
is obtained from
$$
(u_{n+1}^{0}, \cdots, u_{n+1}^{k}, u_{n+1}^{k+1})
$$ 
by the exact same algorithm reversing the time. In fact, only the
consideration (and the storage) of the last two iterates (in terms of
``parareal iterates'') $(u_n^{k}, u_{n}^{k+1})$ is required to perform
the iterations, and the above argument.

Third, if the coarse propagator happens to be identical to the fine
one (which is of course not supposed to be the case for efficiency of
the parareal integrator!), then the symmetrized form
(\ref{eq:sym-para-k}) reads 
$$
u_{n+1/2}^{k+1} = \mathcal{G}_{-\Delta T/2}^{-1}(u_{n}^{k+1}),
\quad
u_{n+1}^{k+1} = \mathcal{G}_{\Delta T/2}(u_{n+1/2}^{k+1}),
$$
and it thus coincides with the standard symmetrized version of the
coarse propagator.

Finally, formally taking the limit $k \rightarrow +\infty$ in
(\ref{eq:sym-para-k}) yields 
$$
u_{n+1}^{\infty} = \mathcal{F}_{\Delta T/2} \circ 
(\mathcal{F}_{-\Delta T/2})^{-1} (u_n^{\infty}).
$$
This shows that the limit of the symmetric parareal algorithm in
terms of parareal iterations coincides with a standard symmetrized
form of the map $\mathcal{F}_{\Delta T}$. Note also that
this latter algorithm is {\em not} the symmetrized version of the fine
propagator, since symmetrization does not occur after each time step, but
after $\Delta T/(2\delta t)$ time steps. 

These observations show the formal consistency of our notion
of symmetrization for parareal-type integrators with the
classical notion of symmetry.

\begin{remark}
Instead of considering~\eqref{eq:sym-gene} to symmetrize a
given scheme $\Psi_{\Delta T}$, one can alternatively consider the
symmetric scheme
$$
U_{n+1} = \left( \Psi_{-\Delta T/2} \right)^{-1} \circ 
\Psi_{\Delta T/2}(U_{n}).
$$
In the parareal context, this yields an algorithm with similar
properties as the above
algorithm~(\ref{eq:sym-para-0})-(\ref{eq:sym-para-k}), 
and with comparable numerical results.   In
  the sequel, for the sake of brevity, we only consider the
  symmetrization procedure~\eqref{eq:sym-gene}.
\end{remark}

\subsection{Evaluation of the complexity of Algorithm~\ref{algorithm-sps}}
\label{sec:para5ham_complexity} 

We assess the computational complexity of Algorithm~\ref{algorithm-sps}
under the same 
assumptions as in Section~\ref{sec:para4ham_complexity}. We again
assume~(\ref{assumption}), and we again start to compute
$\mathcal{F}_{\Delta T/2}(u_{n+1/2}^{k})$ and 
$\mathcal{F}_{-\Delta T/2}(u_{n+1/2}^{k})$ as soon as possible, e.g. as
soon as $u_{n+1/2}^{k}$ has been computed. 

The complexity of the first coarse propagation scales as $C_{\nabla} \
T/dT$. The complexity of each parareal iteration is of the order
of $C_{\nabla} \ \Delta T/\delta t$ (note that we need to call
twice the fine solver, but only for time intervals of length $\Delta 
T/2$). A final coarse iteration needs to be done, the complexity of which is
neglected in comparison to the cost of a fine propagation. 
Denoting by $K_{\rm SP}$ the number of parareal iterations, we
obtain that the complexity of this scheme is
$$
C_{\rm SP} =  K_{\rm SP} \, C_{\nabla} \, \frac{\Delta T}{\delta t}.
$$ 
The complexity is the same as that of the plain parareal algorithm. The
speed-up, which is equal to $\dps \frac{T}{K_{\rm SP} \Delta T}$, is again
equal to the number of processors divided by the number of iterations.

\subsection{Numerical examples}
\label{sec:symm_num}

Let us now again consider the example of the harmonic
oscillator~\eqref{eq:H-ho}, and integrate this Hamiltonian system with
our newly constructed symmetric parareal
algorithm~(\ref{eq:sym-para-0})-(\ref{eq:sym-para-k}). Note that
system~\eqref{eq:H-ho} is an integrable reversible Hamiltonian system,
and thus belongs to a class of problems for which symmetric integrators
have been shown to preserve energy~\cite[Chap. XI]{hlw}, as has been
recalled in Section~\ref{sec:num_integ}.

\begin{figure}[htbp]
\begin{center}
\includegraphics[width=9cm,angle=0,origin=0]{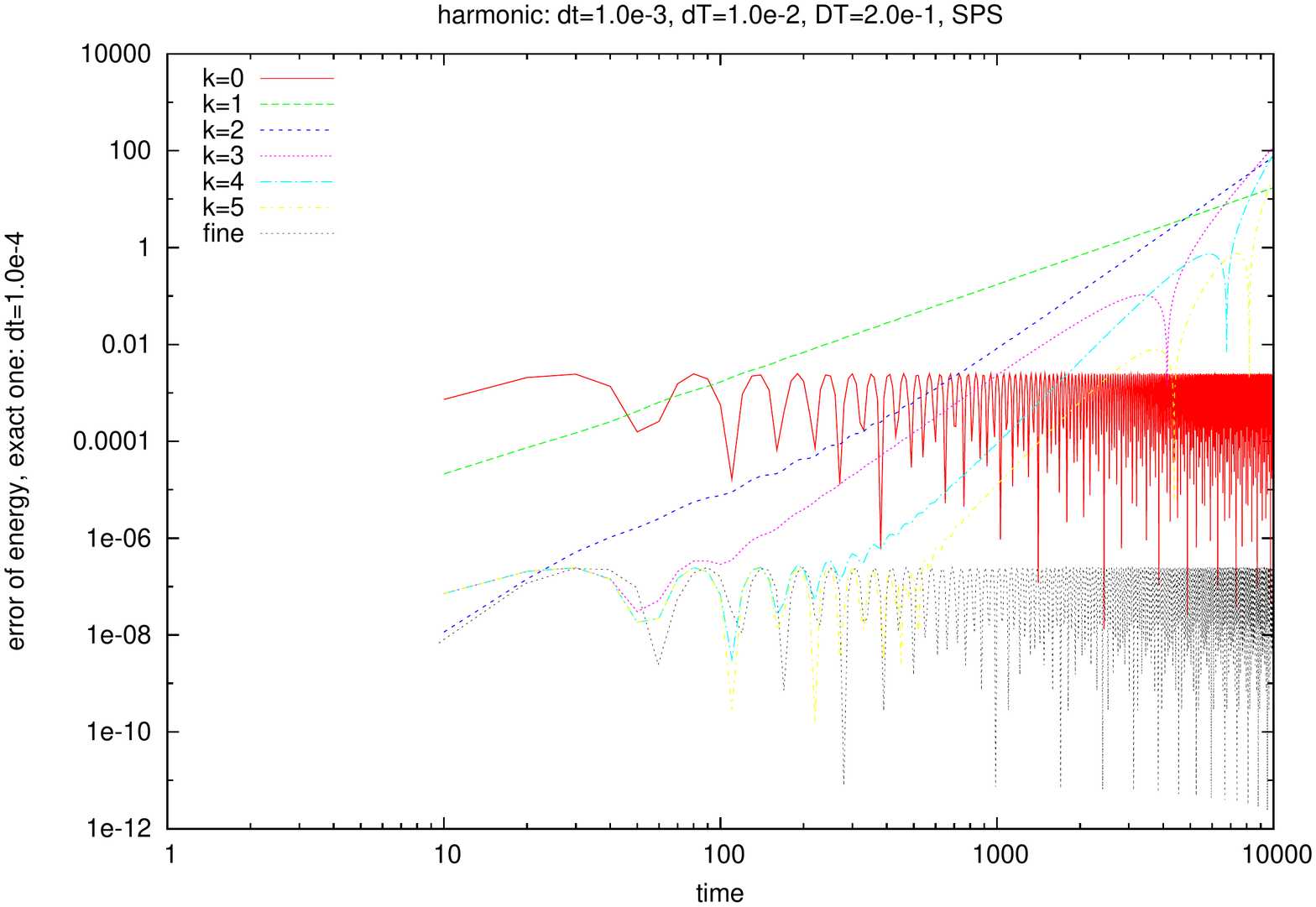}
\includegraphics[width = 9cm, angle=0, origin=0]{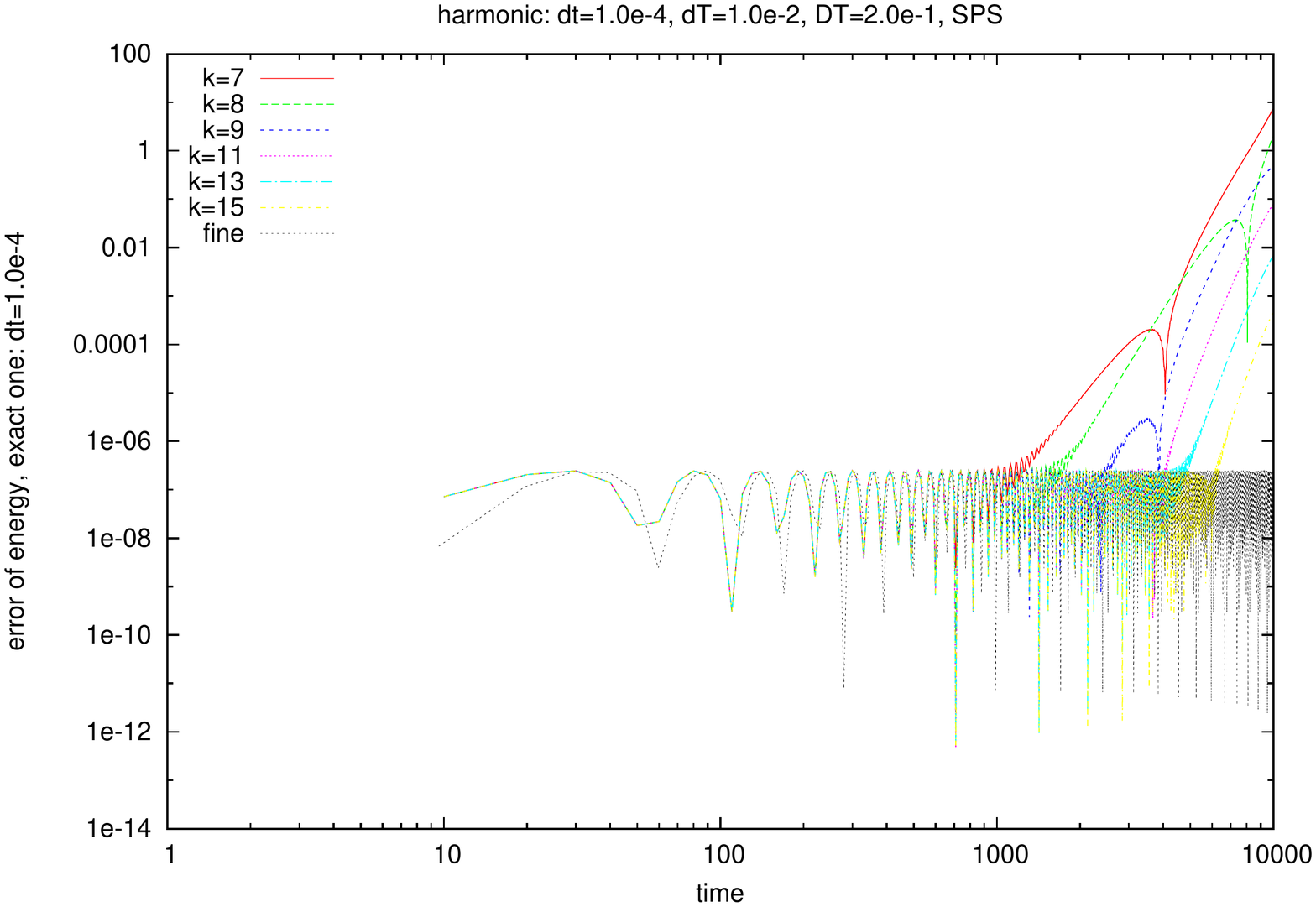}
\caption{ \label{fig:harm-sym-ener}
Errors on the energy for the harmonic oscillator problem, obtained by
Algorithm~\ref{algorithm-sps} with 
$\delta t = 10^{-3}$, $dT = 0.1$, $\Delta T = 0.2$.}
\end{center}
\end{figure}

\begin{figure}[htbp]
\begin{center}
\includegraphics[width=9cm,angle=0,origin=0]{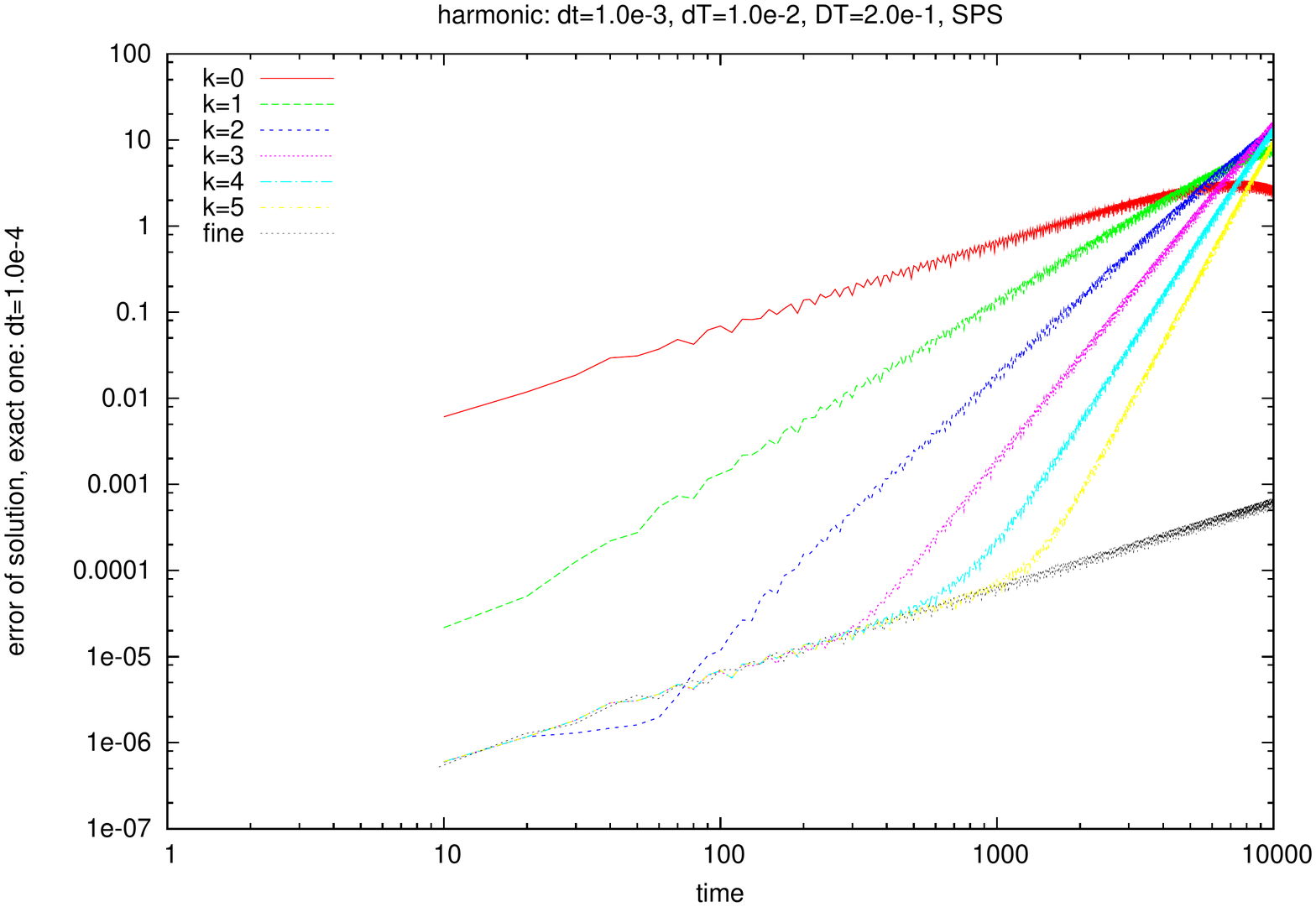}
\includegraphics[width=9cm,angle=0,origin=0]{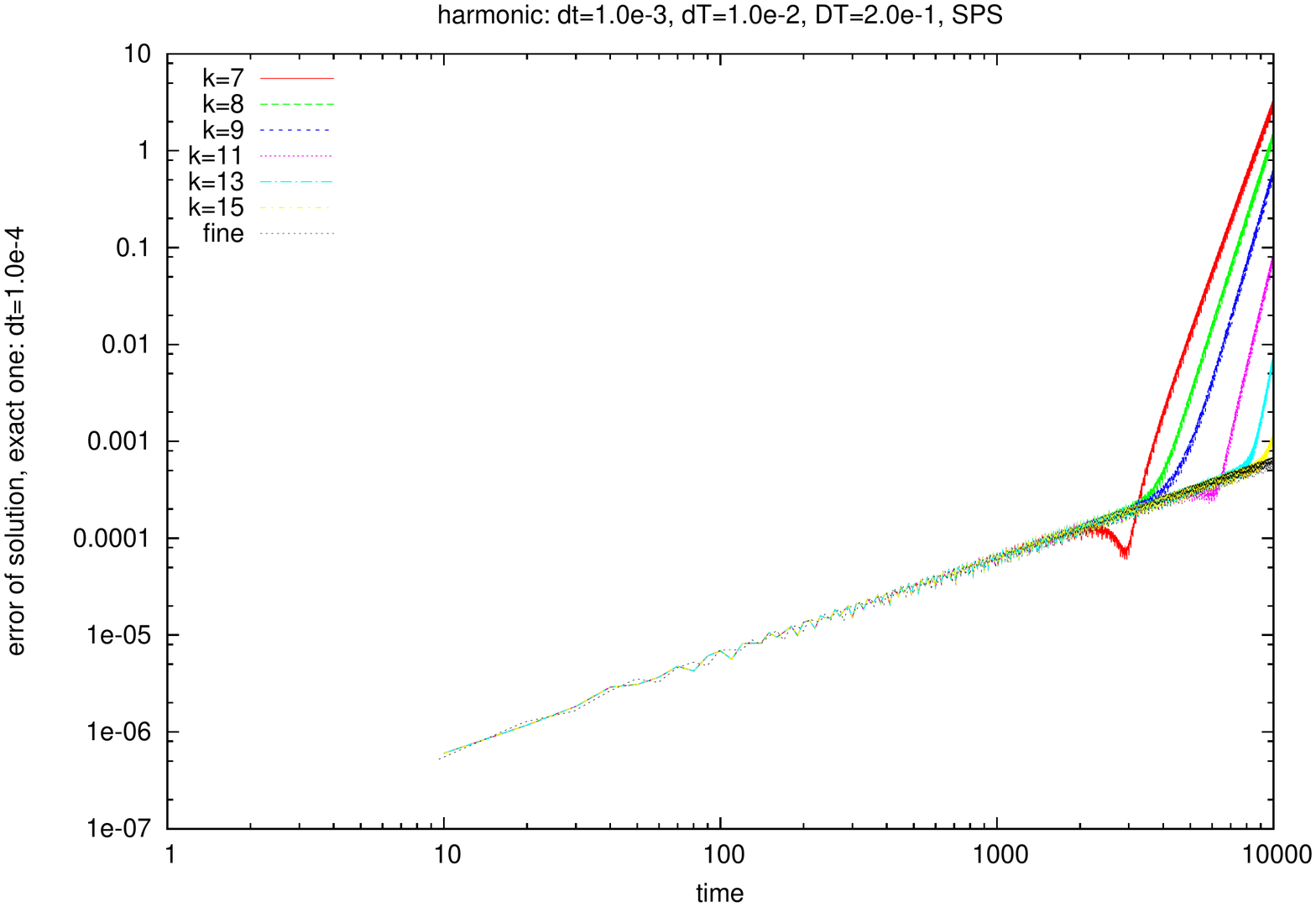}
\caption{ \label{fig:harm-sym-pos}
Errors on the trajectory for the harmonic oscillator problem, obtained
by Algorithm~\ref{algorithm-sps}
with $\delta t = 10^{-3}$, $dT = 0.1$, $\Delta T = 0.2$.}
\end{center}
\end{figure}

The relative errors~\eqref{eq:err_H} on the energy preservation are
shown on Fig.~\ref{fig:harm-sym-ener}, while the
errors~\eqref{eq:err_traj} on the trajectories are shown on
Fig.~\ref{fig:harm-sym-pos}. Parameters have been chosen  
such that the complexity to reach a given time step $n$ at a given iteration
$k$ is equal for the symmetric parareal
algorithm considered here and the standard parareal
algorithm~\eqref{eq:scheme-parareal} discussed 
in the previous sections. Results are disappointing: they are indeed
very similar, both qualitatively and quantitatively, to those
obtained with the plain parareal scheme~(\ref{eq:scheme-parareal}) (see
Section~\ref{sec:para4ham}, Figs.~\ref{fig:harm-ener}
and~\ref{fig:harm-pos}). 

Similar results have been obtained for the Kepler problem, which is
also an integrable reversible Hamiltonian system. 

\medskip

It is useful and instructive to now explain such a poor behavior, using
the following formal 
elements of analysis. We first observe that a parareal
integrator may be seen, at parareal iteration $k$, as an integrator of
a system consisting of $k+1$ identical replicas of the original system under
consideration. From (\ref{eq:sym-para-0}), we know that the first replica
is integrated by a symmetric algorithm. If the system under study is an
integrable reversible Hamiltonian system, its energy is thus preserved
in the long time by the simulation at parareal iteration $k=0$ (this is
confirmed by numerical experiments, see {\it
  e.g.} Fig.~\ref{fig:harm-sym-ener}, curve $k=0$). Next, since the
replicas are noninteracting, the system of replicas is evidently an
integrable reversible system, and has an energy that is equal to $k+1$
times the energy of 
the original system. Assume now that the symmetric propagator
(\ref{eq:sym-para-0})-(\ref{eq:sym-para-k}) applied to the system of
$k+1$ identical replicas 
conserves the energy of the global system (i.e. \emph{the sum} 
$\sum_{\ell=0}^k H(q^\ell,p^\ell)$ of the energies of the replicated
systems), as we could expect from a 
symmetric scheme applied to an integrable reversible Hamiltonian
system. Under this assumption, it follows, by induction on $k$, that the
energy is preserved along the trajectory, at {\em each} parareal
iteration $k$. This is clearly in contradiction with the observed
numerical results! The flaw in the above argument is that, for a
symmetric scheme to preserve the energy of an integrable reversible
system, we have recalled  that, among other conditions, the frequencies present
in the system have to satisfy a non-resonant diophantine
condition. This is precisely not the case here for the system of
replicas, by replication of the original frequencies!

The theoretical argument showing energy preservation does not apply, and numerical
results confirm that energy is indeed not well-preserved, in the
long-time limit. The symmetric parareal
scheme~(\ref{eq:sym-para-0})-(\ref{eq:sym-para-k}) hence needs to be
somehow amended.

\subsection{Symmetric parareal algorithm with frequency perturbation}

To prevent the different replicas from being resonant, a
possibility is to consider, at each parareal iteration $k$, a system
slightly different from the original system. For instance, in the case
of the harmonic oscillator, we may consider at each iteration $k$ a harmonic
oscillator with a specific  frequency $\omega_k$:
$$
H_k(q,p) = \frac{1}{2} p^2 + \frac{1}{2} \omega_k^2 \, q^2.
$$
The unperturbed case corresponds to $\omega_k = \omega_{\rm exact} = 1$
in the above energy.  
Provided the $\omega_k$ are all different from one another (and
non-resonant), the system of replicas is non-resonant. 
In the following,
the shift is chosen such that it vanishes when $k \to \infty$, i.e. $\lim_{k \to \infty} \omega_k = \omega_{\rm exact}$.

More precisely, we introduce a fine propagator
$\mathcal{F}^{(k)}_{\Delta T}$ and a coarse propagator
$\mathcal{G}^{(k)}_{\Delta T}$ for the Hamiltonian dynamics associated
to the Hamiltonian $H_k$. Rather than symmetrizing the parareal
algorithm~\eqref{eq:scheme-parareal}, we consider the scheme
\begin{equation}
\label{eq:scheme-parareal-shift}
\begin{array}{rcl}
u_{n+1}^0 &=& \mathcal{G}^{(0)}_{\Delta T}(u_{n}^0),
\\
u_{n+1}^{k+1} &=& \mathcal{G}^{(k+1)}_{\Delta T}(u_{n}^{k+1}) +  
\mathcal{F}^{(k+1)}_{\Delta T}(u_{n}^{k})
- \mathcal{G}^{(k+1)}_{\Delta T}(u_{n}^{k}),
\end{array}
\end{equation}
(note the upper indices $(k+1)$) and symmetrize this scheme along the lines of
Section~\ref{sec:design_parasym}. We thus obtain, at iteration $k=0$,
the algorithm
\begin{equation}
\label{eq:sym-para-0-shift}
u_{n+1/2}^{0} = \left( \mathcal{G}^{(0)}_{-\Delta T/2} \right)^{-1}(u_n^0), \quad
u_{n+1}^0 = \mathcal{G}^{(0)}_{\Delta T/2}(u_{n+1/2}^0),
\end{equation}
and, for $k \geq 1$, 
\begin{equation}
\label{eq:sym-para-k-shift}
\hspace{-1mm}
\left\{
\hspace{-2mm}
\begin{array}{rcl} 
u_{n+1/2}^{k+1} \! \! &=& \! \!
\left( \mathcal{G}^{(k+1)}_{-\Delta T/2} \right)^{-1}\left[u_{n}^{k+1} -
\mathcal{F}^{(k+1)}_{-\Delta T/2}(u_{n+1/2}^{k}) + 
\mathcal{G}^{(k+1)}_{-\Delta T/2}(u_{n+1/2}^{k})\right] 
\\
u_{n+1}^{k+1} \! \! &=& \! \! \mathcal{G}^{(k+1)}_{\Delta T/2}(u_{n+1/2}^{k+1})
+ \mathcal{F}^{(k+1)}_{\Delta T/2}(u_{n+1/2}^{k}) - 
\mathcal{G}^{(k+1)}_{\Delta T/2}(u_{n+1/2}^{k}).
\end{array}
\right.
\end{equation}

Let us briefly discuss the consistency of this scheme. First, at
iteration $k+1$, if the coarse propagator $\mathcal{G}^{(k+1)}$ is
identical to the fine propagator $\mathcal{F}^{(k+1)}$, then
(\ref{eq:sym-para-k-shift}) reads
$$
u_{n+1/2}^{k+1} = \left( \mathcal{F}^{(k+1)}_{-\Delta T/2} \right)^{-1}(u_{n}^{k+1}),
\quad
u_{n+1}^{k+1} = \mathcal{F}^{(k+1)}_{\Delta T/2}(u_{n+1/2}^{k+1}),
$$
which is a scheme consistent with the Hamiltonian dynamics driven by
$H_{k+1}$. Since the perturbation can be neglected when $k \to \infty$, we
recover, in the limit of large $k$, a scheme consistent with the
original dynamics. 

In addition, formally taking the limit $k \rightarrow +\infty$ in
(\ref{eq:sym-para-k-shift}) yields 
$$
u_{n+1}^{\infty} = \mathcal{F}_{\Delta T/2} \circ 
(\mathcal{F}_{-\Delta T/2})^{-1} (u_n^{\infty}),
$$
where $\mathcal{F}$ is the fine propagator associated to the original
dynamics. This shows that the limit of the
algorithm~\eqref{eq:sym-para-0-shift}-\eqref{eq:sym-para-k-shift} in
terms of parareal iterations coincides with a standard symmetrized
form of the fine propagator of the original dynamics. 

\begin{figure}[htbp]
\begin{center}
\includegraphics[width=9cm,angle=0,origin=0]{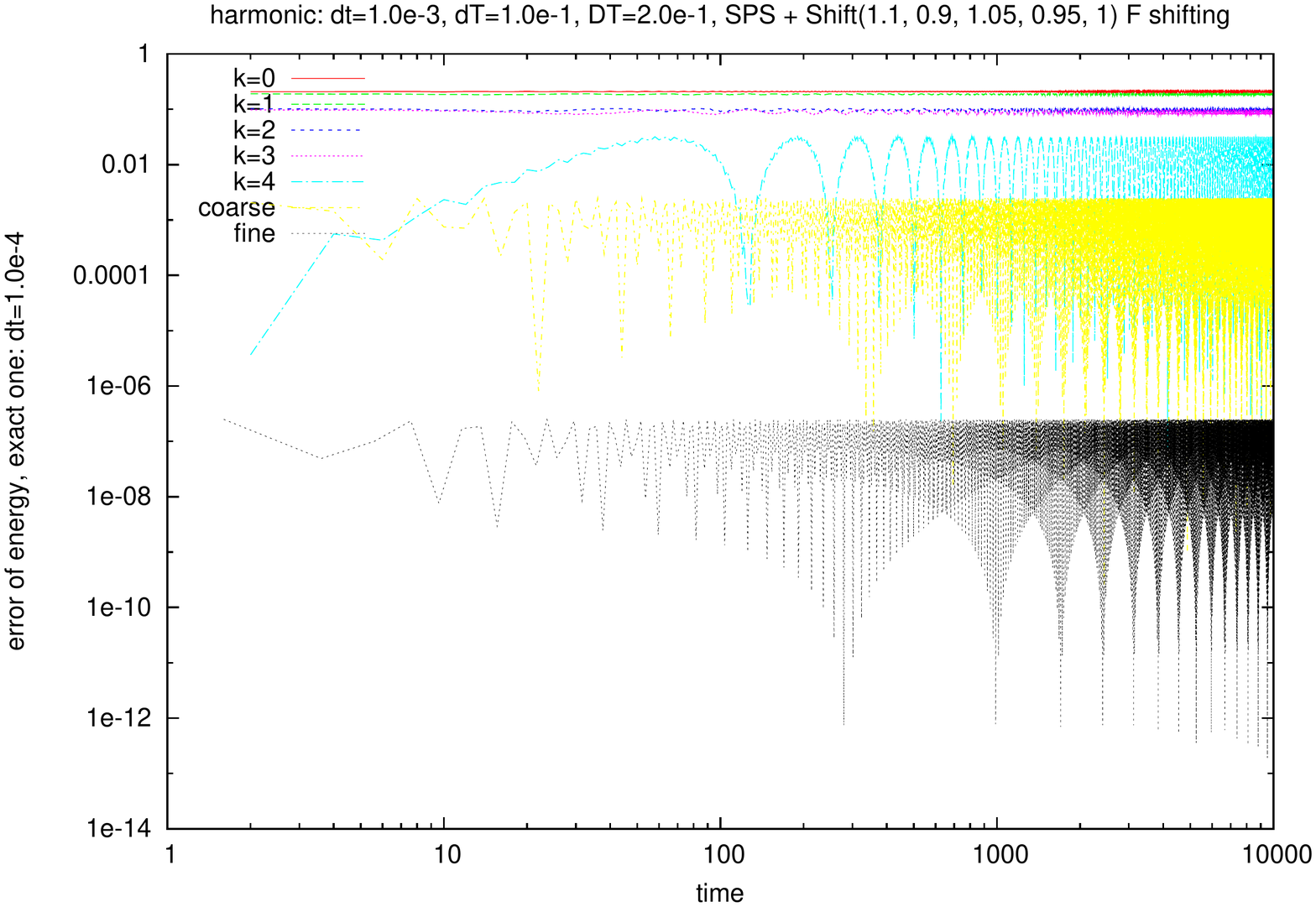}
\includegraphics[width=9cm,angle=0,origin=0]{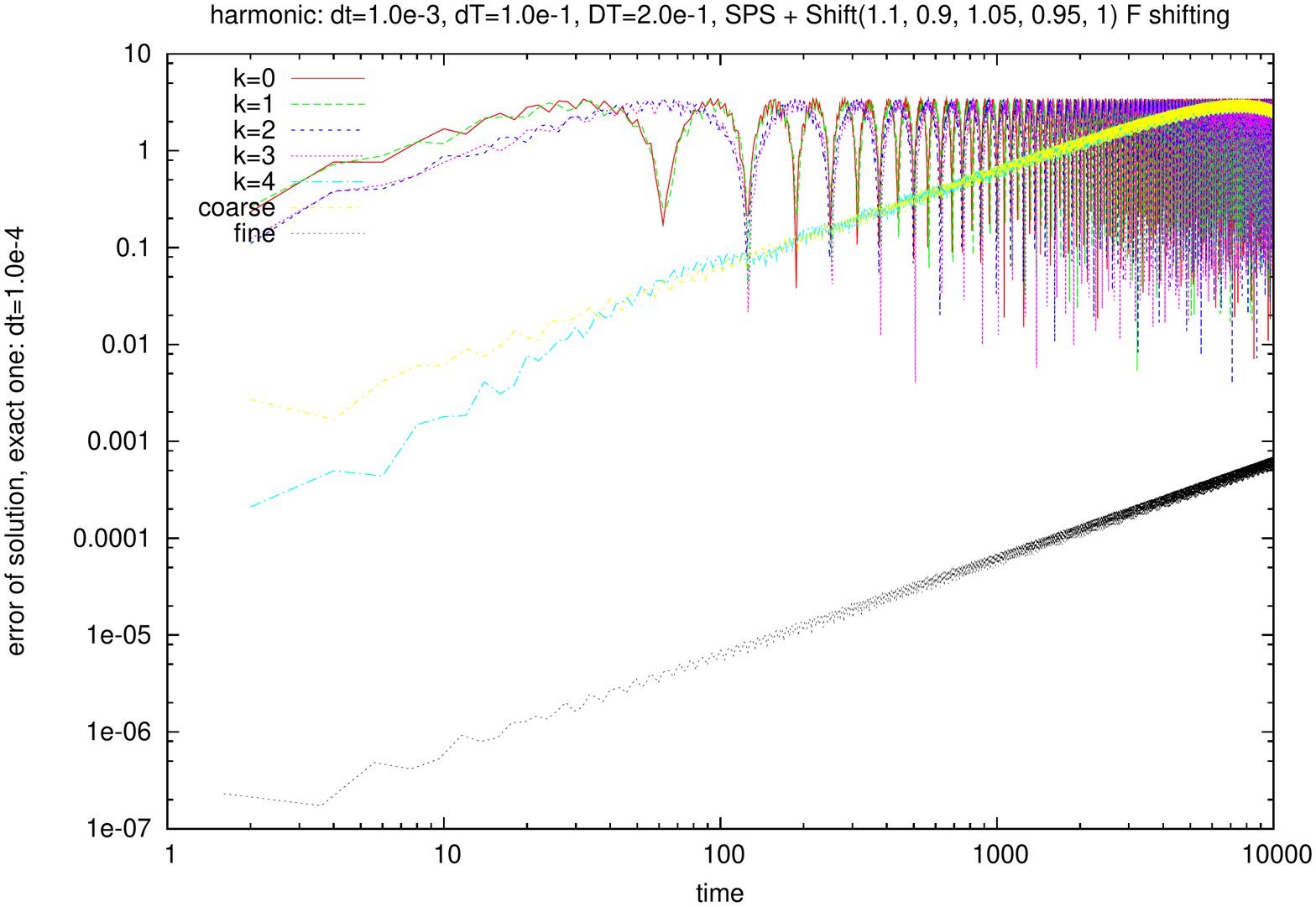}
\caption{ \label{fig:harm-shift}
Errors on the energy (left) and errors on the trajectory (right) for the
harmonic oscillator problem, obtained with the symmetric parareal
method~(\ref{eq:sym-para-0-shift})-(\ref{eq:sym-para-k-shift}), with frequency
perturbation. The frequencies are $\omega_{k=0}=1.1$,
$\omega_{k=1}=0.9$, 
$\omega_{k=2}=1.05$, $\omega_{k=3}=0.95$ and $\omega_{k=4}=\omega_{\rm exact}=1$ 
($\delta t = 10^{-3}$, $dT = 0.1$, $\Delta T = 0.2$).
} 
\end{center}
\end{figure}

On Fig.~\ref{fig:harm-shift}, we plot the errors~\eqref{eq:err_H} on the
energy and the errors~\eqref{eq:err_traj} on the trajectory, obtained
with the algorithm
(\ref{eq:sym-para-0-shift})-(\ref{eq:sym-para-k-shift}) applied to the 
harmonic oscillator problem. We focus on the results obtained at the
last parareal iteration, the only one for which the perturbation
vanishes. We observe that the energy does not drift,
while the trajectory is not extremely accurate (actually, both
observations hold also for intermediate values of $k$). We are hence in the
regime of geometric integration, with a good accuracy on energy being
{\em not} a consequence of a good accuracy on the trajectory. This is an
improvement in comparison to the previous algorithms 
(parareal algorithm~\eqref{eq:scheme-parareal} and symmetric parareal
algorithm~(\ref{eq:sym-para-0})-(\ref{eq:sym-para-k})). 

\begin{figure}[htbp]
\begin{center}
\includegraphics[width=9cm,angle=0,origin=0]{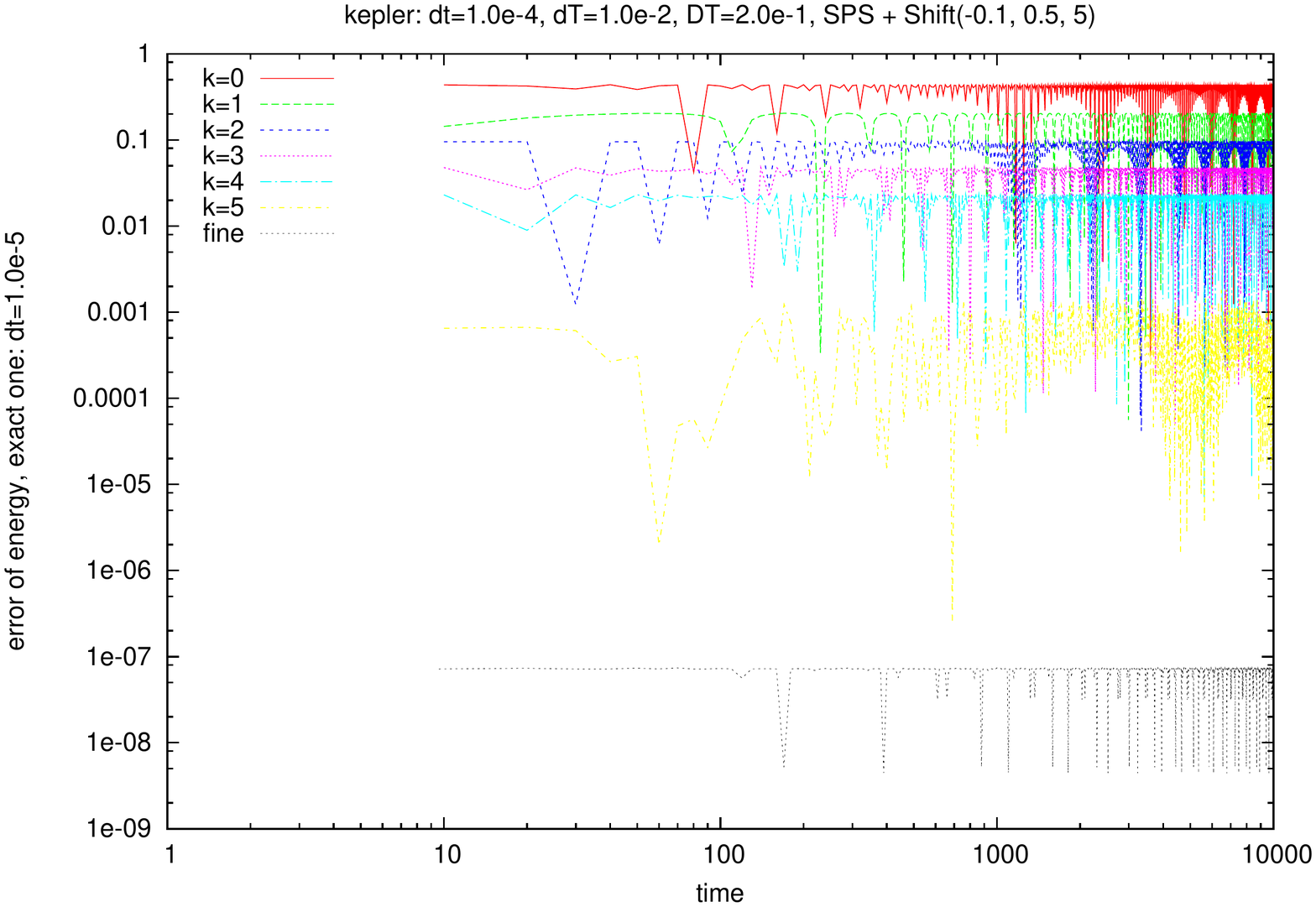}
\includegraphics[width=9cm,angle=0,origin=0]{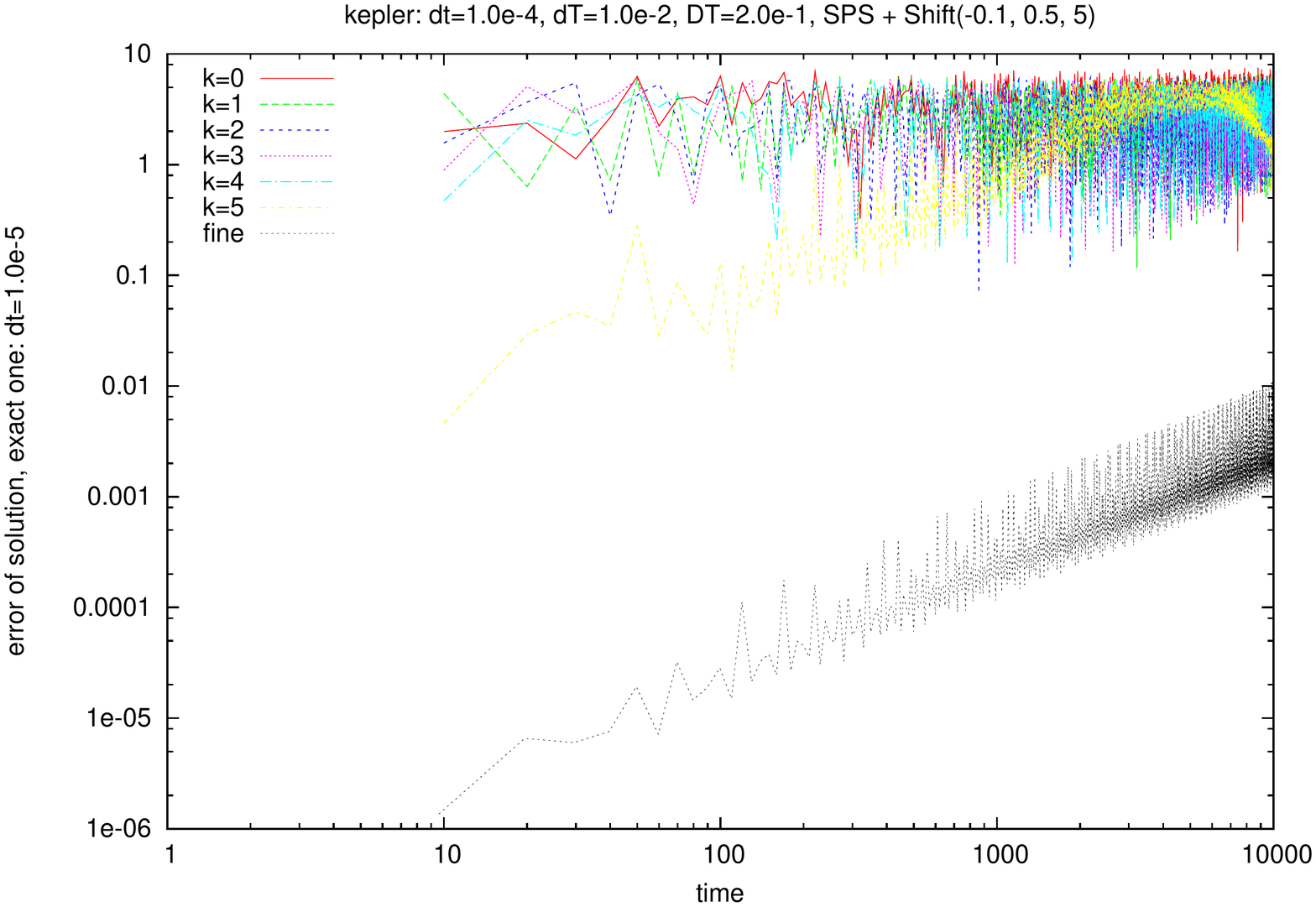}
\caption{ \label{fig:kepler-shift}
Error on the energy (left) and errors on the trajectory (right) for 
Kepler problem, obtained with the symmetric parareal
method~(\ref{eq:sym-para-0-shift})-(\ref{eq:sym-para-k-shift}), with frequency
perturbation. The coefficient $\alpha_k$ in~\eqref{eq:kepler-pert} is
$\alpha_{k=0}=0.9$,
$\alpha_{k=1}=0.95$, $\alpha_{k=2}=0.975$, $\alpha_{k=3}=0.9875$,
$\alpha_{k=4}=0.99375$ and $\alpha_{k=5}=\alpha_{\rm exact} = 1$ 
($\delta t = 10^{-4}$, $dT = 10^{-2}$, $\Delta T = 0.2$).
} 
\end{center}
\end{figure}

We next turn to the Kepler problem, where we introduce a perturbation by
considering 
\begin{equation}
\label{eq:kepler-pert}
H_k(q,p) = \frac{1}{2} p^T p - \frac{\alpha_k}{\| q \|},
\end{equation}
where the unperturbed case corresponds to $\alpha_k = \alpha_{\rm exact}
= 1$. Results are 
shown on Fig.~\ref{fig:kepler-shift}. The same conclusions as for the
harmonic oscillator test case hold. 

We hence numerically observe, in both cases considered, that the
energy does not drift when we use the
scheme~(\ref{eq:sym-para-0-shift})-(\ref{eq:sym-para-k-shift}), and that
energy preservation does not owe to trajectory accuracy. This formally
validates our understanding: the energy drifts observed with the 
symmetric parareal algorithm~(\ref{eq:sym-para-0})-(\ref{eq:sym-para-k})
are due to resonances. As soon as the replicas are made non-resonant, we
recover a behavior in terms of energy preservation that is typical of
geometric algorithms.  

Nevertheless, even though energy does not drift, we also see that it is
not extremely well preserved. It is often the case that, at the final
iteration, energy preservation is actually not really
better than if obtained using the coarse
propagator on the original dynamics. We also observe that the error on
the trajectory remains 
quite large, for the practical values of $k$ that we consider. So the
algorithm~(\ref{eq:sym-para-0-shift})-(\ref{eq:sym-para-k-shift}),
although better than~\eqref{eq:scheme-parareal}
and~(\ref{eq:sym-para-0})-(\ref{eq:sym-para-k}), is not 
satisfactory. 

\begin{remark}
Another possibility to integrate at each parareal iteration systems
that are non resonant from one another is to use, at each iteration $k$, a
different time step $dT_k$ for the coarse propagator. This strategy,
that is more general than the one discussed above (indeed, ``shifting the
frequency'' may not be easy if the system is not explicitly an harmonic
oscillator), yields results that are qualitatively similar to the ones
reported here.
\end{remark}

\begin{remark}
The idea of ``shifting the frequency'' is independent from the fact that
the algorithm is symmetric. It is hence possible to apply the {\em standard} parareal
method~\eqref{eq:scheme-parareal} on systems that are slightly different
from one another at each parareal iteration (in
(\ref{eq:sym-para-k-shift}), we have used the {\em symmetric} version of the
parareal algorithm on such shifted systems). Results are slightly better
when this strategy is used with the symmetric version of the algorithm.
\end{remark}

\section{Parareal algorithm with projection}
\label{sec:proj}

As an alternative to symmetrization described in Section~\ref{sec:sym},
we consider in this section a different idea. We couple the plain
parareal method with a projection of the trajectory on a specific
manifold, defined by the preservation of some invariants of the system
(namely the energy, and possibly other ones). 

Observe that many systems have actually more than one preserved
quantity: this is the case for the two-dimensional Kepler problem,
where the angular momentum $L = q \times p$ and the Runge-Lenz-Pauli vector 
$\dps{ A = p \times L - \frac{q}{\| q \|} }$ are also preserved,
or for any completely integrable system. However, the generic situation
is that some invariant quantities are known (among which the energy),
and that there may be other quantities that are preserved by the
dynamics, but may not have been identified. As a consequence, we first
consider a generic method, based on projecting on the constant
energy manifold 
$$
{\cal M} = \left\{ (q,p) \in \R^{2d}; \ H(q,p) = H_0 \right\},
$$
for the Kepler problem. Note that the harmonic oscillator
  is not a discriminating test case in this respect, since  the
energy is then the {\em unique} invariant. Next, for the sake of
completeness, again for the Kepler problem, we consider a method based
on projecting on the manifold where {\em both} $H$ and $L$ are kept constant
(see Section~\ref{sec:2inv}). 

\begin{remark}
Symplectic algorithms, such as SHAKE~\cite{shake} or
RATTLE~\cite{rattle1,rattle2},  
have been proposed for the numerical integration
of Hamiltonian systems restricted to some manifold defined by holonomic
constraints. Note however that here, the constraint at hand, $H(q,p) =
H_0$, is {\em not} holonomic. 
\end{remark}

Following~\cite[Chap. IV]{hlw}, we consider the following
definition of the projection $\pi_{\cal M}$ onto the manifold ${\cal
  M}$. For any 
$\widetilde{y} = (\widetilde{q},\widetilde{p}) \in \R^{2d}$, we define 
$(q,p) = y = \pi_{\cal M}(\widetilde{y}) \in \R^{2d}$
by 
\begin{equation}
\label{eq:proj-eq}
y = \pi_{\cal M}(\widetilde{y})
= \widetilde{y} + \lambda \nabla H(\widetilde{y}),
\quad \text{where $\lambda \in \R$ is such that} \quad
H(y) = H_0.
\end{equation}

\subsection{Parareal algorithm with projection} 

The coarse propagator being the velocity Verlet algorithm,
which yields an acceptable energy preservation, we do not need to
project on the constant energy manifold at iteration $k=0$. But, as
shown by our numerical results of the previous sections, the
trajectory departs from this manifold for the iterations $k \geq
1$, in the long time limit. This motivates the consideration of the
following algorithm.  

\begin{figure}[htbp]
\begin{center}
\includegraphics[width = 8cm,angle =0]{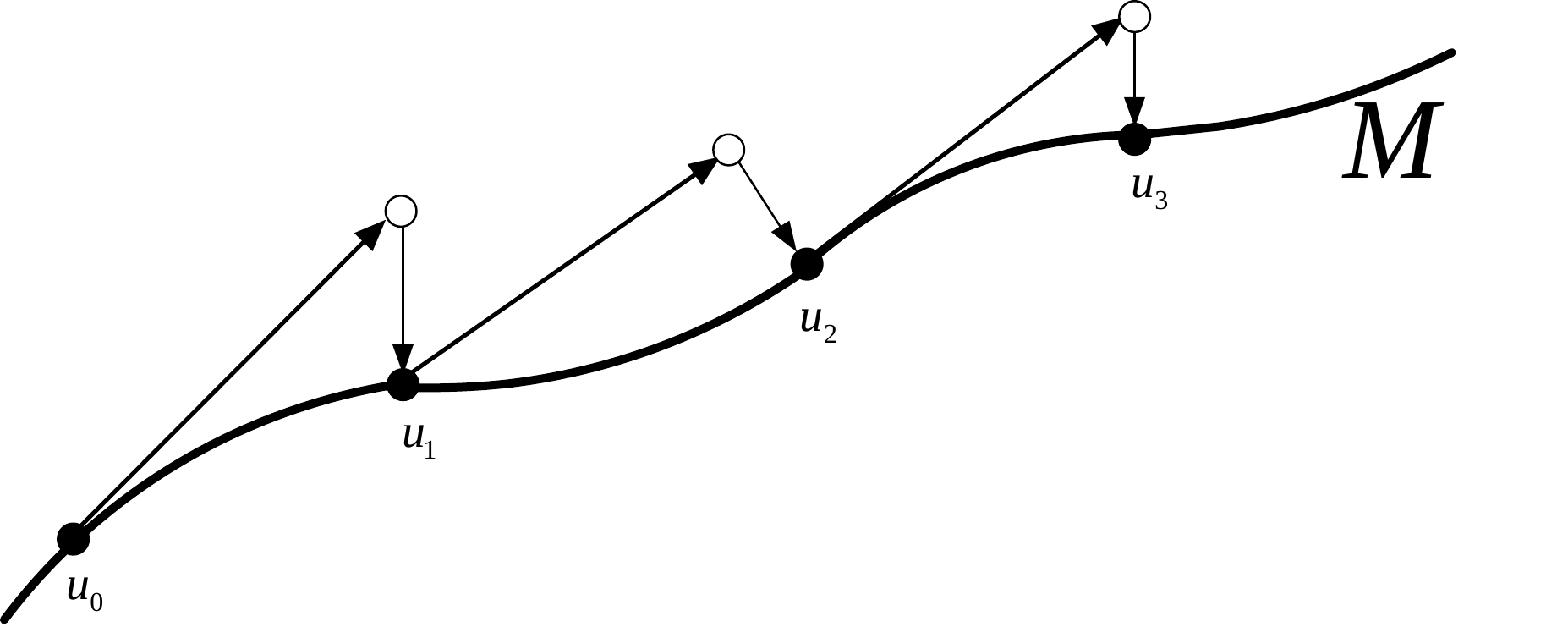}
\caption{\label{fig:proj} Schematic representation of the standard
  projection.}
\end{center}
\end{figure}

At iteration $k=0$, we set 
\begin{equation}
\label{eq:para-proj-0}
u_{n+1}^{k=0} = \mathcal{G}_{\Delta T}(u_{n}^{k=0}).
\end{equation}
For the next iterations, we simply couple the parareal
algorithm~\eqref{eq:scheme-parareal} with a projection step, as
schematically represented on Fig.~\ref{fig:proj}. We hence define the
solution at iteration $k+1$ from the solution at iteration $k$ by
\begin{equation}
\label{eq:para-proj-k}
u_{n+1}^{k+1} = \pi_{\cal M}\left[ \mathcal{G}_{\Delta T}(u_{n}^{k+1})
+ \mathcal{F}_{\Delta T}( u_{n}^{k}) 
- \mathcal{G}_{\Delta T}(u_{n}^{k}) \right].
\end{equation}

This leads to the following algorithm.
\begin{algorithm}
\label{algorithm-ps-proj}
Let $u_0$ be the initial condition. 
\begin{enumerate}
\item Initialization: set $u_{0}^{0} = u_{0}$, and sequentially compute 
$\left\{ u_{n}^{0} \right\}_{1 \leq n \leq N}$ by
$u_{n+1}^{0} = \mathcal{G}_{\Delta T} (u_{n}^0)$. 
\item Assume that, for some $k \geq 0$, the sequence $\left\{ u_{n}^{k}
  \right\}_{0 \leq n \leq N}$ is known. Compute 
$\left\{ u_{n}^{k+1} \right\}_{0 \leq n \leq N}$ at iteration $k+1$ by
the following steps: 
\begin{enumerate}
\item For all $0 \leq n \leq N-1$, compute 
$\mathcal{F}_{\Delta T}(u_{n}^{k}) - \mathcal{G}_{\Delta T}(u_{n}^{k})$ in parallel;
\item Set $u_{0}^{k+1} = u_{0}$;
\item Compute $\left\{ u_{n+1}^{k+1} \right\}_{0 \leq n \leq N-1}$
  sequentially by 
$$
u_{n+1}^{k+1} = \pi_{\cal M} \left[ 
\mathcal{G}_{\Delta T} (u_{n}^{k+1}) 
+ \mathcal{F}_{\Delta T} (u_{n}^{k}) - \mathcal{G}_{\Delta T}
(u_{n}^{k}) \right].
$$
\end{enumerate}
\end{enumerate}
\end{algorithm}

The additional complexity of the above algorithm, in comparison to the parareal
algorithm~(\ref{eq:scheme-parareal}), comes from solving the nonlinear
projection step~(\ref{eq:proj-eq}). Note however
that we expect to have a good initial guess for this problem, since the
solution at the previous parareal iteration is expected to be a coarse
approximation of the exact solution. We hence expect to solve the
nonlinear projection step with a few iterations of, say, a Newton
algorithm. We will see in the sequel that this is indeed the case. 

\subsection{Evaluation of the complexity of
  Algorithm~\ref{algorithm-ps-proj}}

We assess here the complexity of Algorithm~\ref{algorithm-ps-proj},
under the same assumptions as in Sections~\ref{sec:para4ham_complexity}
and~\ref{sec:para5ham_complexity}. Again, the complexity of the first
coarse propagation scales as $C_{\nabla} \ T /dT$.

Let us denote by $C_{\rm proj}$ the complexity of each projection
step. It depends on the number $m_{\rm proj}$ of 
iterations needed to solve this nonlinear problem. Each iteration
requires to evaluate the gradient of the energy, a task the complexity
of which is $C_{\nabla}$. Hence,
$$
C_{\rm proj} =  C_{\nabla} \, m_{\rm proj}.
$$ 
In all the simulations we have performed, $m_{\rm proj}$ is rather
small. This extra cost $C_{\rm proj}$ can be neglected with respect to
the cost of the fine propagation.

The complexity of each parareal iteration is thus again of the order of
$C_{\nabla} \, \Delta T/\delta t$. Denoting by $K_{\rm P/proj}$ the
number of parareal iterations, we obtain that the complexity of this
scheme is of the order of 
$$
C_{\rm P/proj} =  K_{\rm P/proj} \, C_{\nabla} \, \frac{\Delta T}{\delta t},
$$ 
as for the parareal scheme~\eqref{eq:scheme-parareal} and its symmetric
variant~(\ref{eq:sym-para-0})-(\ref{eq:sym-para-k}).
The speed-up $\dps \frac{T}{K_{\rm P/proj} \Delta T}$ is again equal to
the number of processors divided by the number of iterations, as long as
$m_{\rm proj}$ remains small. 

\subsection{Numerical results}
\label{sec:num_proj}

\begin{figure}[htbp]
\begin{center}
\includegraphics[width=9cm,angle=0,origin=0]{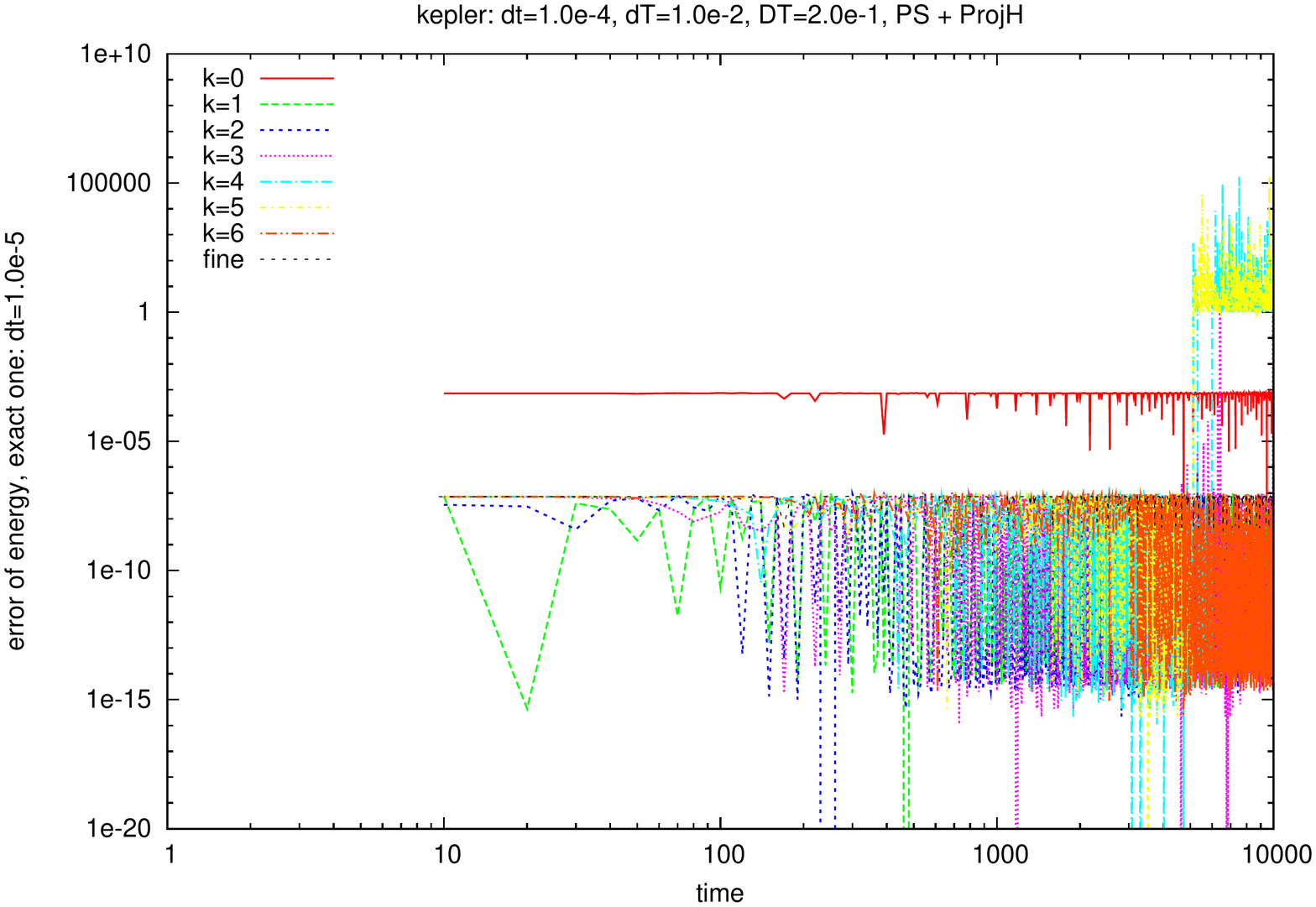}
\includegraphics[width=9cm,angle=0,origin=0]{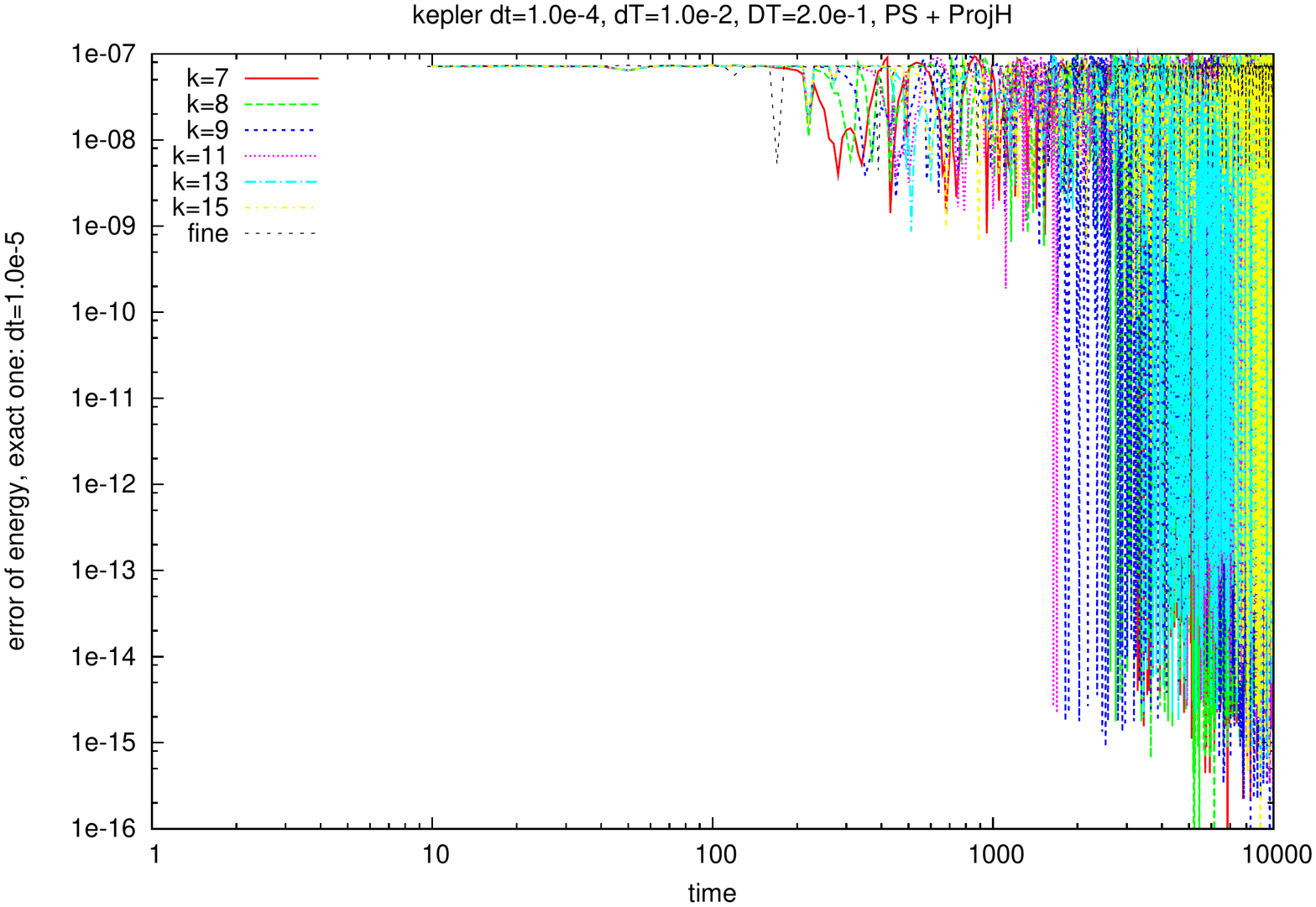}
\caption{ \label{fig:kepler-projH-ener}
Errors on the energy for Kepler problem, obtained by
Algorithm~\ref{algorithm-ps-proj} ($\delta t = 10^{-4}$, $dT = 0.01$,
$\Delta T = 0.2$).} 
\end{center}
\end{figure}

\begin{figure}[htbp]
\begin{center}
\includegraphics[width=9cm,angle=0,origin=0]{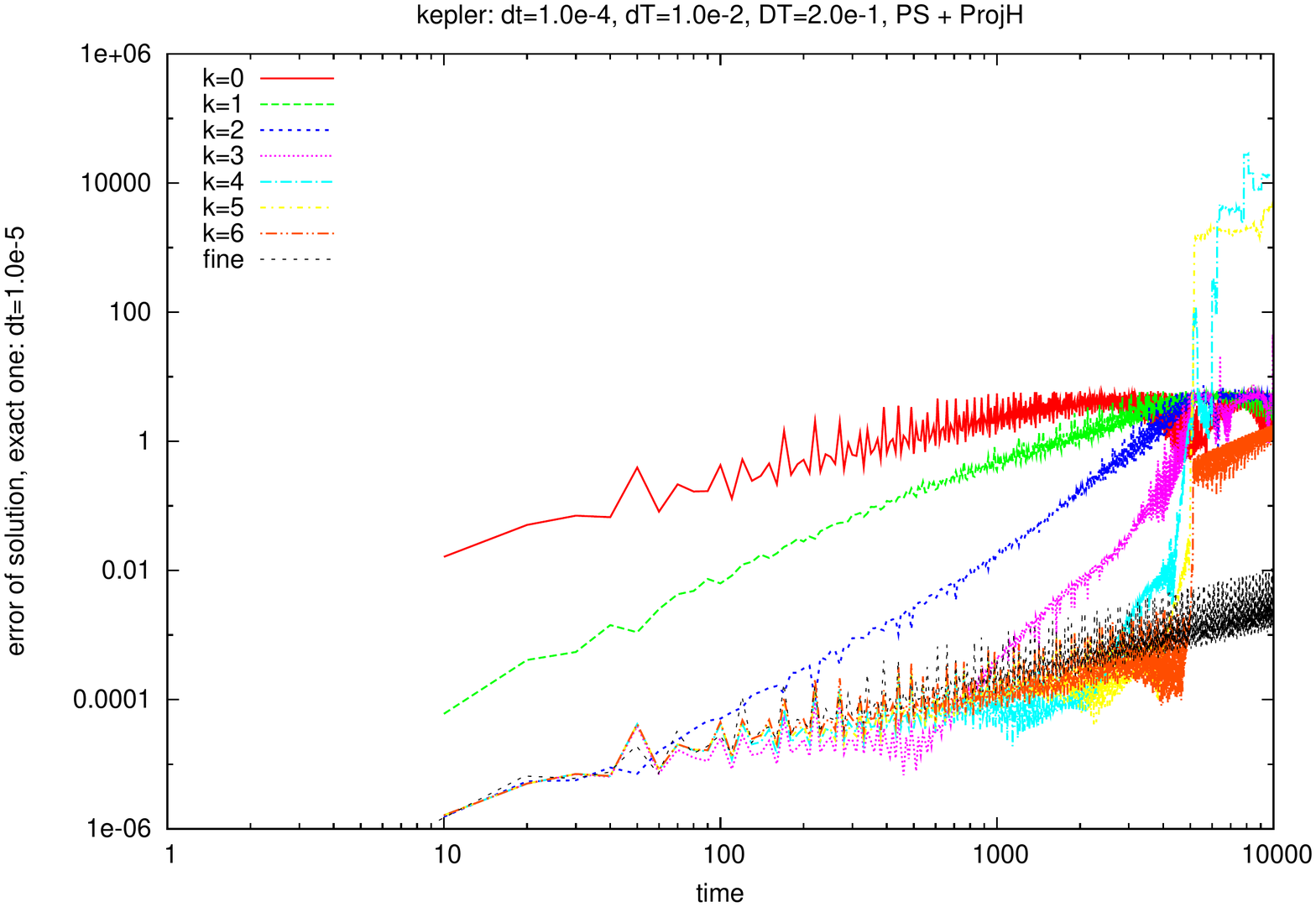}
\includegraphics[width=9cm,angle=0,origin=0]{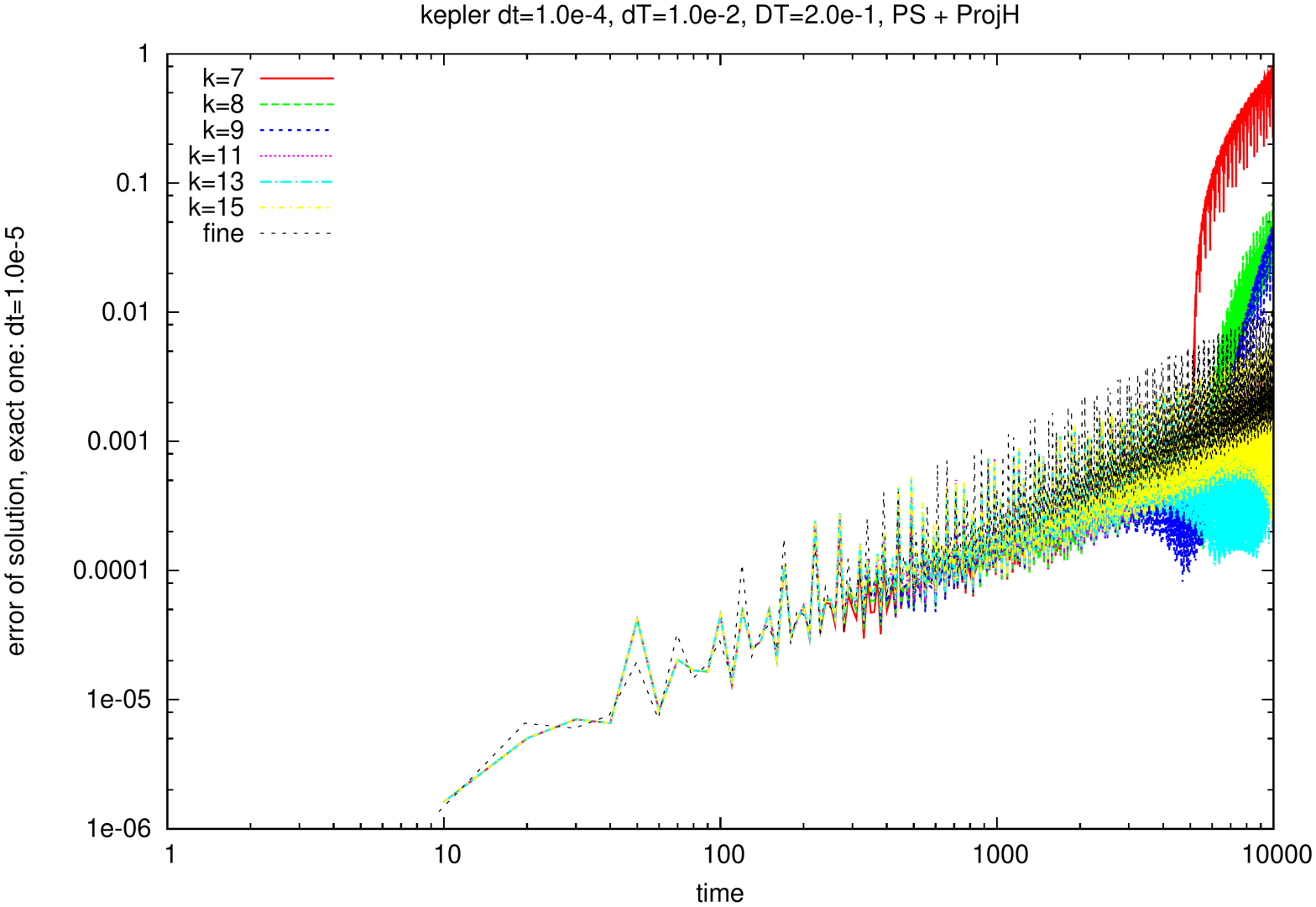}
\caption{ \label{fig:kepler-projH-pos}
Errors on the trajectory for Kepler problem, obtained by 
Algorithm~\ref{algorithm-ps-proj} ($\delta t = 10^{-4}$, $dT = 0.01$, $\Delta
T = 0.2$).} 
\end{center}
\end{figure}

We consider the Kepler problem associated to the
energy~\eqref{eq:H-kepler}, and we plot the relative
errors~\eqref{eq:err_H} on the energy preservation and the
errors~\eqref{eq:err_traj} on the trajectory 
on Figs.~\ref{fig:kepler-projH-ener} and~\ref{fig:kepler-projH-pos} 
respectively. 

The projection step~(\ref{eq:proj-eq}) is solved using a Newton
algorithm. Note yet that the energy is not exactly preserved, because
only a few steps of the Newton algorithm are performed. Indeed, 
the algorithm is terminated as soon as one of the following three
convergence criteria is satisfied:
\begin{description}
\item {C1:} the quantity {\tt err} is smaller than the tolerance {\tt tol},
\item {C2:} the maximum number $N_{\rm iter}^{\rm max}$ of iterations
  has been reached, 
\item {C3:} the error on the energy does not decrease,
\end{description}
where ${\tt err}$ is here the relative error on the
energy~\eqref{eq:err_H}. We have chosen to work with $N_{\rm iter}^{\rm
  max} = 2$. It turns out that the best choice for {\tt tol} is the
relative error on the energy preservation of the fine scheme. In the
current case, with $\delta t = 10^{-4}$, this corresponds to taking
${\tt tol} = 10^{-7}$.  

We observe an excellent
energy preservation, although we possibly stop the Newton algorithm
before convergence: the prescribed tolerance is reached over the complete
time range for all parareal iterations $k \geq 7$. This is
confirmed in Table~\ref{tab:occ},
where we show that, in the Newton algorithm, the criteria C1 (the
requested tolerance has been reached) is satisfied before the criterias
C2 and C3. In addition, for $k  
\geq 7$, the error~\eqref{eq:err_L} on the angular momentum preservation is
always smaller than $10^{-2}$, and it decreases down to $10^{-4}$ for
the iterations $k \geq 11$. 

\begin{table}[htbp]
\begin{center}
\begin{tabular}{|c|c|c|}
\hline
Criteria C1 & Criteria C2 & Criteria C3 \\
\hline 91.6 \% & 6.8 \% & 1.6 \% \\
\hline
\end{tabular}
\caption{To stop the projection procedure, we use three criteria. We gather
  here how often each criteria triggered the projection procedure to
  stop in Algorithm~\ref{algorithm-ps-proj}.
\label{tab:occ}}
\end{center}
\end{table}

We also observe that only 11 iterations are needed to obtain convergence of
the trajectory (that is, the trajectory obtained at iteration $k=11$ is
as accurate as the one given by the fine
scheme) on the time range $[0,10^4]$. 

With this algorithm, we are hence able to recover accurately the
trajectory on the complete time range $[0,10^4]$, for a moderate number
of iterations. Results are much better than with the other modifications
of the parareal algorithm that we have described so far (for the same
value of the parameters $\delta t$, $dT$ and $\Delta T$). 

\begin{figure}[htbp]
\begin{center}
\includegraphics[width=9cm,angle=0,origin=0]{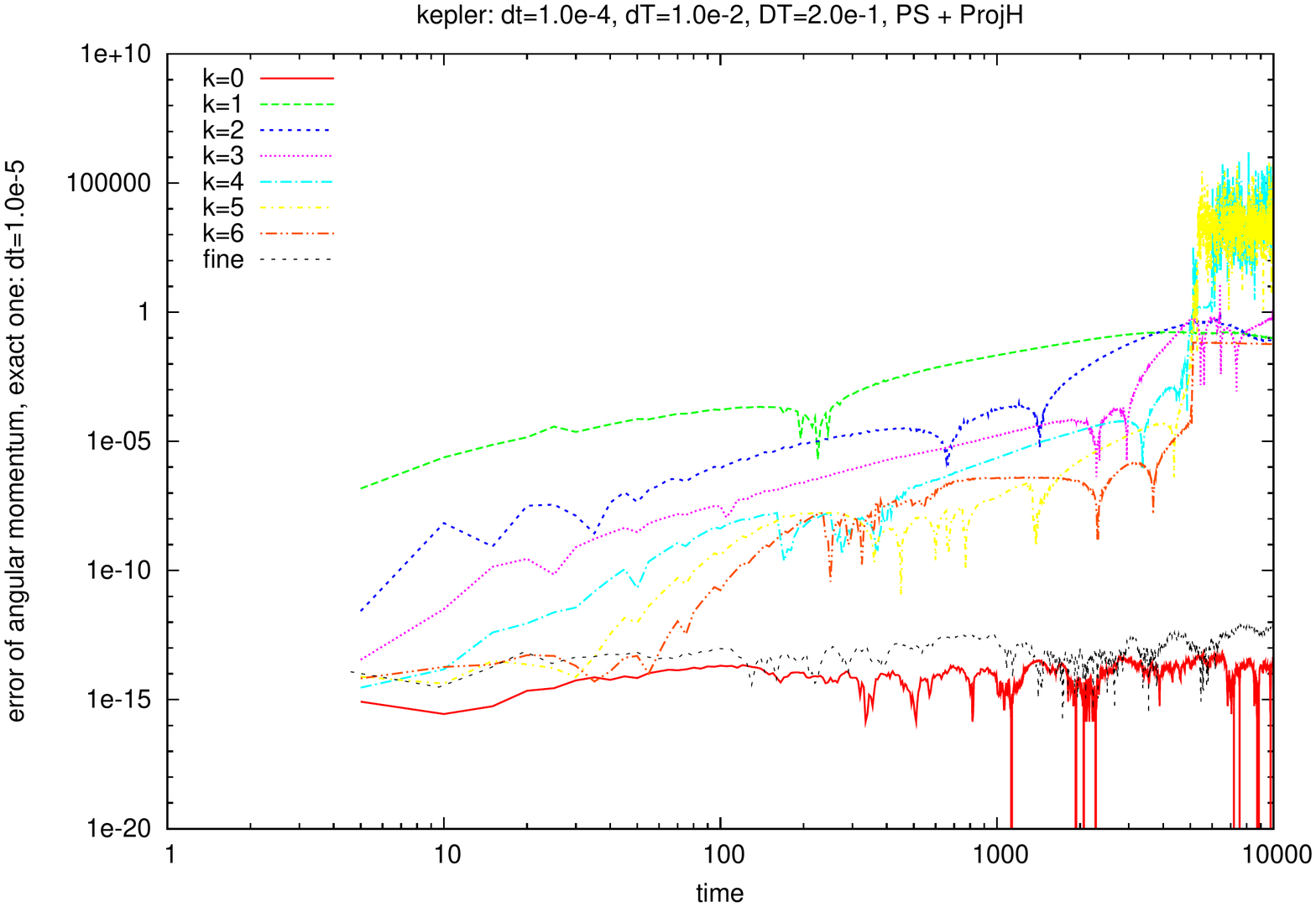}
\includegraphics[width=9cm,angle=0,origin=0]{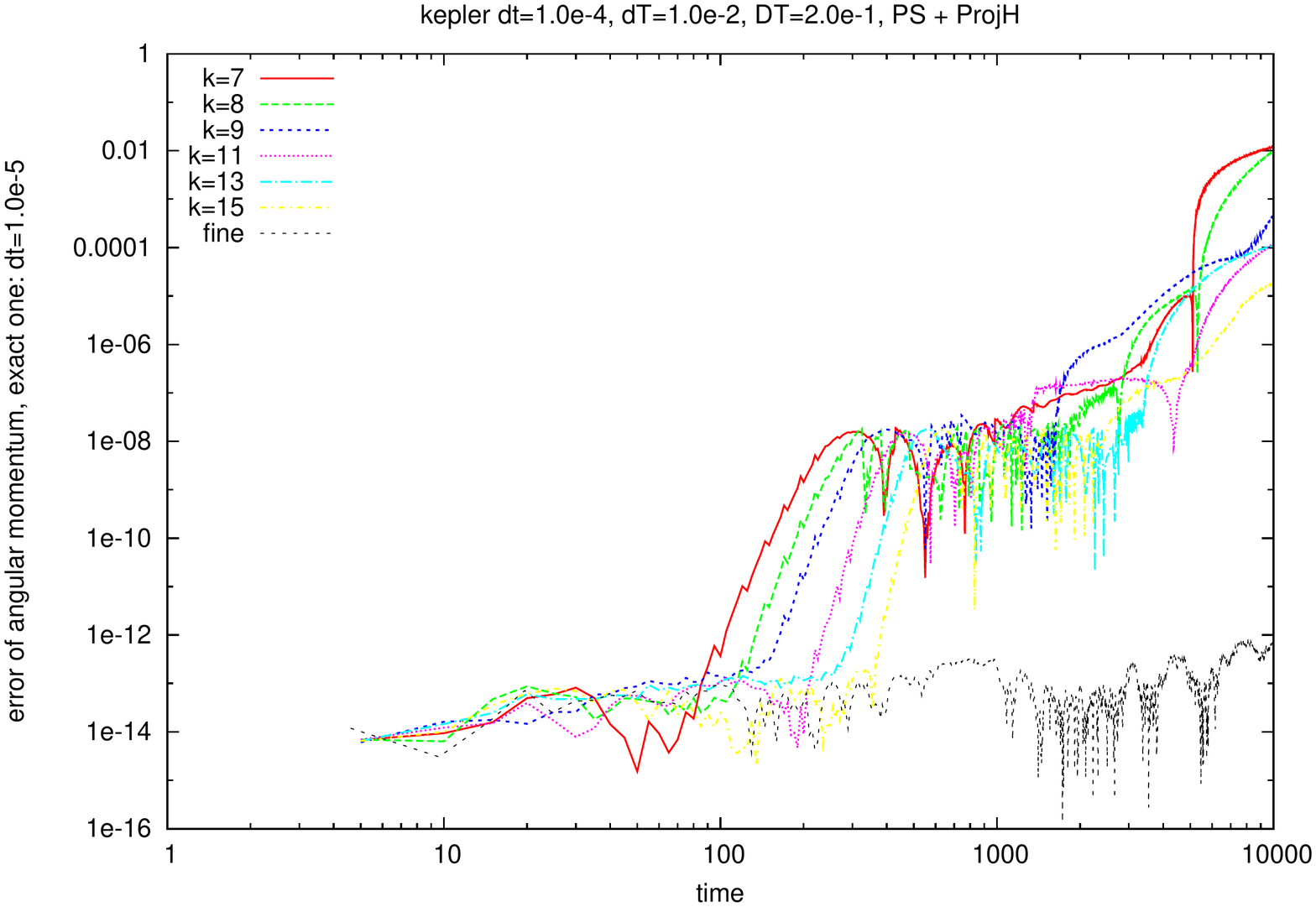}
\caption{ \label{fig:kepler-projH-L}
Errors on the angular momentum  for Kepler problem, obtained by 
Algorithm~\ref{algorithm-ps-proj} 
($\delta t = 10^{-4}$, $dT = 0.01$, $\Delta T = 0.2$).} 
\end{center}
\end{figure}

On Fig.~\ref{fig:kepler-projH-L}, we plot the relative 
errors on the angular momentum $L(q,p) = q \times p$, which is another
invariant of the Kepler problem. The relative error is defined by 
\begin{equation}
\label{eq:err_L}
{\rm err}_n^k = \frac{\left| L(q_n^k,p_n^k) - L(q_0,p_0) \right|}
{\left| L(q_0,p_0) \right|}.
\end{equation}
This quantity does not blow-up. We yet observe that it is not preserved
as well as would be the case with a geometric integrator. The rather
good behavior of this invariant is a consequence of the accuracy of
the trajectory. 

\begin{remark}
Note that the parameters chosen for the numerical simulation of the
harmonic oscillator and the Kepler system do not satisfy assumption
(\ref{assumption}). Our aim with these two toy-problems is indeed more
to understand the basic features of the variants of the parareal
algorithm than to be as efficient as possible. Our viewpoint will be
different when we consider a system of higher dimensionality, in
Section~\ref {sec:outer} below.
\end{remark}

\subsection{Considering more than one invariant}
\label{sec:2inv}

In this section, we again consider the case of the Kepler problem, and
now take into account that, for this specific system, we know another
invariant besides the energy, namely the angular momentum $L = q \times
p$. Define the manifold
$$
{\cal M}_2 = \left\{ (q,p) \in \R^{2d}; \ H(q,p) = H_0, \ L(q,p) = L_0 \right\},
$$
and consider the projection $\pi_{{\cal M}_2}$ onto that manifold, defined
by
$$
y = \pi_{{\cal M }_2}(\widetilde{y})
= \widetilde{y} + \lambda_1 \nabla H(\widetilde{y}) 
+ \lambda_2 \nabla L(\widetilde{y}),
$$
where $(\lambda_1,\lambda_2) \in \R^2$ is such that $H(y) = H_0$ and
$L(y) = L_0$. 

We couple this projection step with the parareal algorithm as follows
. We first define $u_{n}^{k=0}$ at iteration $k=0$
by~\eqref{eq:para-proj-0}: 
$$
u_{n+1}^{k=0} = \mathcal{G}_{\Delta T}(u_{n}^{k=0}).
$$
We next define the solution at iteration $k+1$ from the solution at iteration $k$ by
\begin{equation}
\label{eq:para-proj2-k}
u_{n+1}^{k+1} = \pi_{{\cal M}_2}\left[ \mathcal{G}_{\Delta T}(u_{n}^{k+1})
+ \mathcal{F}_{\Delta T}( u_{n}^{k}) 
- \mathcal{G}_{\Delta T}(u_{n}^{k}) \right].
\end{equation}
The algorithm is similar to Algorithm~\ref{algorithm-ps-proj} where
$\pi_{\cal M}$ is replaced by $\pi_{{\cal M}_2}$. 

\begin{figure}[htbp]
\begin{center}
\includegraphics[width=9cm,angle=0,origin=0]{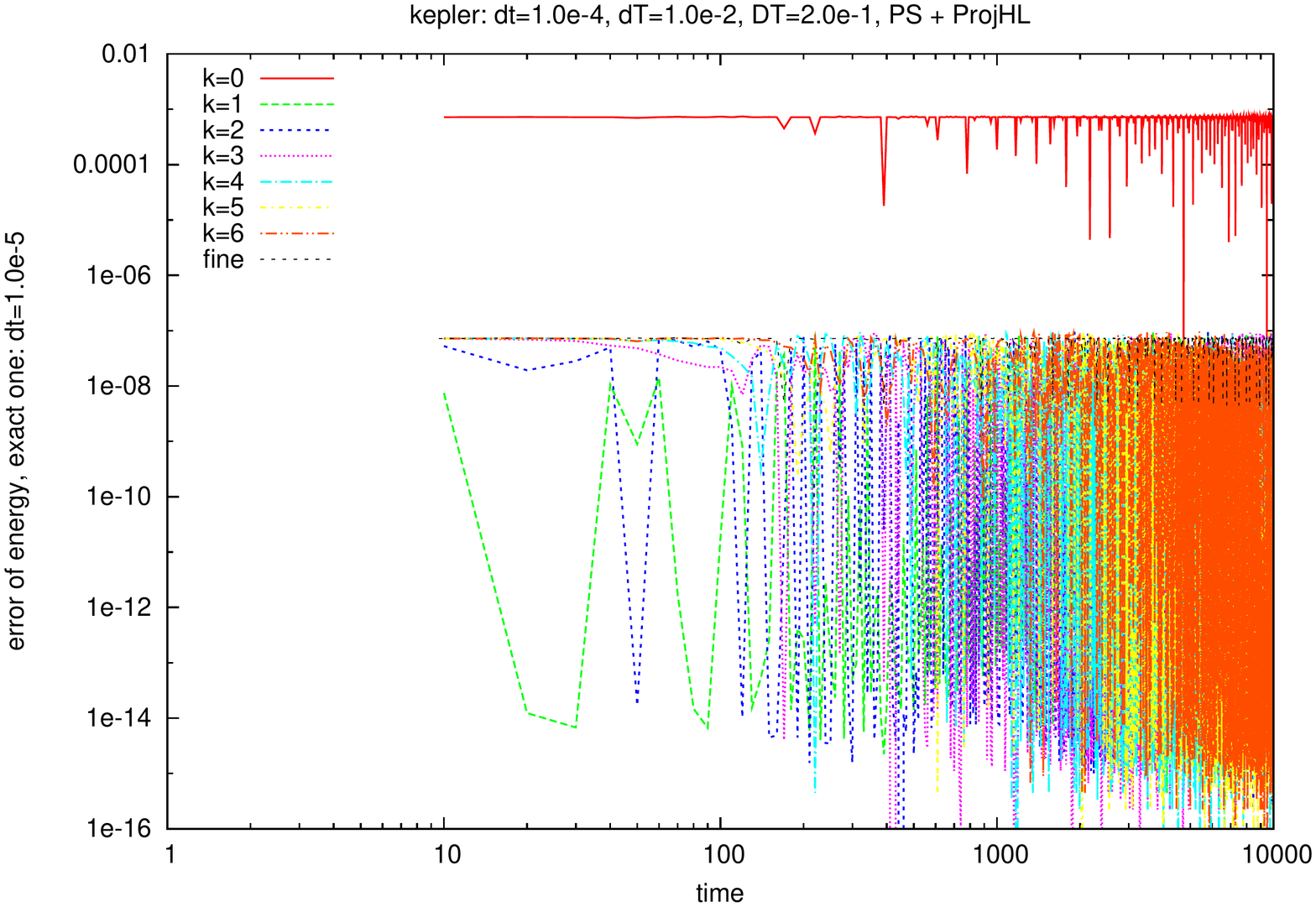}
\includegraphics[width=9cm,angle=0,origin=0]{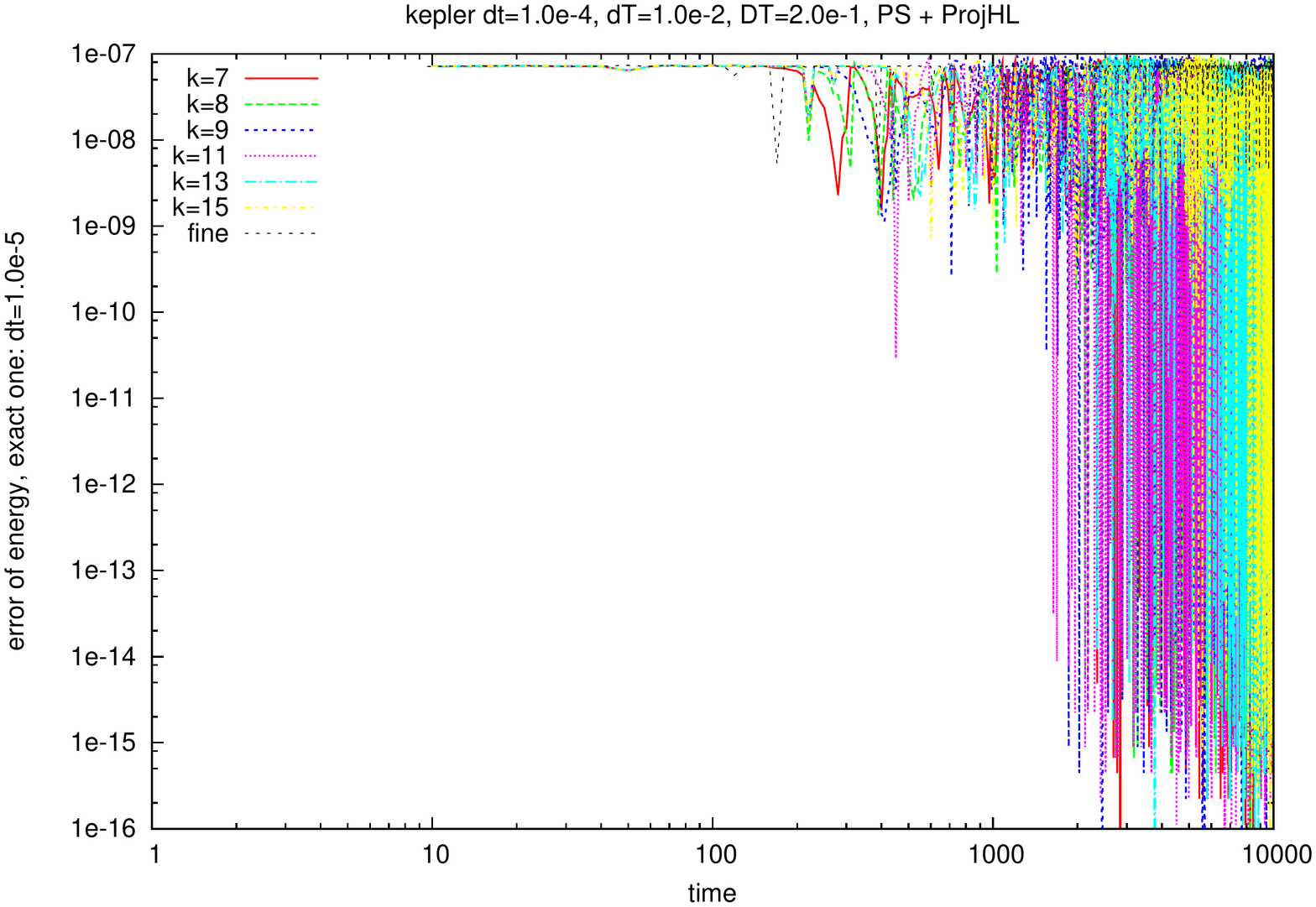}
\caption{ \label{fig:kepler-projHL-ener}
Errors on the energy for Kepler problem, obtained by the parareal
method~\eqref{eq:para-proj-0}-\eqref{eq:para-proj2-k} with projection on
the manifold ${\cal M}_2$ of constant energy and angular momentum
($\delta t = 10^{-4}$, $dT = 0.01$, $\Delta T = 0.2$).} 
\end{center}
\end{figure}

\begin{figure}[htbp]
\begin{center}
\includegraphics[width=9cm,angle=0,origin=0]{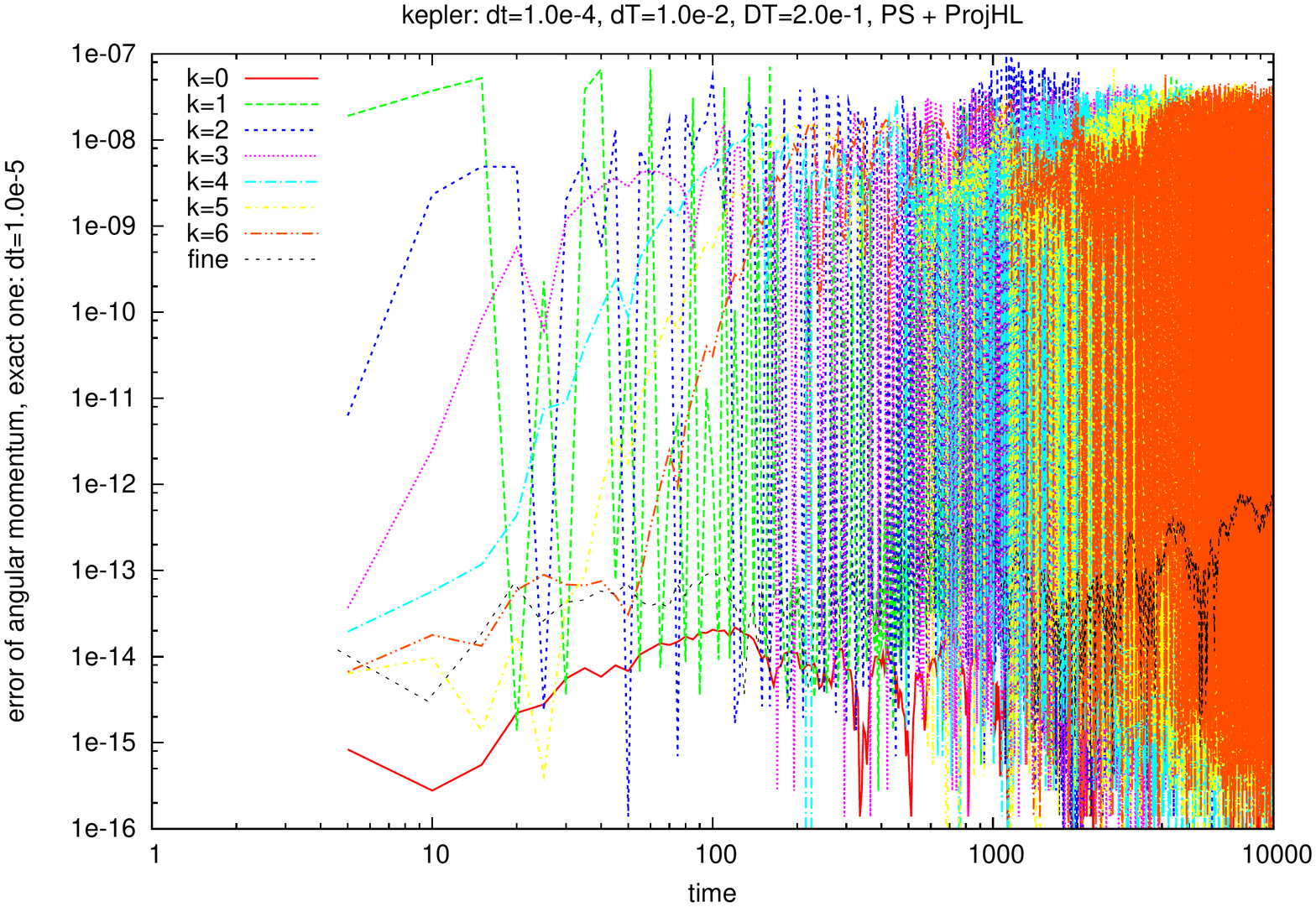}
\includegraphics[width=9cm,angle=0,origin=0]{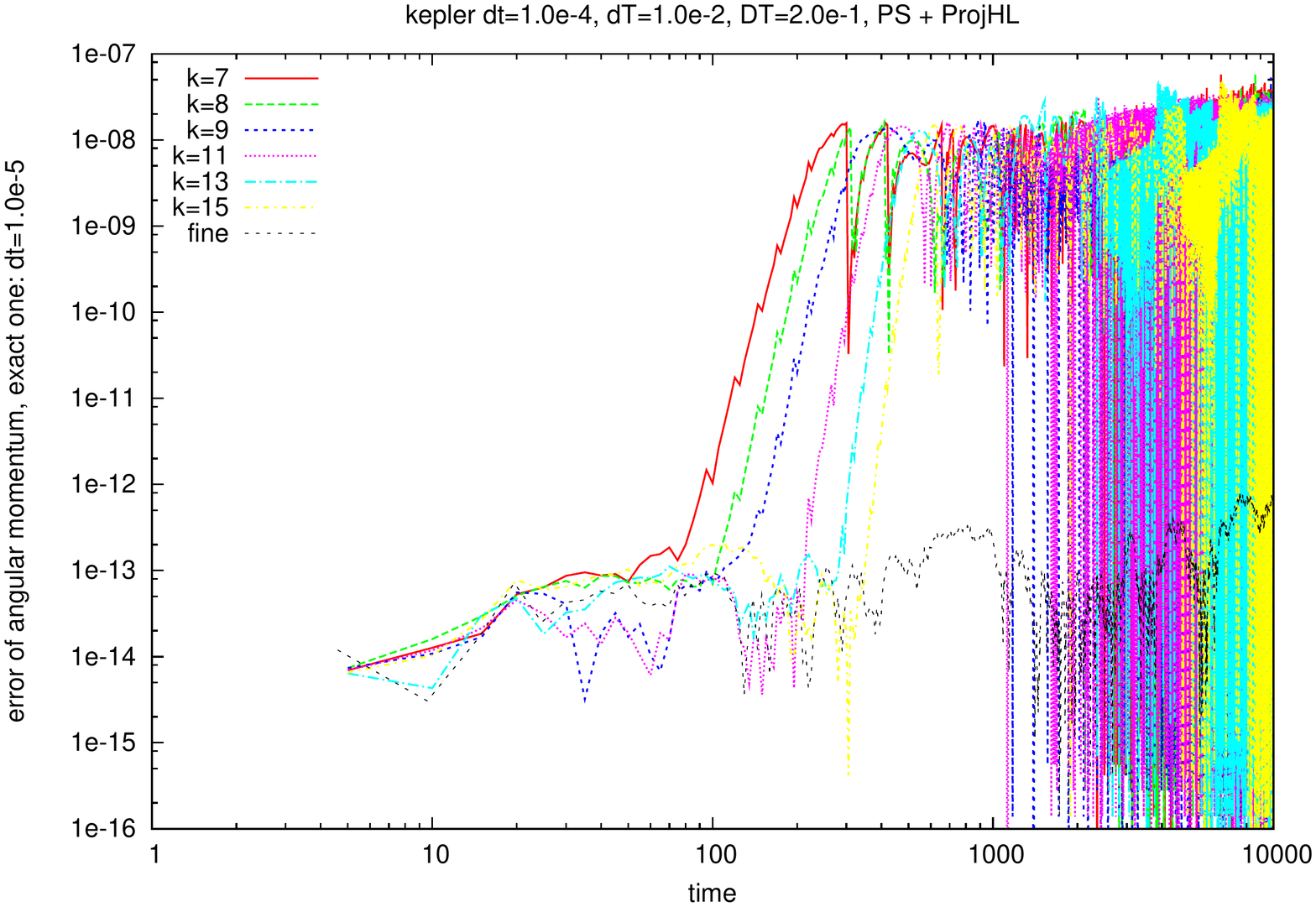}
\caption{ \label{fig:kepler-projHL-L}
Errors on the angular momentum  for Kepler problem, obtained by the parareal
method~\eqref{eq:para-proj-0}-\eqref{eq:para-proj2-k} with projection on
the manifold ${\cal M}_2$ of
constant energy and angular momentum ($\delta t = 10^{-4}$, $dT = 0.01$,
$\Delta T = 0.2$).}
\end{center}
\end{figure}

\begin{figure}[htbp]
\begin{center}
\includegraphics[width=9cm,angle=0,origin=0]{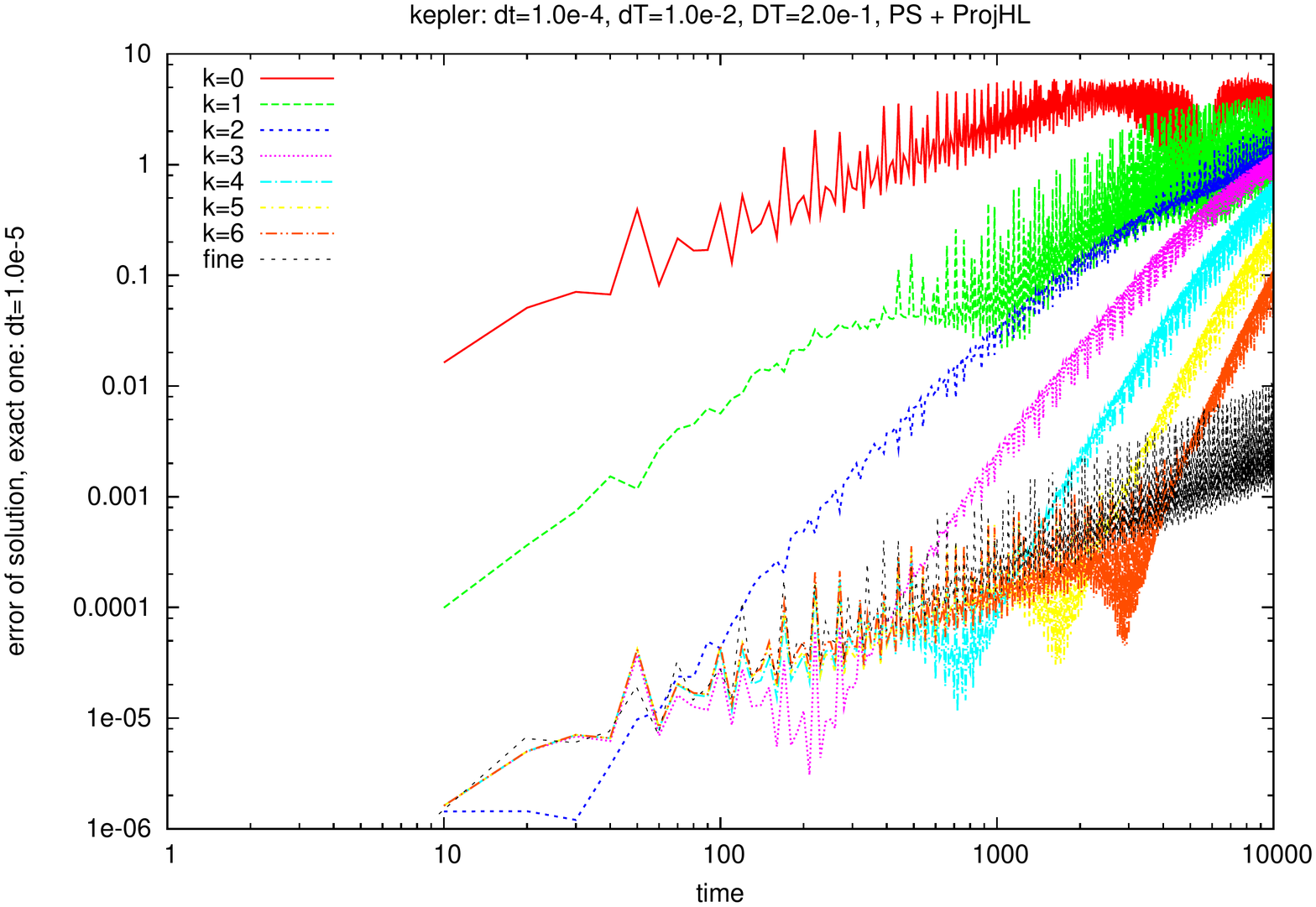}
\includegraphics[width=9cm,angle=0,origin=0]{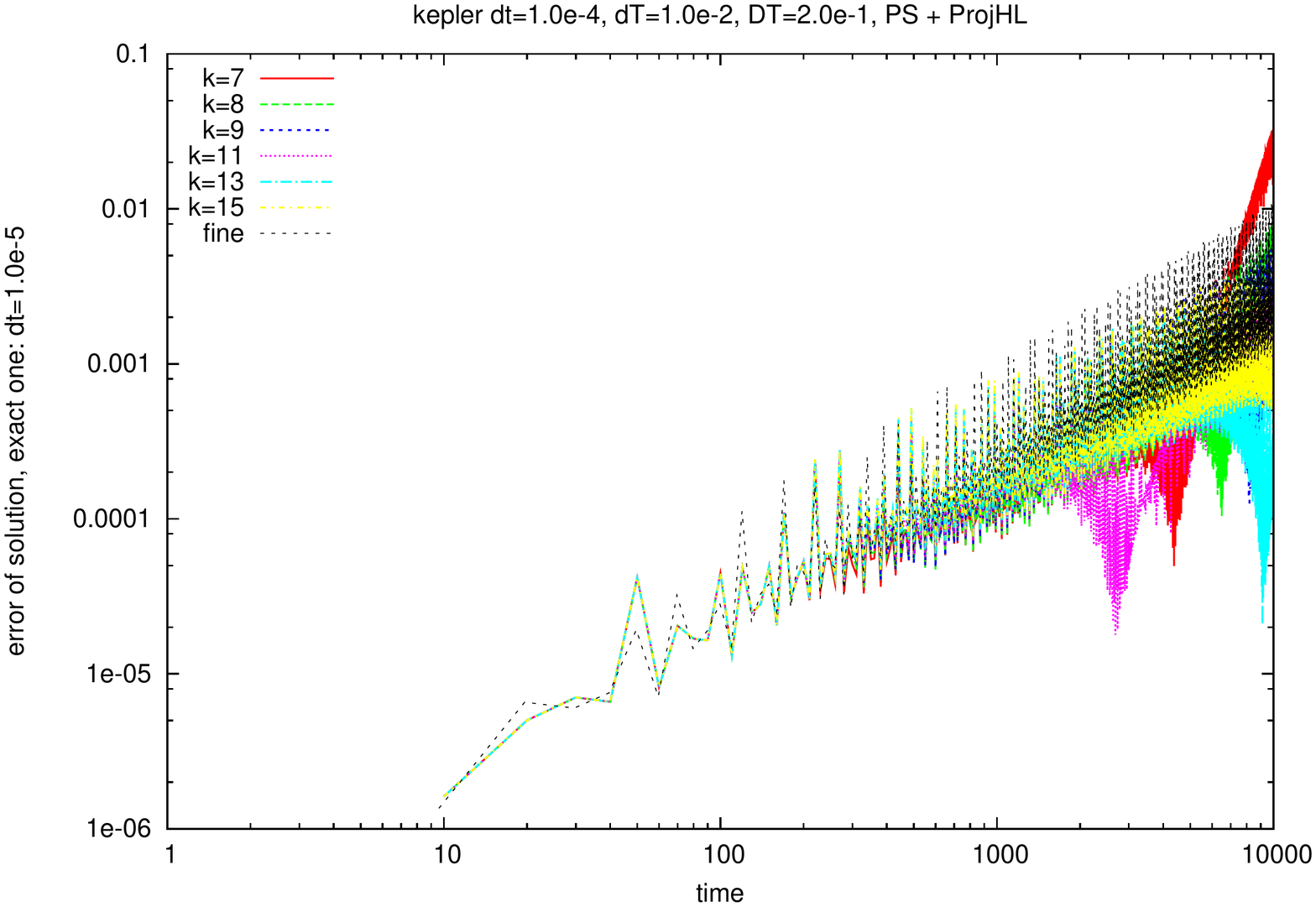}
\caption{ \label{fig:kepler-projHL-pos}
Errors on the trajectory for Kepler problem, obtained by the parareal
method~\eqref{eq:para-proj-0}-\eqref{eq:para-proj2-k} with projection on
the manifold ${\cal M}_2$ of constant energy and angular momentum
($\delta t = 10^{-4}$, $dT = 0.01$, $\Delta T = 0.2$).} 
\end{center}
\end{figure}

Let us now integrate the Kepler dynamics with this algorithm (we again
use a Newton algorithm to compute the Lagrange multipliers, with the
same three stopping criteria as in the previous section). We plot the relative
errors~\eqref{eq:err_H} on the energy preservation, the relative
errors~\eqref{eq:err_L} on the angular momentum preservation and the
errors~\eqref{eq:err_traj} on the trajectory, 
on
Figs.~\ref{fig:kepler-projHL-ener}, \ref{fig:kepler-projHL-L}
and~\ref{fig:kepler-projHL-pos} respectively.

We observe qualitatively the same behavior as for
Algorithm~\ref{algorithm-ps-proj}, which is based
on projection on the constant energy manifold. Only 8 iterations are
needed to obtain convergence of 
the trajectory, rather than 11 before. Of course, we observe an excellent
energy and angular momentum preservations: the error for both quantities
is lower than the prescribed tolerance for all parareal iterations $k
\geq 1$.  

In summary, this algorithm, that projects the trajectory on the manifold
${\cal M}_2$ where {\em both} energy {\em and} angular momentum are
constant, yields results that are slightly more accurate than those
obtained with Algorithm~\ref{algorithm-ps-proj}, where the trajectory is
projected on the manifold ${\cal M}$ where {\em only} the energy is constant. 
However, the behavior of the trajectory is not much
improved. From this prototypical analysis, it does not seem worth to
project the trajectory on the manifold where several invariants (besides
the energy) are constant, especially if we bear in mind that  the
determination of these other invariants might be a difficult task in the
general case. 

\section{Symmetric parareal algorithm with symmetric projection}
\label{sec:sym_proj}

We have introduced in Section~\ref{sec:proj} a parareal algorithm using
some projection step. Of course, just as the original parareal
algorithm~\eqref{eq:scheme-parareal}, this algorithm is neither symplectic nor
symmetric (see Section~\ref{sec:para4ham}). 

Yet, it is well known that it is possible to couple a symmetric
algorithm with an appropriate projection step, while still keeping the algorithm
symmetric. We first briefly recall this idea before using it in our
context.

\subsection{Symmetric projection}
\label{sec:sym_proj_gene}

\begin{figure}[htbp]
\begin{center}
\includegraphics[width = 6cm,angle =0]{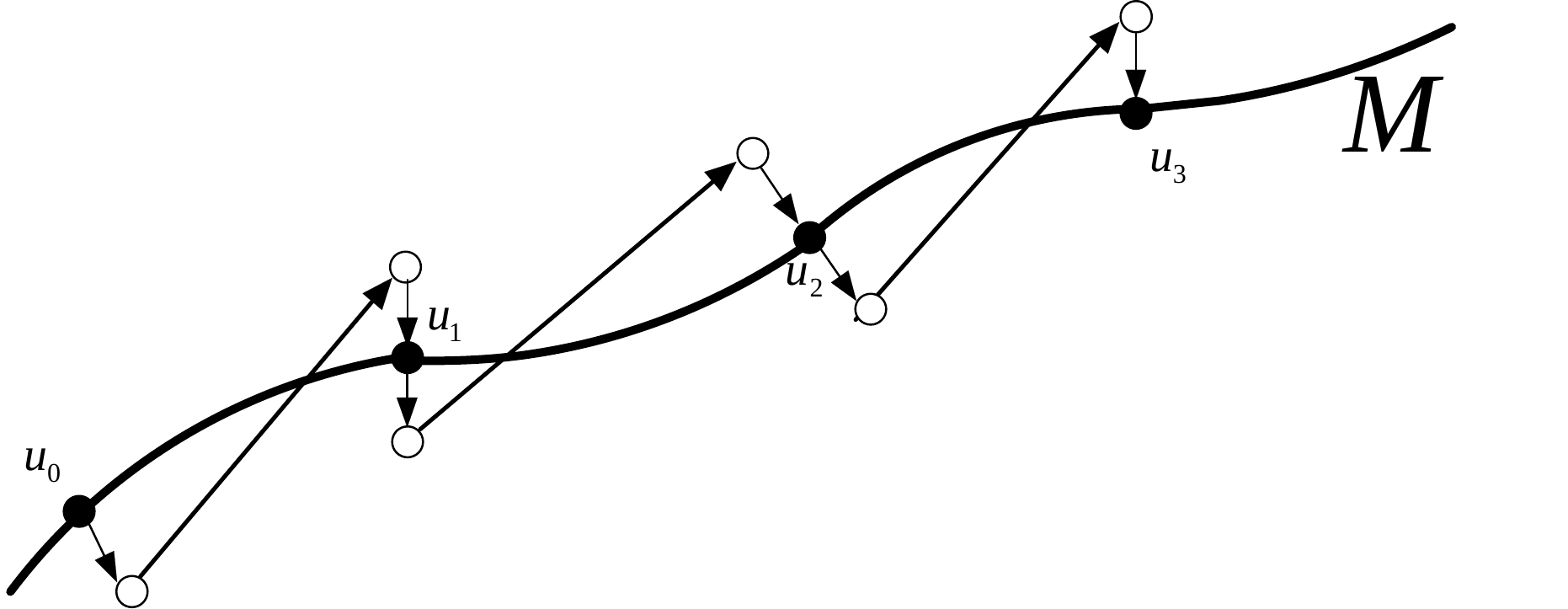}
\caption{\label{fig:symm-proj} Symmetric projection on the manifold ${\cal
    M}$.}
 \end{center}
\end{figure}

Consider the constant energy manifold
$$
{\cal M} = \left\{ y = (q,p) \in \R^{2d}; \ H(q,p) = H_0 \right\},
$$
and assume that we have at hand a symmetric algorithm with numerical
flow $\Psi_{\Delta T}$. 
Consider now the following algorithm (see~Fig.~\ref{fig:symm-proj} for a
schematic representation). Assume that $y_n \in {\cal M}$. The
approximation $y_{n+1}$ is next defined by the set of equations 
\begin{equation}
\label{eq:symm-proj}
\left\{
\begin{array}{rcl}
\widetilde{y}_n &=& y_n +  \mu \nabla H(y_n),
\\
\widetilde{y}_{n+1} &=& \Psi_{\Delta T}(\widetilde{y}_n),
\\
y_{n+1} &=& \widetilde{y}_{n+1} + \mu \nabla H(y_{n+1}),
\end{array}
\right.
\end{equation}
with $\mu$ chosen such that $y_{n+1} \in {\cal M}$. 
Consequently, $\mu$ and $y_{n+1}$ satisfy the equations
\begin{equation}
\label{eq:symm-proj2}
\left\{
\begin{array}{rcl}
y_{n+1} - \Psi_{\Delta T}(y_n + \mu \nabla H(y_n)) - \mu \nabla H(y_{n+1}) 
& = & 0,
\\ 
H(y_{n+1}) & = & H_0.
\end{array}
\right.
\end{equation}
We denote by $\pi_{{\cal M}}^{\rm sym}$ this algorithm:
$$
y_{n+1} = \pi_{{\cal M}}^{\rm sym}(y_n) \quad \text{if and only
  if~\eqref{eq:symm-proj} is satisfied}.
$$
As the {\em same} Lagrange multiplier $\mu$ is used in the first and
third lines of~\eqref{eq:symm-proj}, and as $\Psi_{\Delta T}$ is
symmetric, one can easily see that $\pi_{{\cal M}}^{\rm sym}$ is a
symmetric algorithm.

To solve the nonlinear equations~(\ref{eq:symm-proj2}), we first recast
them in the form ${\cal S}(y_{n+1},\mu) = 0$, with 
$$
{\cal S}(y_{n+1},\mu) = 
\left(
\begin{array}{c}
{\cal S}_1(y_{n+1},\mu)
\\ 
{\cal S}_2(y_{n+1},\mu)
\end{array}
\right)
=
\left(
\begin{array}{c}
y_{n+1} - \Psi_{\Delta T}(y_n + \mu \nabla H(y_n)) - \mu \nabla H(y_{n+1}) 
\\ 
H(\Psi_{\Delta T}(y_n + \mu \nabla H(y_n)) + \mu \nabla H(y_{n+1})) - H_0
\end{array}
\right).
$$
We use a Newton-like algorithm, where the jacobian matrix of ${\cal S}$
is approximated by
\begin{equation}
\label{eq:jac_matrix}
\text{Jac } {\cal S}(y_{n+1},\mu) \approx
\left(
\begin{array}{cc}
I & - \nabla H(y_n) - \nabla H(y_{n+1})
\\
0 & \nabla H(\hat{y}_{n+1})^T (\nabla H(y_n) + \nabla H(y_{n+1}))
\end{array}
\right),
\end{equation}
with $\hat{y}_{n+1} = \Psi_{\Delta T}(y_n + \mu \nabla H(y_n)) + \mu
\nabla H(y_{n+1})$. 
To stop the Newton algorithm, we again use the three criteria C1, C2,
and C3 presented above, where {\tt err} is here the relative residu,
$$
{\tt err} = \frac{\left| S_1(y_{n+1},\mu) \right|}{|y_{n+1}|} 
+ 
\frac{\left| S_2(y_{n+1},\mu) \right|}{|H_0|}.
$$

\bigskip

An alternative to the algorithm~\eqref{eq:symm-proj} is to use $\nabla
H(\widetilde{y}_{n+1})$ rather than $\nabla H(y_{n+1})$ in the third
line of~\eqref{eq:symm-proj}. This yields the following algorithm: 
\begin{equation}
\label{eq:q-symm-proj}
\left\{
\begin{array}{rcl}
\widetilde{y}_n &=& y_n +  \mu \nabla H(y_n),
\\
\widetilde{y}_{n+1} &=& \Psi_{\Delta T}(\widetilde{y}_n),
\\
y_{n+1} &=& \widetilde{y}_{n+1} + \mu \nabla H(\widetilde{y}_{n+1}),
\end{array}
\right.
\end{equation}
with again $\mu$ chosen such that $y_{n+1} \in {\cal M}$. The advantage
in comparison to~\eqref{eq:symm-proj} is that the nonlinear system to
be solved is easier. Indeed, once $\mu$ is identified, the approximation
$y_{n+1}$ can be obtained in an explicit fashion. The Lagrange
multiplier $\mu$ now solves the equation ${\cal S}(\mu) = 0$, with
$$
{\cal S}(\mu) 
=
H \left[ \Psi_{\Delta T} \left( y_n +  \mu \nabla H(y_n) \right) + \mu \nabla H
\left( y_n +  \mu \nabla H(y_n) \right) \right] - H_0,
$$
which is a {\em scalar} nonlinear equation, to be solved for a {\em
  scalar} unknown.  
We denote by $\pi_{{\cal M}}^{\rm q-sym}$ this algorithm, that we call
the {\em quasi-symmetric projection algorithm}:
$$
y_{n+1} = \pi_{{\cal M}}^{\rm q-sym}(y_n) \quad \text{if and only
  if~\eqref{eq:q-symm-proj} is satisfied}.
$$
Note that this algorithm is {\em not} symmetric. It will yet turn out
to be inexpensive, and to yield interesting results. 

In practice, the nonlinear problem ${\cal S}(\mu) = 0$ is solved using a
Newton algorithm, where the derivative of ${\cal S}$ is approximated by
$$
{\cal S}'(\mu) \approx 
\nabla H(y_{n+1})^T (\nabla H(y_n) + \nabla H (\widetilde{y}_n)).
$$
The algorithm is stopped as soon as one of the three criteria C1, C2, C3
is satisfied (in C1, ${\tt err}$ is the relative error on the energy). 

\subsection{Application to the parareal setting}

We now wish to use the above two projection procedures, namely the
symmetric projection algorithm 
$\pi_{{\cal M}}^{\rm sym}$ and the quasi-symmetric projection algorithm
$\pi_{{\cal M}}^{\rm q-sym}$, in the parareal context. To this end,
we return to the formalism used in
Section~\ref{sec:design_parasym}.

We first choose a number $K$ of iterations for the
parareal algorithm, and set
$$
U_n := (u_n^0, u_n^1, \cdots, u_n^K).
$$
The symmetric parareal algorithm is given by~\eqref{eq:sym-para}, that we
here write as
$$
U_{n+1} = \Psi^{\rm sym}_{\Delta T}(U_n),
$$ 
where the map $\Psi^{\rm sym}_{\Delta T}$ is symmetric in the classical
sense. We next introduce $K$ energies, namely, for any $1 \leq k \leq K$,
$$
H_k(U_n) := H(u_n^k)
$$
and we consider a symmetric algorithm, based on the symmetric
map $\Psi^{\rm sym}_{\Delta T}$ and on a symmetric projection on the manifold where all
the energies $H_k$ are preserved. Following~\eqref{eq:symm-proj}, this
algorithm writes
\begin{eqnarray*}
\widetilde{U}_n &=& U_n + \sum_{k=1}^K \mu_k \nabla H_k(U_n)
\\
\widetilde{U}_{n+1} &=& \Psi^{\rm sym}_{\Delta T}(\widetilde{U}_n)
\\
U_{n+1} &=& \widetilde{U}_{n+1} + \sum_{k=1}^K \mu_k \nabla H_k(U_{n+1}),
\end{eqnarray*}
where the Lagrange multipliers $\mu_k$ are such that $H_k(U_{n+1}) =
H_0$ for any $1 \leq k \leq K$. 

Since $\nabla H_k(U_n) = (0,\ldots,0,\nabla H(u_n^k),0,\ldots,0)$, the
above equations read
\begin{eqnarray*}
\widetilde{u}_n^0 &=& u_n^0 \quad \text{and} \quad
\widetilde{u}_n^k = u_n^k + \mu_k \nabla H(u_n^k) \text{ for any $1
  \leq k \leq K$},
\\
\widetilde{U}_{n+1} &=& \Psi^{\rm sym}_{\Delta T}(\widetilde{U}_n),
\\
u_{n+1}^0 &=& \widetilde{u}_{n+1}^0 \quad \text{and} \quad
u_{n+1}^k = \widetilde{u}_{n+1}^k + \mu_k \nabla H(u_{n+1}^k) \text{ for any $1
  \leq k \leq K$}.
\end{eqnarray*}

This yields the following algorithm, that we call \emph {symmetric parareal
algorithm with symmetric projection} in the sequel:
\begin{algorithm}
\label{algorithm-sps-symm-proj}
Let $u_0$ be the initial condition. 
\begin{enumerate}
\item Initialization: set $u_{0}^{0} = u_{0}$, and sequentially compute
$\left\{ u_{n+1/2}^0 \right\}_{0 \leq n \leq N-1}$ and 
$\left\{ u_{n+1}^0 \right\}_{0 \leq n \leq N-1}$ as
$$
u_{n+1/2}^{0} = \mathcal{G}_{-\Delta T/2}^{-1}(u_{n}^0), 
\quad
u_{n+1}^{0} = \mathcal{G}_{\Delta T/2}(u_{n+1/2}^0).
$$
Set $\widetilde{u}_{n+1/2}^0 = u_{n+1/2}^0$. 
\item Assume that, for some $k \geq 0$, the sequence
$\left\{ \widetilde{u}_{n+1/2}^k \right\}_{0 \leq n \leq N-1}$ is known. Compute
the sequences $\left\{ u_{n}^{k+1} \right\}_{0 \leq n \leq N}$ and 
$\left\{ \widetilde{u}_{n+1/2}^{k+1} \right\}_{0 \leq n \leq N-1}$ at iteration
$k+1$ by the following steps:
\begin{enumerate}
 \item For all $0 \leq n \leq N-1$, compute in parallel 
$\mathcal{F}_{-\Delta T/2}(\widetilde{u}_{n+1/2}^{k})$, 
$\mathcal{F}_{\Delta T/2}(\widetilde{u}_{n+1/2}^{k})$,
$\mathcal{G}_{-\Delta T/2}(\widetilde{u}_{n+1/2}^{k})$ and 
$\mathcal{G}_{\Delta T/2}(\widetilde{u}_{n+1/2}^{k})$;
\item Set $u_{0}^{k+1} = u_{0}$;
\item For $0 \leq n \leq N-1$, sequentially compute
  $\widetilde{u}_{n+1/2}^{k+1}$ and $u_{n+1}^{k+1}$ as
\begin{eqnarray*}
\widetilde{u}_n^{k+1} &=& u_n^{k+1} + \mu_{k+1} \nabla H(u_n^{k+1})
\\
\widetilde{u}_{n+1/2}^{k+1} &=&
\mathcal{G}_{-\Delta T/2}^{-1}\left[\widetilde{u}_{n}^{k+1} -
\mathcal{F}_{-\Delta T/2}(\widetilde{u}_{n+1/2}^{k}) + 
\mathcal{G}_{-\Delta T/2}(\widetilde{u}_{n+1/2}^{k})\right]
\\
\widetilde{u}_{n+1}^{k+1} &=& \mathcal{G}_{\Delta T/2}(\widetilde{u}_{n+1/2}^{k+1})
+ \mathcal{F}_{\Delta T/2}(\widetilde{u}_{n+1/2}^{k}) - 
\mathcal{G}_{\Delta T/2}(\widetilde{u}_{n+1/2}^{k})
\\
u_{n+1}^{k+1} &=& \widetilde{u}_{n+1}^{k+1} + \mu_{k+1} \nabla H(u_{n+1}^{k+1})
\end{eqnarray*}
where $\mu_{k+1}$ is such that $H(u_{n+1}^{k+1}) = H_0$. 
\end{enumerate}
\end{enumerate}
\end{algorithm}

We observe that all the expensive computations (involving the fine
propagator) can again be performed in parallel. Like in the symmetric
parareal algorithm, 
the map that defines $u_{n+1}^{k+1}$ from $u_n^{k+1}$ is not symmetric
in the classical sense,
but the map that defines $(u_{n+1}^0, u_{n+1}^1, \cdots, u_{n+1}^K)$
from $(u_n^0, u_n^1, \cdots, u_n^K)$ is symmetric. 

\begin{remark}
\label{rem:quasi-sym}
Instead of using the symmetric projection
scheme~\eqref{eq:symm-proj}, we can use the quasi-symmetric
projection scheme~\eqref{eq:q-symm-proj}. In the above algorithm, it
amounts to replacing the last line in 2c), namely 
$$
u_{n+1}^{k+1} = \widetilde{u}_{n+1}^{k+1} + \mu_{k+1} \nabla H(u_{n+1}^{k+1}),
$$
by 
$$
u_{n+1}^{k+1} = \widetilde{u}_{n+1}^{k+1} + \mu_{k+1} \nabla H(\widetilde{u}_{n+1}^{k+1}),
$$
where again $\mu_{k+1}$ is such that $H(u_{n+1}^{k+1}) = H_0$. This
algorithm is called {\em symmetric parareal algorithm with quasi-symmetric
projection}. 
\end{remark}
 
\subsection{Evaluation of the complexity of Algorithm~\ref{algorithm-sps-symm-proj}}

Under the same hypothesis as in Section~\ref{sec:para4ham_complexity},
and the same implementation following~\cite{batchelor}, we obtain that
the complexity of Algorithm~\ref{algorithm-sps-symm-proj} is of the
order of 
$$
C_{\rm SP/sym-proj} =  K_{\rm SP/sym-proj} \ C_{\nabla} \, 
\frac{\Delta T}{\delta t},
$$ 
where $K_{\rm SP/sym-proj}$ is the number of parareal iterations. The
speed-up is hence equal to $\dps \frac{T}{K_{\rm SP/sym-proj} \Delta
  T}$, which is the number of processors divided by the number of
iterations. A similar analysis holds for the symmetric parareal scheme
with {\em quasi-symmetric} projection. 

\subsection{Numerical results for the Kepler problem }

\begin{figure}[htbp]
\begin{center}
\includegraphics[width=9cm,angle=0,origin=0]{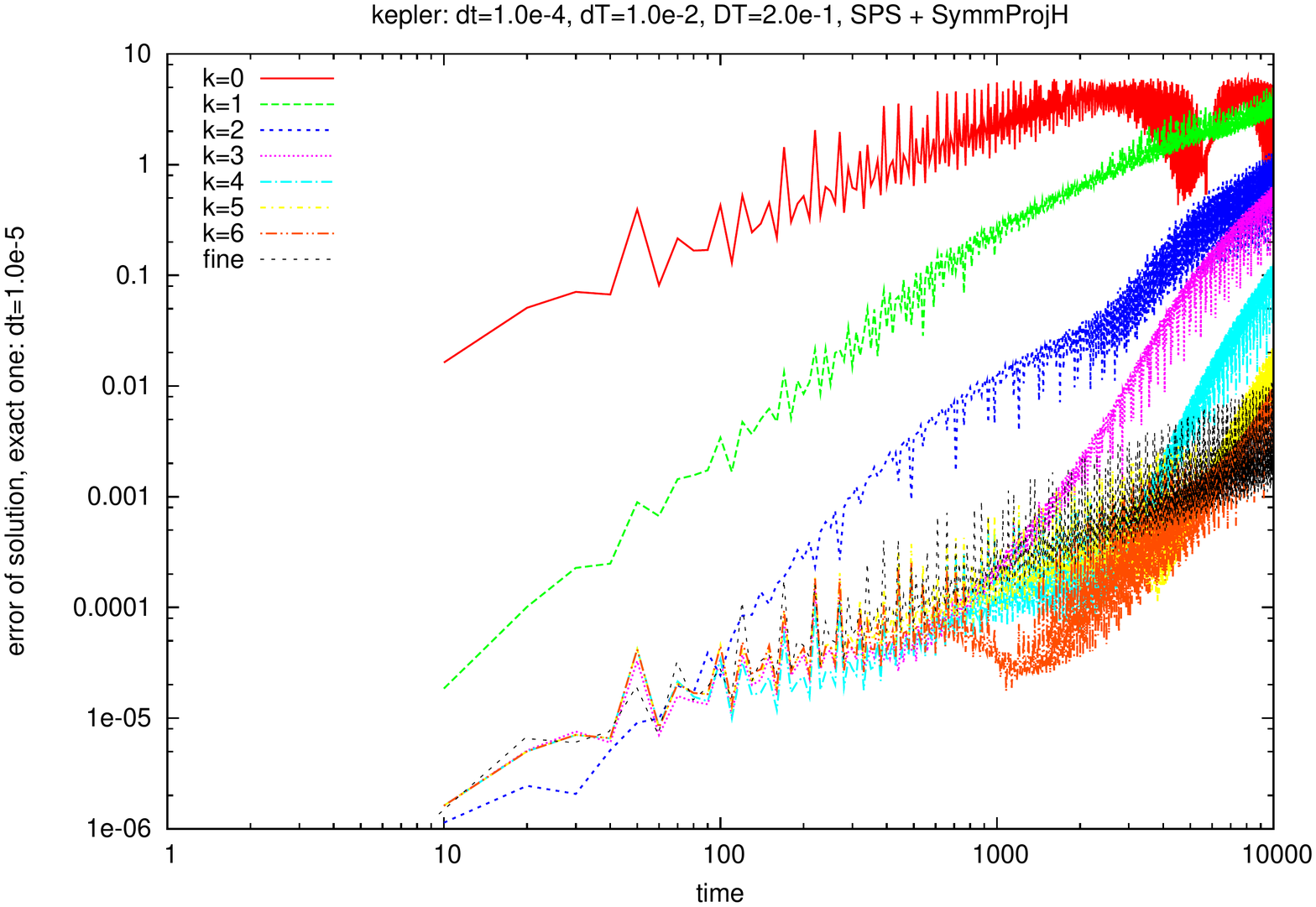}
\includegraphics[width=9cm,angle=0,origin=0]{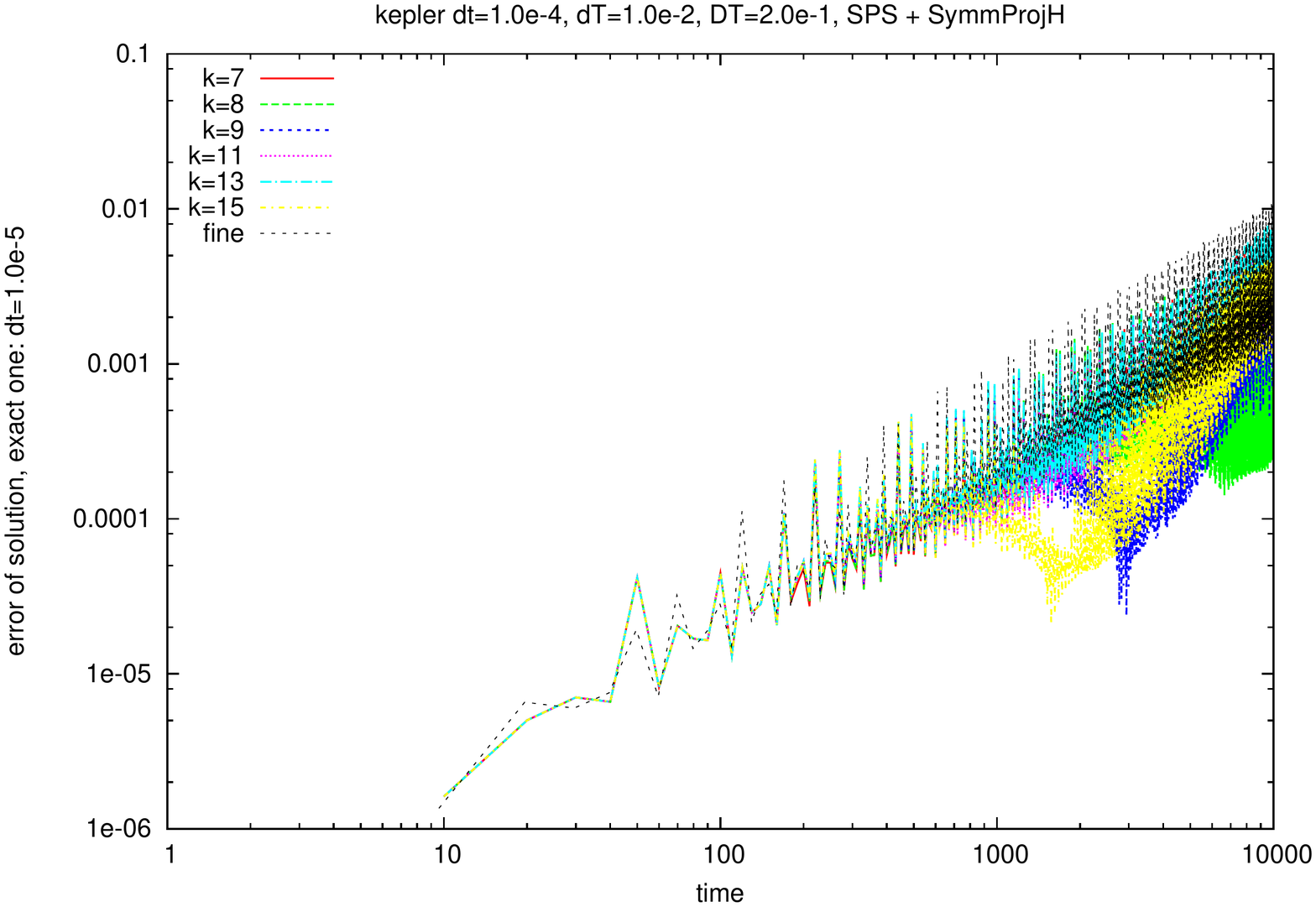}
\caption{ \label{fig:kepler-pos-symsym}
Errors on the trajectory for Kepler problem obtained by
Algorithm~\ref{algorithm-sps-symm-proj} ($\delta t = 10^{-4}$, $dT =
0.01$, $\Delta T = 0.2$).}  
\end{center}
\end{figure}

\begin{figure}[htbp]
\begin{center}
\includegraphics[width=9cm,angle=0,origin=0]{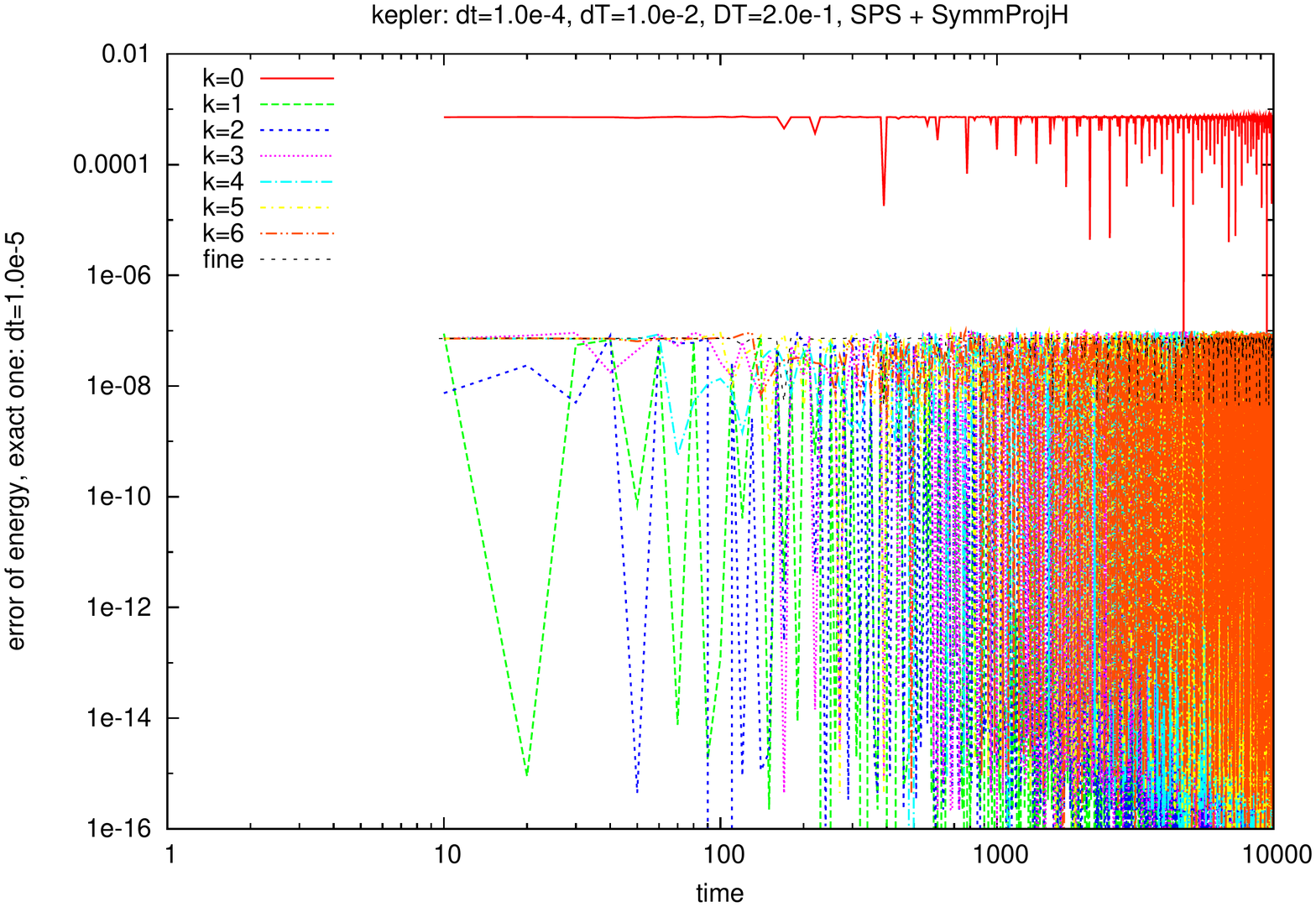}
\includegraphics[width=9cm,angle=0,origin=0]{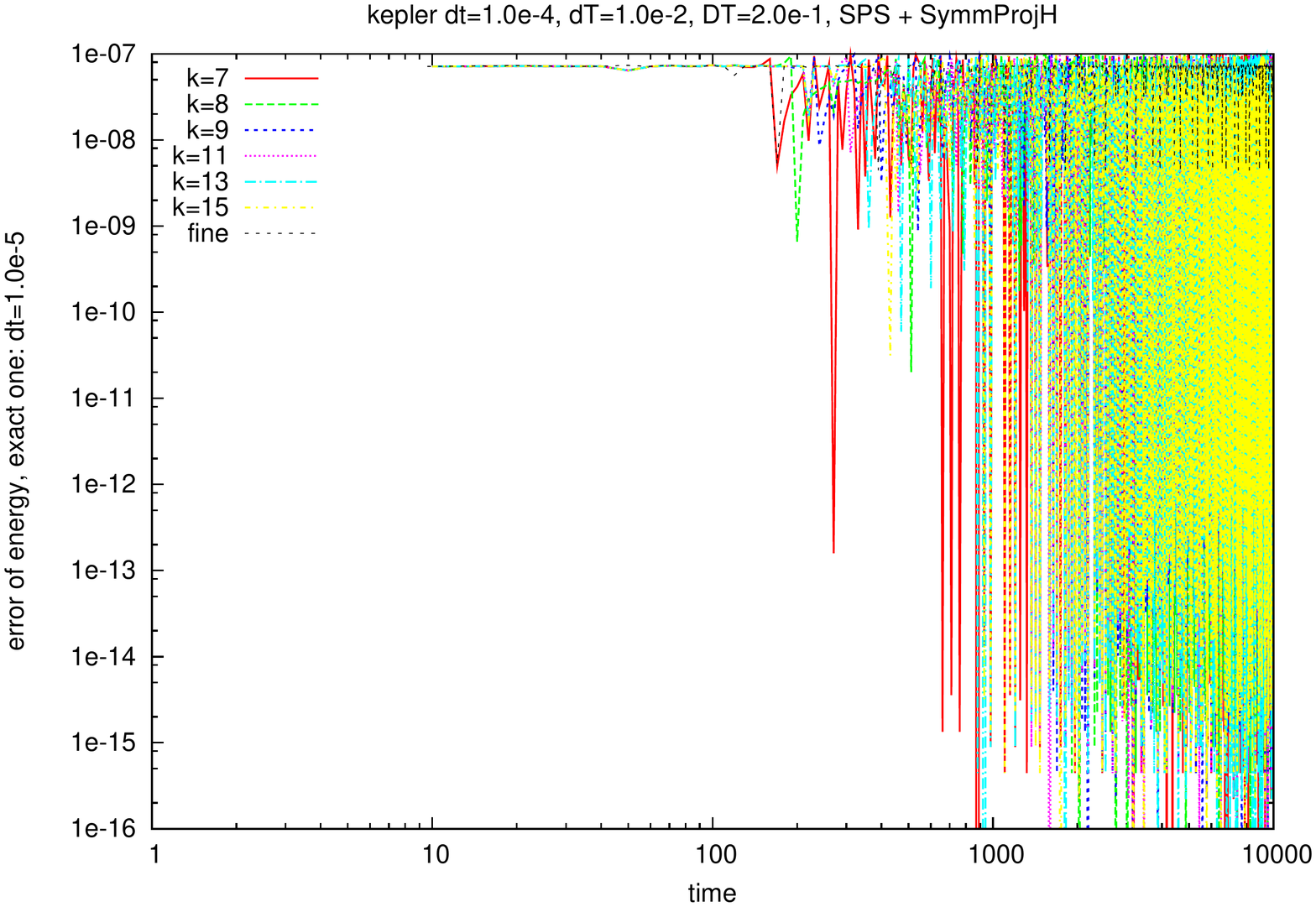}
\caption{ \label{fig:kepler-ener-symsym}
Errors on the energy for Kepler problem obtained by
Algorithm~\ref{algorithm-sps-symm-proj} 
($\delta t = 10^{-4}$, $dT = 0.01$, $\Delta T = 0.2$).} 
\end{center}
\end{figure}

\begin{figure}[htbp]
\begin{center}
\includegraphics[width=9cm,angle=0,origin=0]{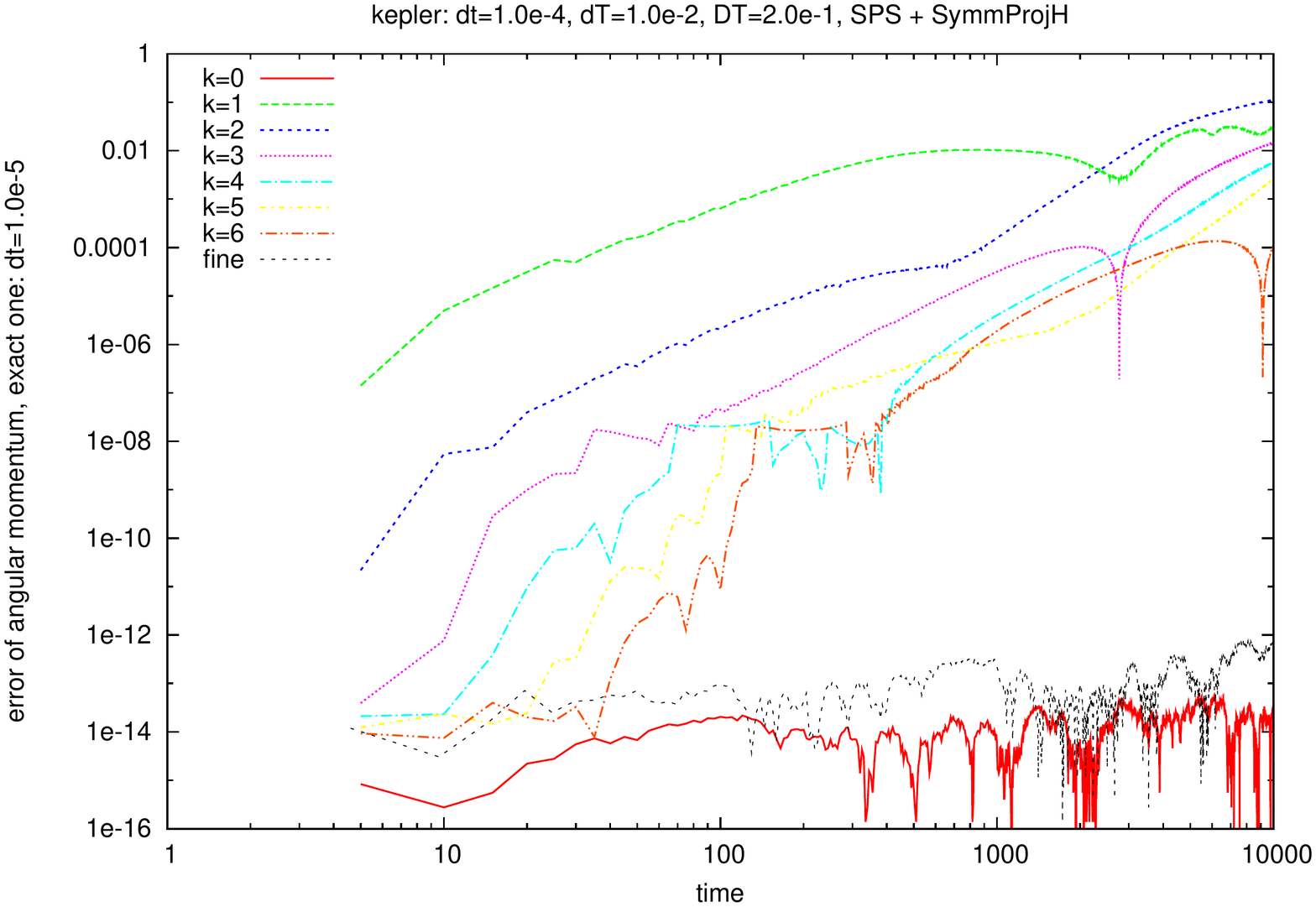}
\includegraphics[width=9cm,angle=0,origin=0]{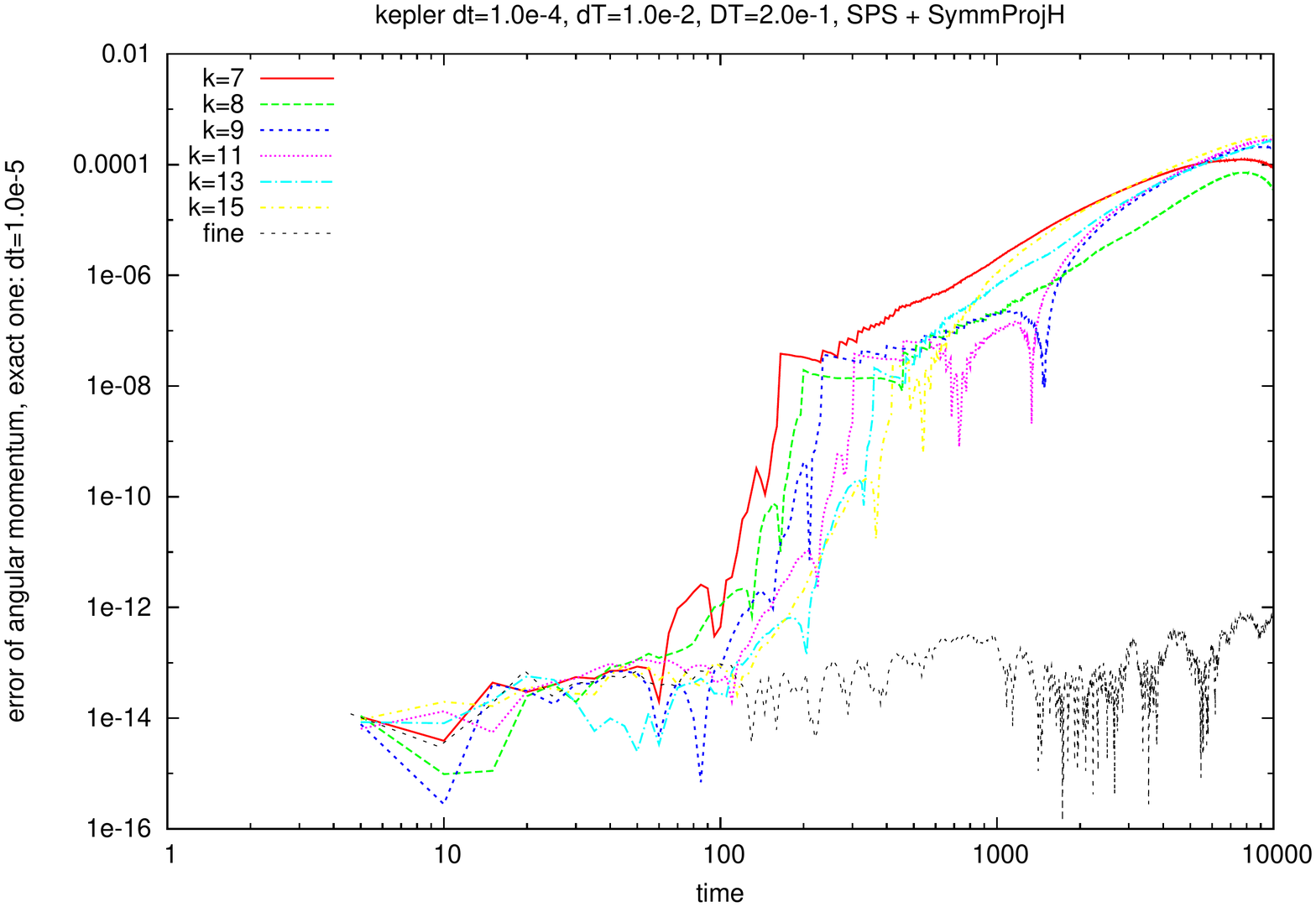}
\caption{ \label{fig:kepler-L-symsym}
Errors on the angular momentum for Kepler problem obtained by
Algorithm~\ref{algorithm-sps-symm-proj} 
($\delta t = 10^{-4}$, $dT = 0.01$, $\Delta T = 0.2$).} 
\end{center}
\end{figure}

We first consider Algorithm~\ref{algorithm-sps-symm-proj}, namely the
symmetric projection algorithm, with a projection on the constant energy
manifold. The iterative projection procedure is stopped using the three
criteria described above, with the parameters ${\tt tol} = 10^{-7}$ and
$N_{\rm iter}^{\rm max}=2$. Results are shown on 
Figs.~\ref{fig:kepler-pos-symsym}, \ref{fig:kepler-ener-symsym}
and~\ref{fig:kepler-L-symsym}, for the errors~\eqref{eq:err_traj} on the
trajectory, the relative errors~\eqref{eq:err_H} on the energy
preservation and the relative errors~\eqref{eq:err_L} on the angular
momentum preservation, respectively. 

We observe that, after only 5 iterations, the errors on the trajectory are
comparable to the errors on the trajectory of the fine scheme used sequentially
on the whole interval $[0,10^4]$. As expected, the energy is extremely
well preserved at any parareal iteration $k \geq 1$, although, as above,
we limit the number of iterations to 
solve the nonlinear projection step. Otherwise stated, convergence in
this nonlinear procedure is reached for a very small number of
iterations. We also observe that, for all parareal iterations $k \geq 7$, the
angular momentum is preserved with a relative accuracy of $5 \times 10^{-4}$
over the complete time range.

These results are better than those reported in
Section~\ref{sec:num_proj}, for Algorithm~\ref{algorithm-ps-proj}
(namely the standard parareal algorithm coupled with a standard, 
{\em not symmetrized}, projection step). First, trajectory convergence
is obtained after 5 parareal 
iterations, rather than 11 (compare Figs.~\ref{fig:kepler-projH-pos}
and~\ref{fig:kepler-pos-symsym}). Second, numerical results illustrate the fact that the nonlinear projection step on the constant energy
manifold converges more rapidly. With our newly constructed 
Algorithm~\ref{algorithm-sps-symm-proj}, the relative error on the
energy preservation is smaller than the 
tolerance $10^{-7}$ for any parareal iteration $k \geq 1$, whereas it
is the case only for $k \geq 6$ for Algorithm~\ref{algorithm-ps-proj}
(compare Figs.~\ref{fig:kepler-projH-ener} and~\ref{fig:kepler-ener-symsym}). 
Similarly, the angular momentum is better preserved, at any parareal
iteration (compare Figs.~\ref{fig:kepler-projH-L}
and~\ref{fig:kepler-L-symsym}). 

\bigskip

\begin{figure}[htbp]
\begin{center}
\includegraphics[width=9cm,angle=0,origin=0]{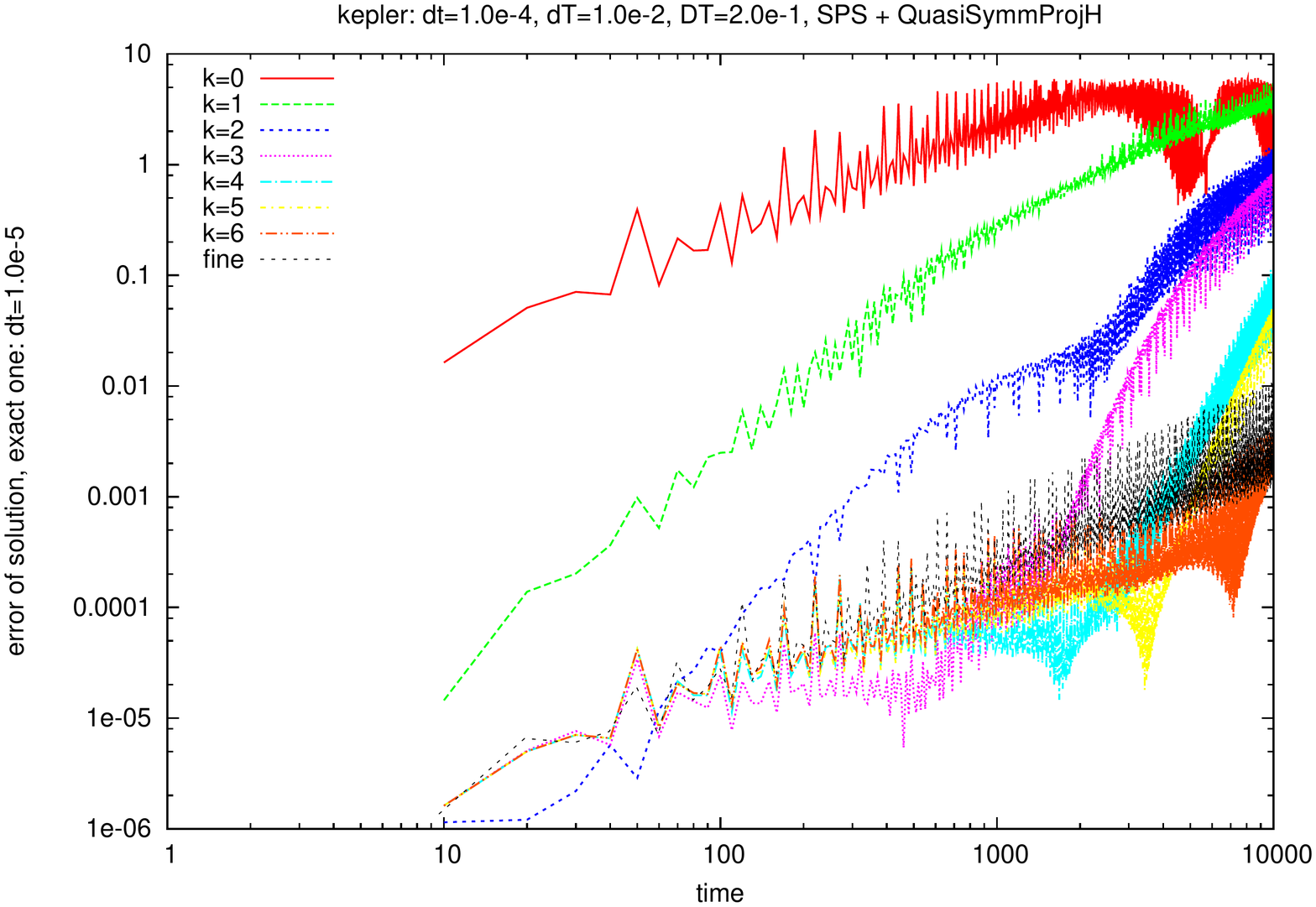}
\includegraphics[width=9cm,angle=0,origin=0]{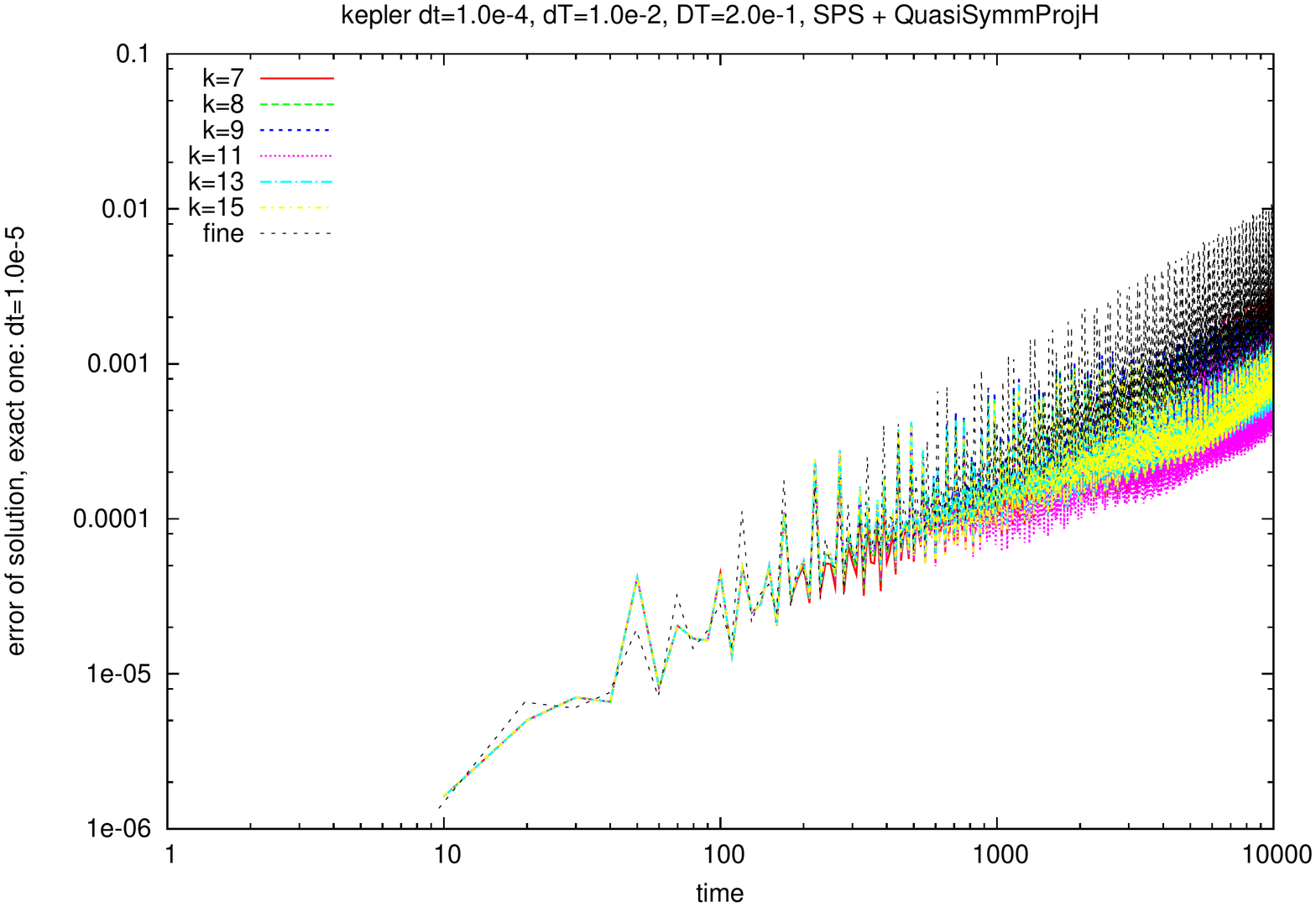}
\caption{ \label{fig:kepler-pos-symqsym}
Errors on the trajectory for Kepler problem obtained by the symmetric parareal
method with quasi-symmetric projection onto the constant energy manifold
($\delta t = 10^{-4}$, $dT = 0.01$, $\Delta T = 0.2$).} 
\end{center}
\end{figure}

\begin{figure}[htbp]
\begin{center}
\includegraphics[width=9cm,angle=0,origin=0]{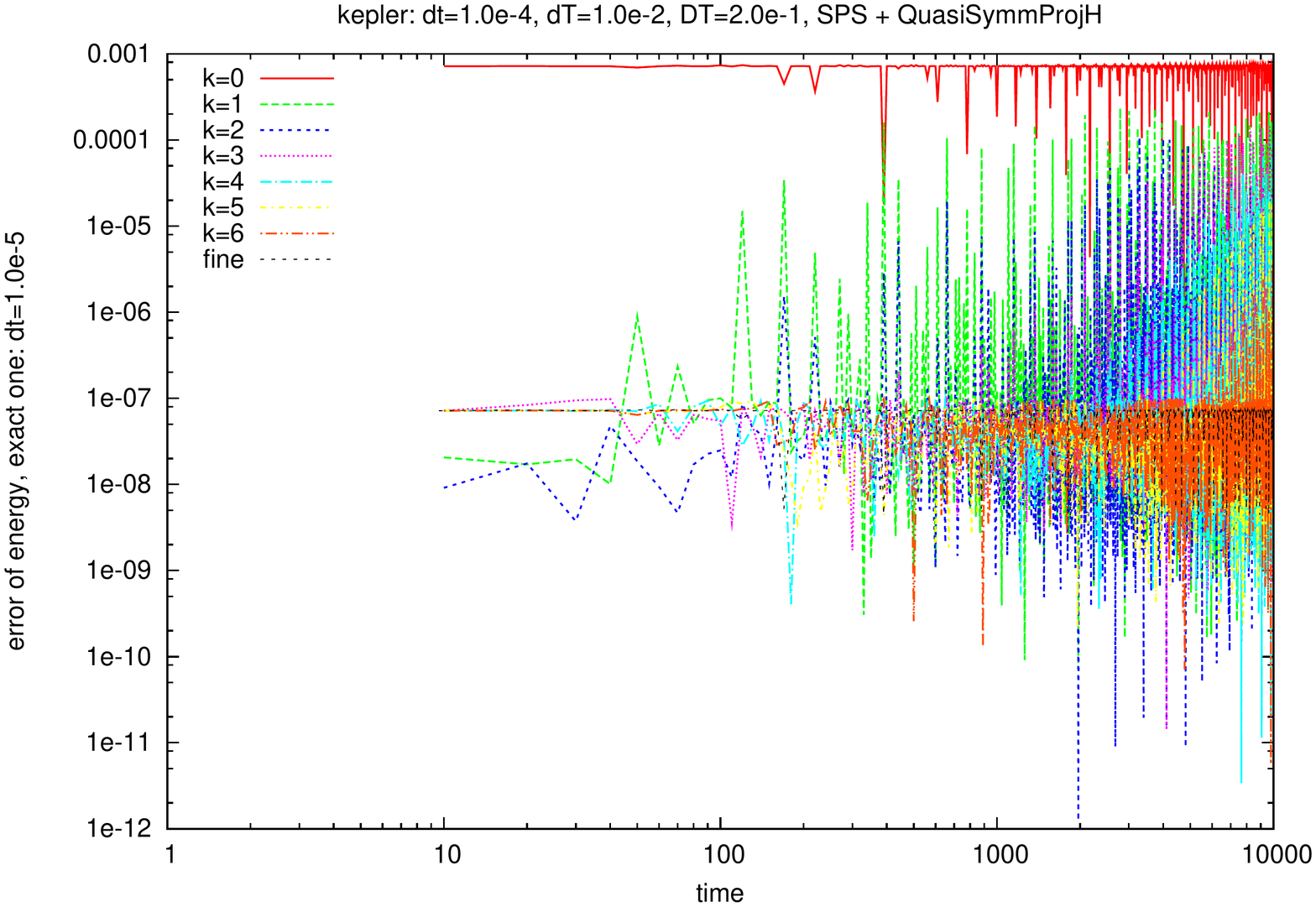}
\includegraphics[width=9cm,angle=0,origin=0]{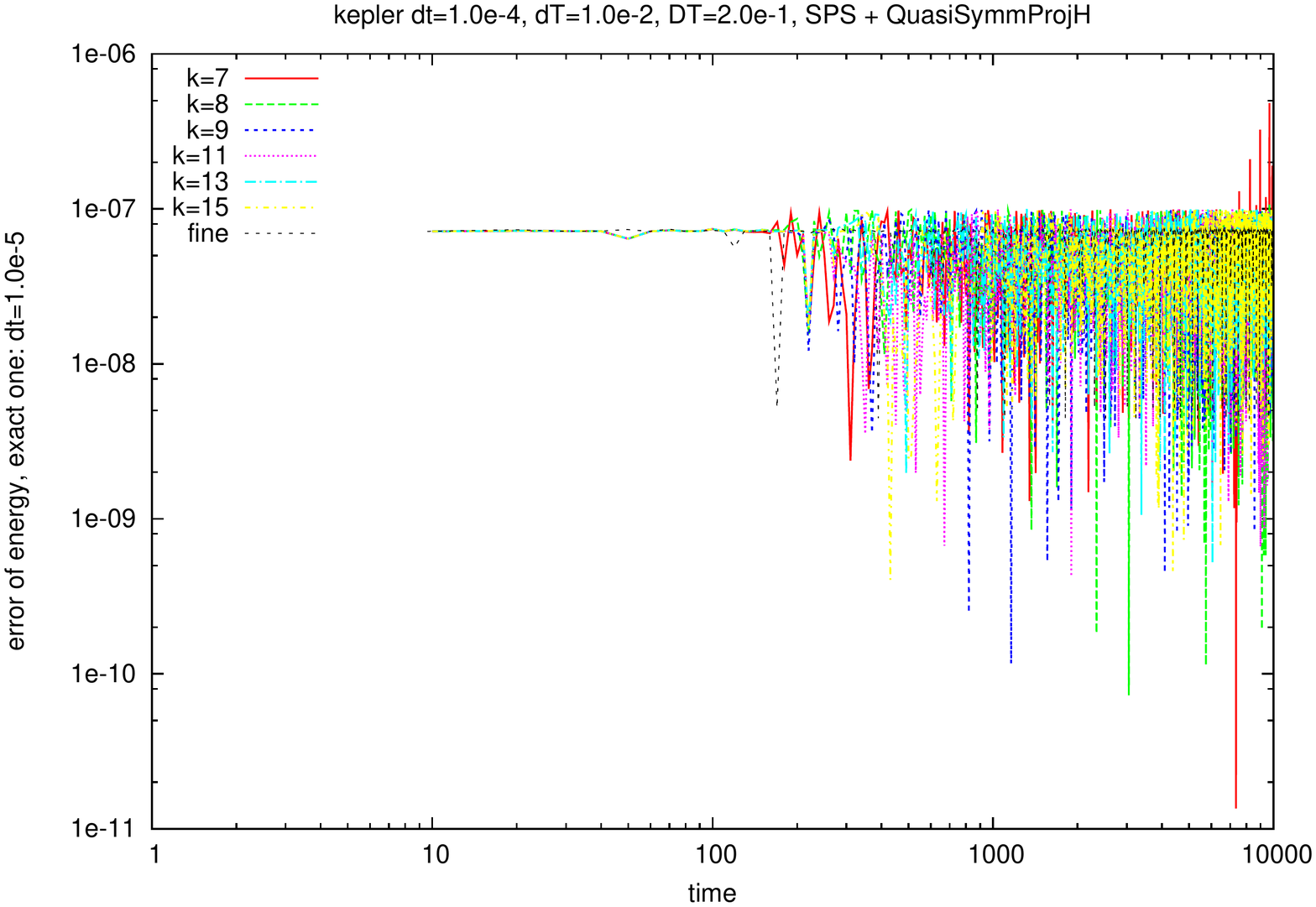}
\caption{ \label{fig:kepler-ener-symqsym}
Errors on the energy for Kepler problem obtained by the symmetric parareal
method with quasi-symmetric projection onto the constant energy manifold
($\delta t = 10^{-4}$, $dT = 0.01$, $\Delta T = 0.2$).}
\end{center}
\end{figure}

\begin{figure}[htbp]
\begin{center}
\includegraphics[width=9cm,angle=0,origin=0]{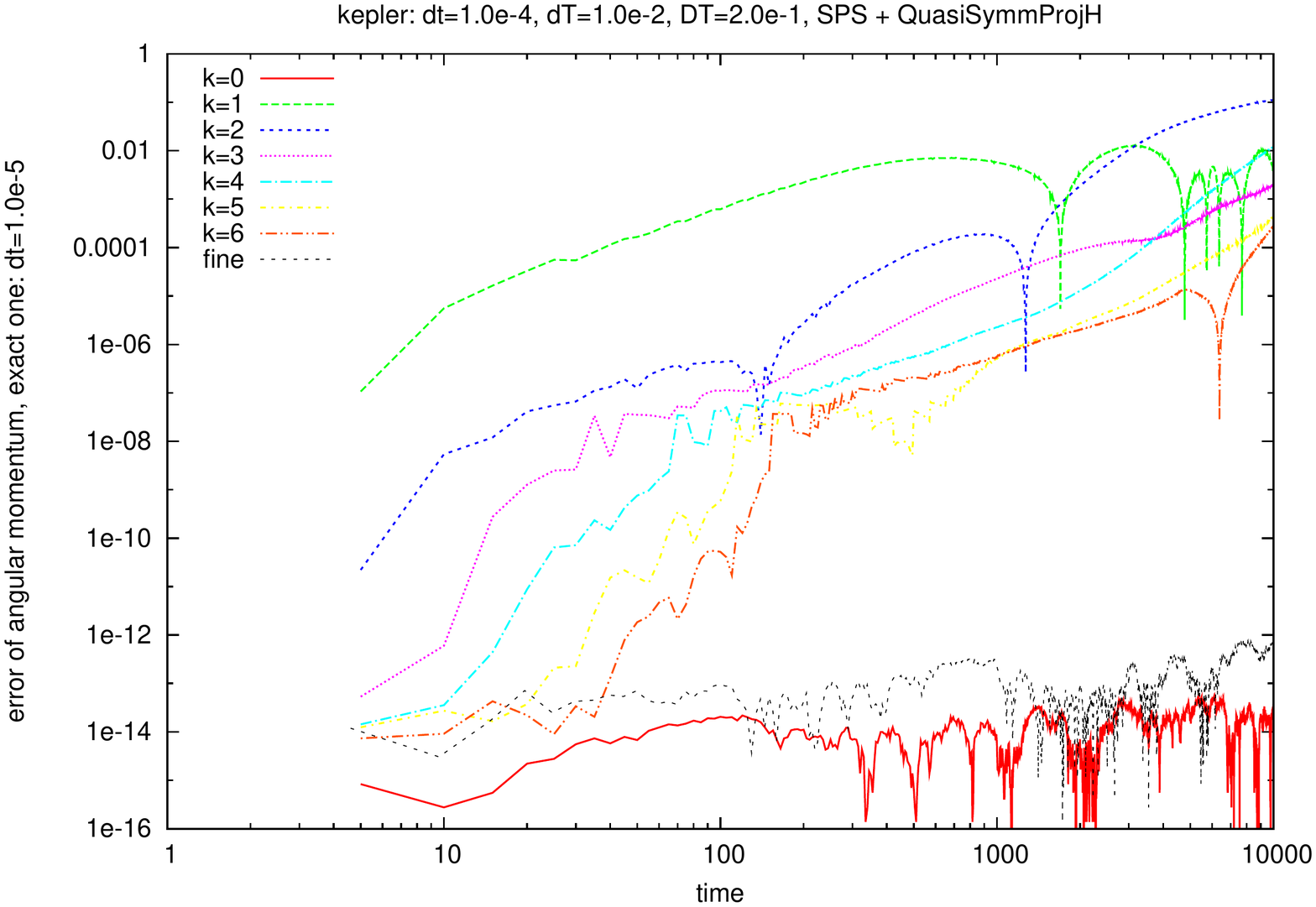}
\includegraphics[width=9cm,angle=0,origin=0]{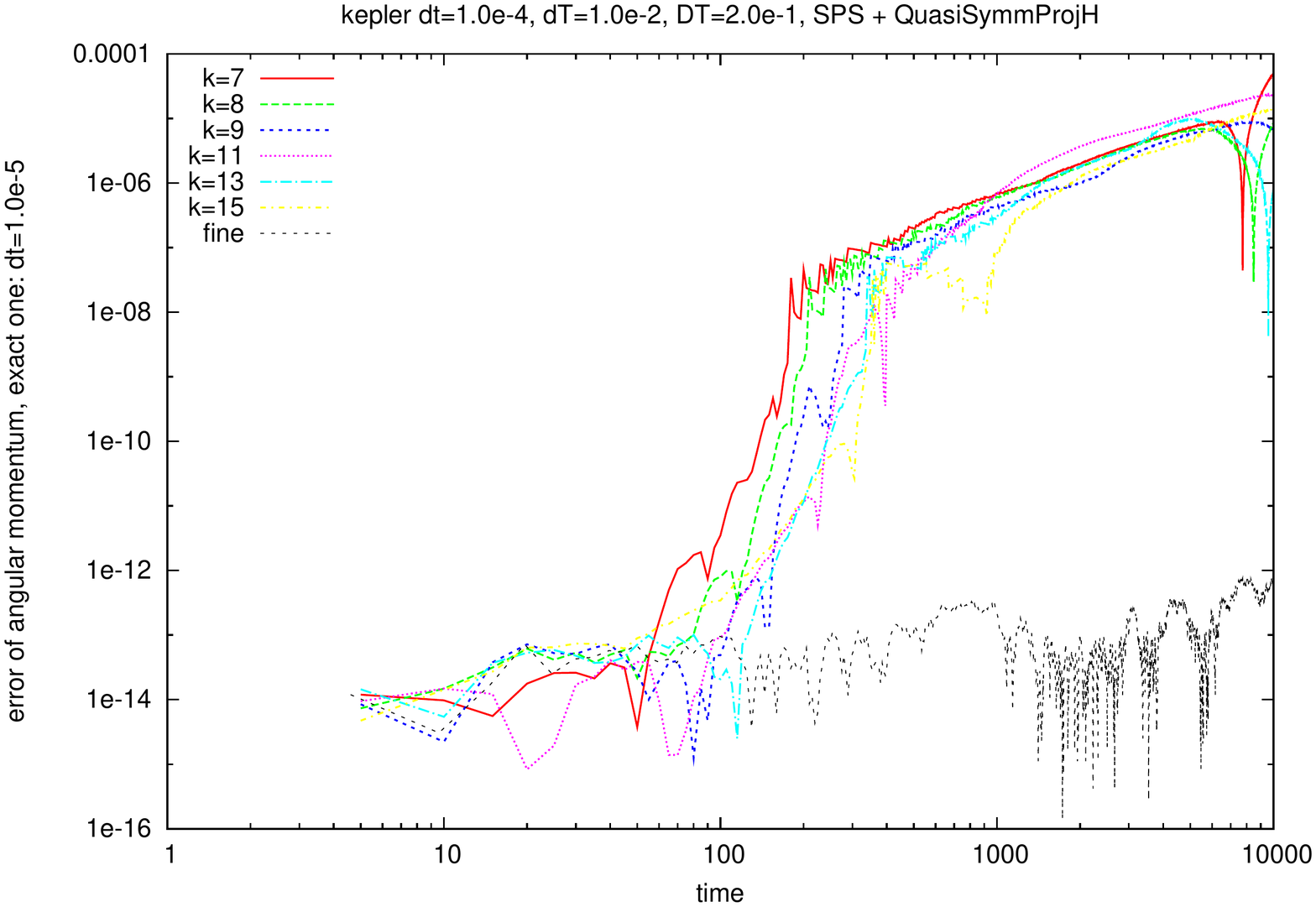}
\caption{ \label{fig:kepler-L-symqsym}
Errors on the angular momentum  for Kepler problem obtained by the symmetric parareal
method with quasi-symmetric projection onto the constant energy manifold
($\delta t = 10^{-4}$, $dT = 0.01$, $\Delta T = 0.2$).}
\end{center}
\end{figure}

We next consider the quasi-symmetric projection
algorithm, with projection on the same manifold, with the same
parameters. Results are shown on 
Figs.~\ref{fig:kepler-pos-symqsym}, \ref{fig:kepler-ener-symqsym}
and~\ref{fig:kepler-L-symqsym}, for the
errors~\eqref{eq:err_traj} on the trajectory, the relative
errors~\eqref{eq:err_H} on the energy preservation and the relative
errors~\eqref{eq:err_L} on the angular momentum preservation,
respectively.

We observe comparable performances in terms of trajectory accuracy and
angular momentum preservation as with
Algorithm~\ref{algorithm-sps-symm-proj} (which uses a {\em symmetric},
rather than quasi-symmetric, projection). However, the preservation of
energy is not as good: the error on the energy is smaller than the
tolerance requested by criteria C1 on the complete time range $[0,T]$
only for parareal iterations $k \geq 8$, whereas it is the case for any
$k \geq 1$ when the symmetric projection is employed (compare
Figs.~\ref{fig:kepler-ener-symsym} and~\ref{fig:kepler-ener-symqsym}).  

\section{An example in higher dimension: the outer solar system}
\label{sec:outer}

We now consider a system in a dimension higher
than those considered previously. We study here the evolution of the
outer solar system, which is composed of 6 three-dimensional particles,
representing the Sun and the planets from Jupiter to Pluto. The
Hamiltonian reads
\begin{equation}
\label{eq:ham_outer_complet}
H(q,p) = \frac12 p^T M^{-1} p + V(q), \quad q \in \R^{18}, \quad p \in \R^{18}, 
\end{equation}
with 
\begin{equation}
\label{eq:V_outer_complet}
V(q) = -\sum_{1 \leq i < j \leq 6} \frac{G \, m_i \, m_j}{\| q_i - q_j \|},
\end{equation}
and where $M$ is the diagonal mass matrix $M = \text{diag}(m_i)$. 
The values for the initial conditions, the gravitational constant $G$
and the masses $m_i$ of the particles are taken
from~\cite[Sec. I.2.4]{hlw}. By convention, the first particle is the
Sun, the mass of which is 1000 times as large as the mass of the heaviest planet
we consider. 

In the sections above, we pointed out that
Algorithm~\ref{algorithm-sps-symm-proj} yields good
results for the Kepler problem, and converges in few iterations to the
results obtained using the fine propagator in a purely sequential manner. We show here 
that (i) this algorithm also behaves well
to simulate the outer solar system, whose dimensionality is higher, and
(ii) that the coarse solver can be driven by a simplified
dynamics without damaging the good properties of the algorithm.  The
advantage of choosing a simplified coarse integrator is that the
complexity of the sequential computations using the coarse propagator is now negligible with
respect to the complexity of the parallel computations using the fine
integrator on the interval $[n \Delta T, (n+1) \Delta T]$, in
opposition to what was the case for the Kepler 
system.

The dynamics governed by the
Hamiltonian~\eqref{eq:ham_outer_complet}-\eqref{eq:V_outer_complet} can indeed
be considered as a perturbation of the dynamics driven by
\begin{equation}
\label{eq:ham_outer_simp}
H(q,p) = \frac12 p^T M^{-1} p + V_{\rm simp}(q), \quad q \in \R^{18}, \quad p \in \R^{18}, 
\end{equation}
with 
\begin{equation}
\label{eq:V_outer_simp}
V_{\rm simp}(q) = -\sum_{2 \leq j \leq 6} \frac{G \, m_1 \, m_j}{\| q_1 - q_j \|}.
\end{equation}
In $V_{\rm simp}$, we only take into account the gravitational
interaction between the Sun (at position $q_1$) and the other planets
(at position $q_j$, $2 \leq j \leq 6$), and we ignore the
interaction between pairs of planets. The fact 
that~\eqref{eq:ham_outer_complet}-\eqref{eq:V_outer_complet}
is a small perturbation
to~\eqref{eq:ham_outer_simp}-\eqref{eq:V_outer_simp} owes to the huge
discrepancy between the mass $m_1$ of Sun and the masses of the
planets. The latter system is less expensive to simulate, since
$V_{\rm simp}$ is a sum of 5 terms, whereas $V$ is a sum of 15 terms.

We set $\delta t = 10^{-2}$, $dT = 50$ and $\Delta T = 200$, and study
the algorithms over the time range $[0,T]$ with $T=2 \times 10^5$. 

\begin{figure}[htbp]
\begin{center}
\includegraphics[width=9cm,angle=0]{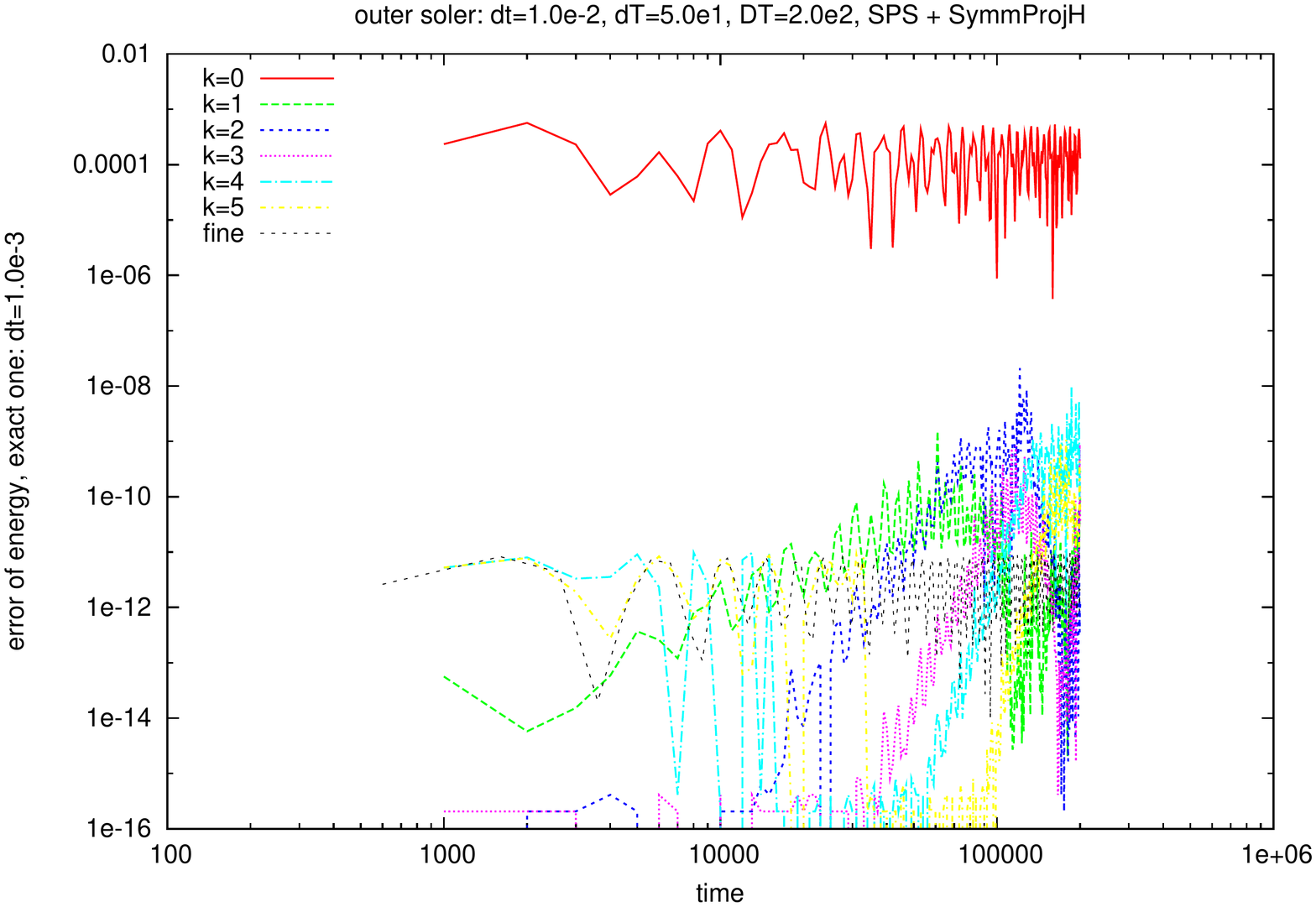}
\includegraphics[width=9cm,angle=0]{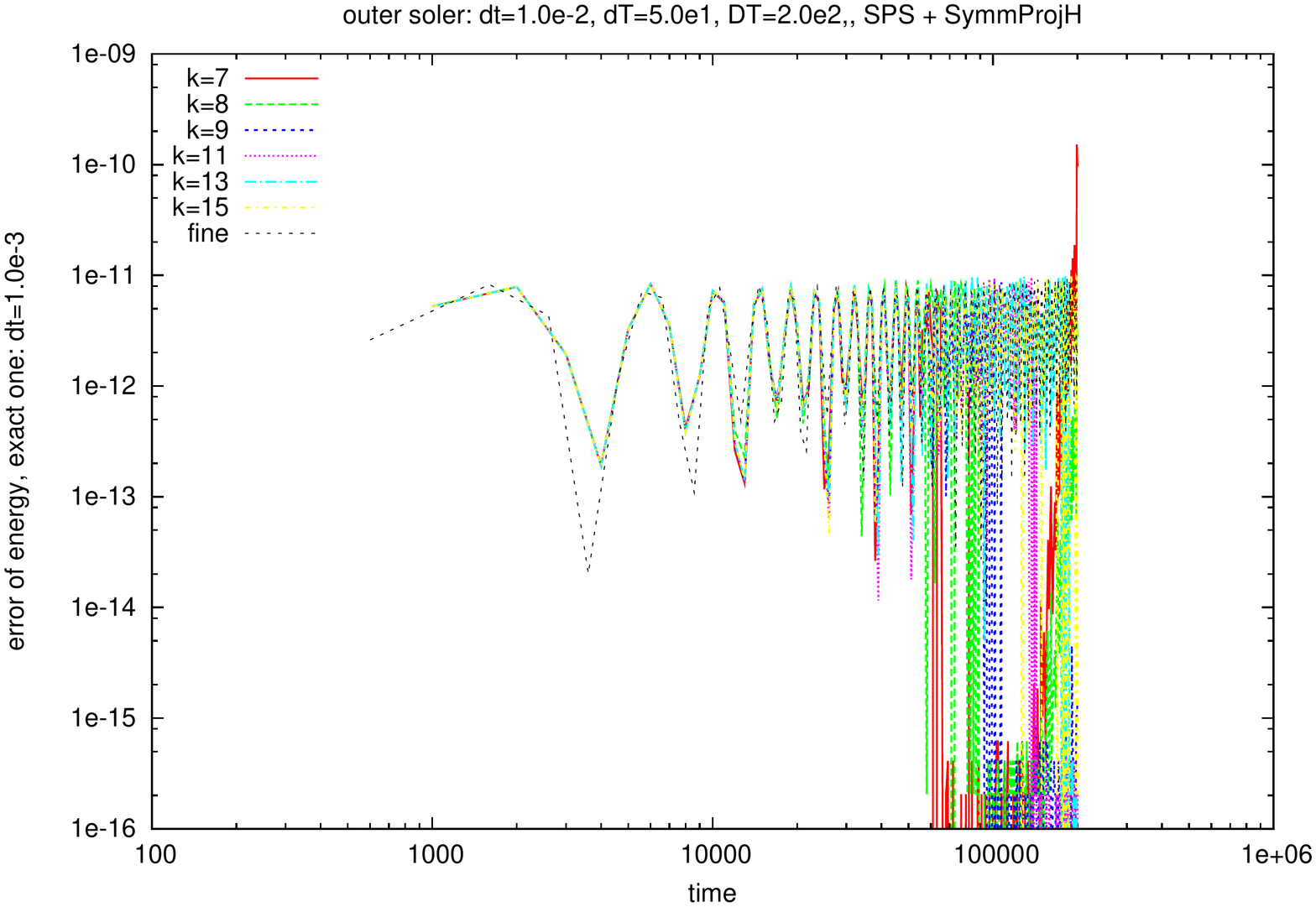}
\caption{ \label{fig-solar-GCS-case2-SPS-SymmProjH-E-0-5}
Errors~\eqref{eq:err_H} on the energy for the outer solar system
obtained using Algorithm~\ref{algorithm-sps-symm-proj} ($\delta
t=10^{-2}$, $dT=50$, $\Delta T = 200$).} 
\end{center}
\end{figure}

\begin{figure}[htbp]
\begin{center}
\includegraphics[width=9cm,angle=0]{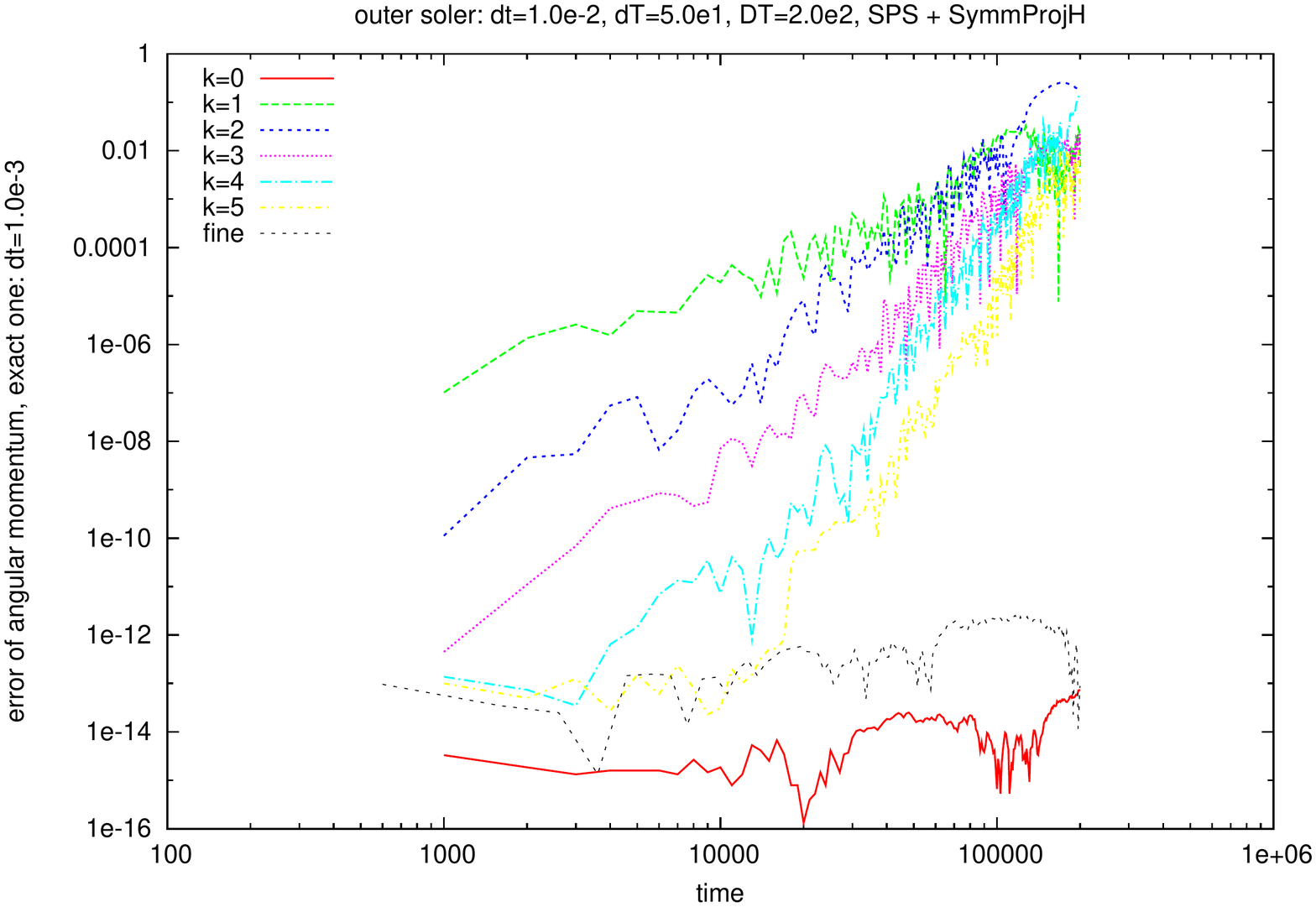}
\includegraphics[width=9cm,angle=0]{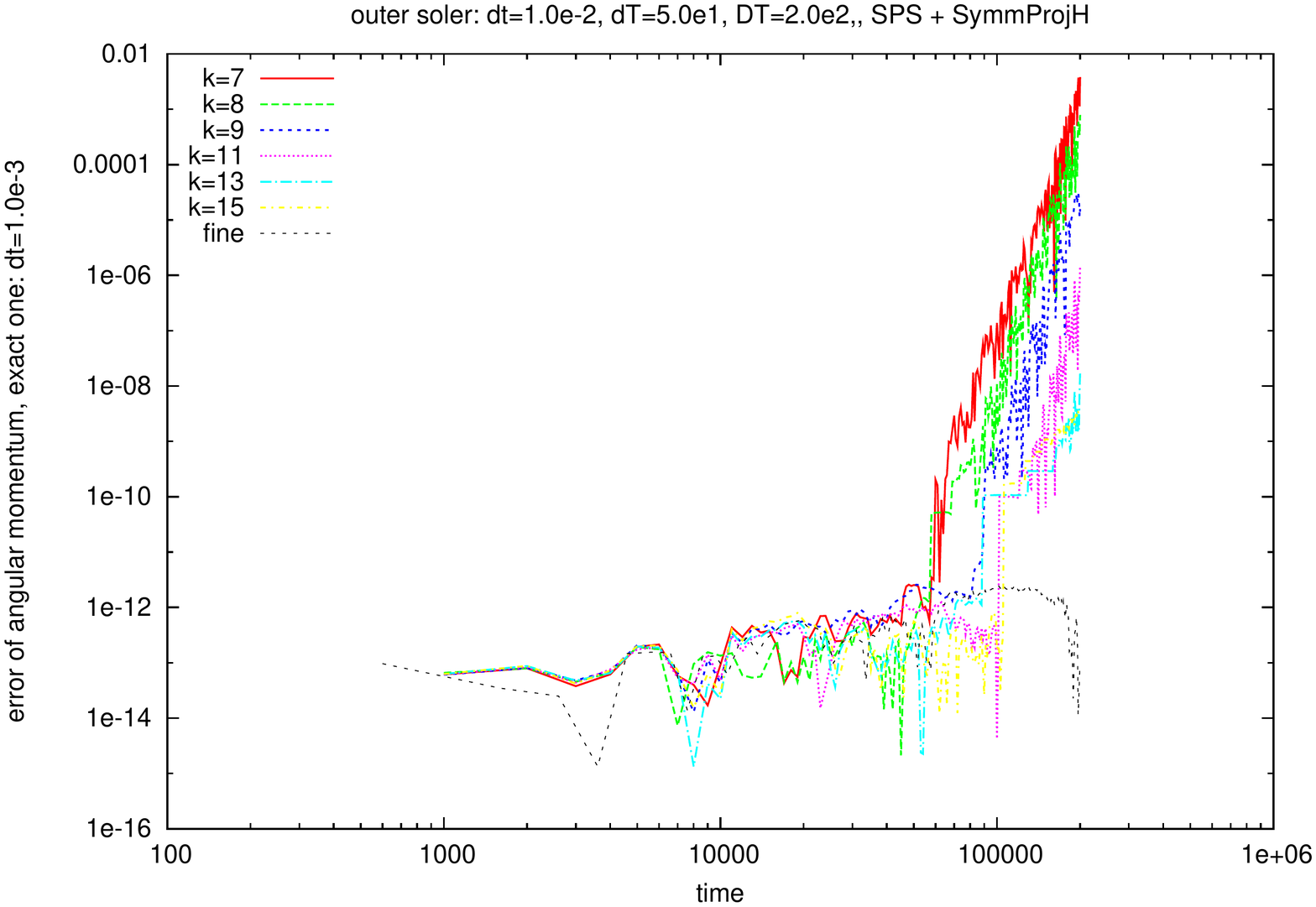}
\caption{ \label{fig-solar-GCS-case2-SPS-SymmProjH-L0-0-5}
Errors~\eqref{eq:err_L} on the angular momentum first component for the
outer solar system obtained using
Algorithm~\ref{algorithm-sps-symm-proj} ($\delta t=10^{-2}$, $dT=50$,
$\Delta T = 200$).} 
\end{center}
\end{figure}

\begin{figure}[htbp]
\begin{center}
\includegraphics[width=9cm,angle=0]{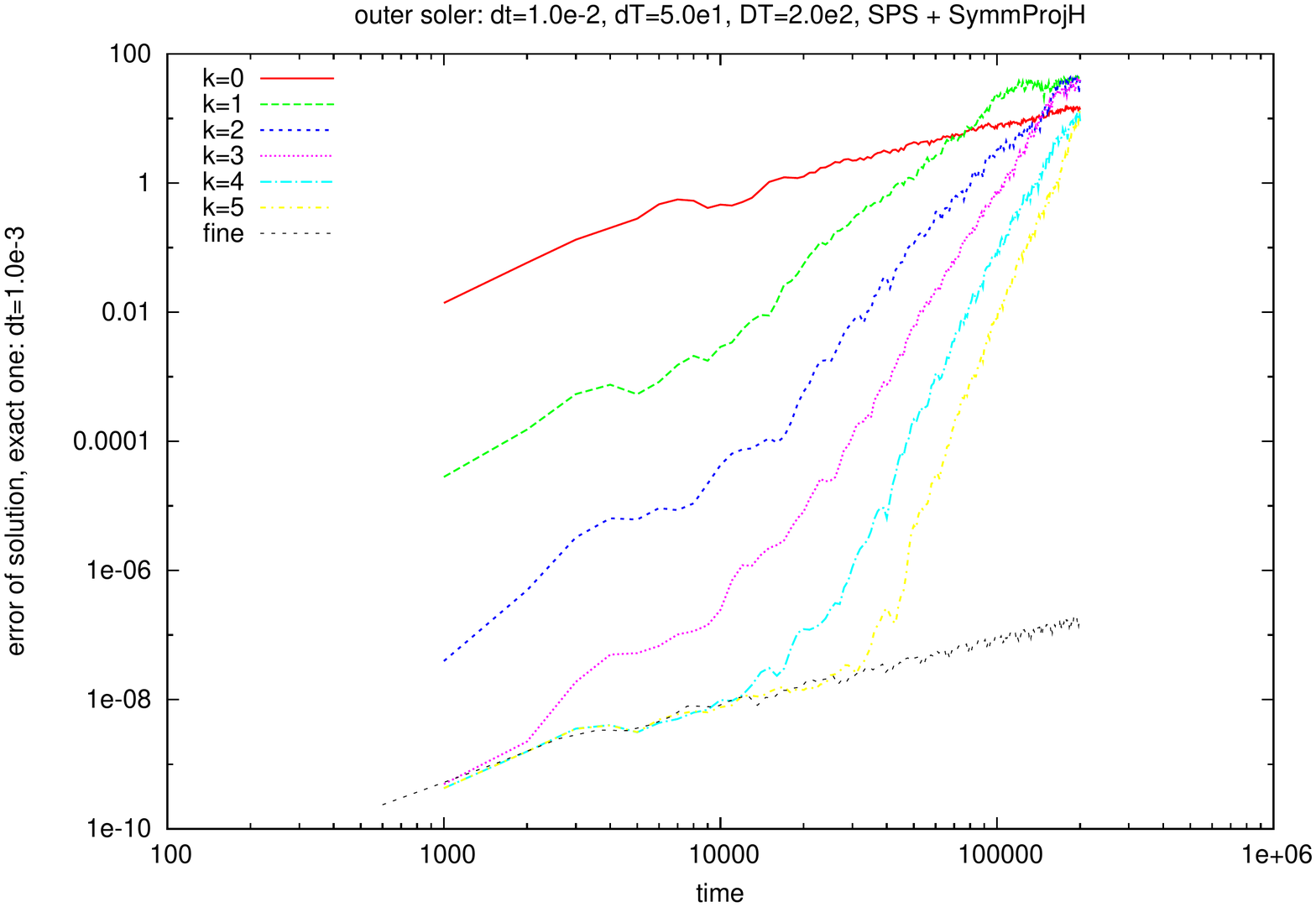}
\includegraphics[width=9cm,angle=0]{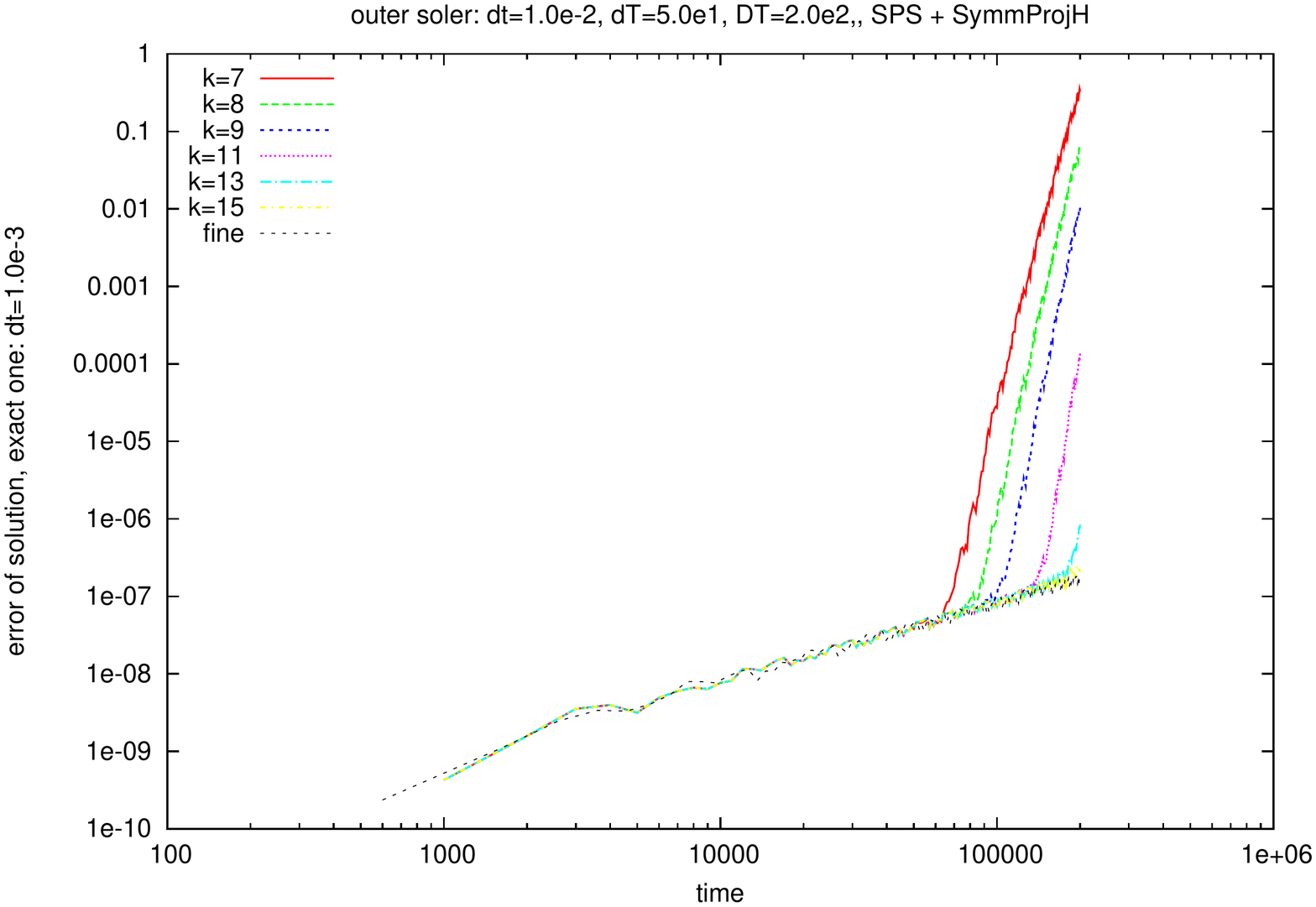}
\caption{ \label{fig-solar-GCS-case2-SPS-SymmProjH-P-0-5}
Errors~\eqref{eq:err_traj} on the trajectory for the outer solar system
obtained using Algorithm~\ref{algorithm-sps-symm-proj} ($\delta
t=10^{-2}$, $dT=50$, $\Delta T = 200$).} 
\end{center}
\end{figure}

Results obtained with Algorithm~\ref{algorithm-sps-symm-proj} are shown on
Figs.~\ref{fig-solar-GCS-case2-SPS-SymmProjH-E-0-5},
\ref{fig-solar-GCS-case2-SPS-SymmProjH-L0-0-5}
and~\ref{fig-solar-GCS-case2-SPS-SymmProjH-P-0-5}. As above, we project
on the manifold ${\cal M}$ of constant energy, and 
we use the three criteria C1, C2, C3 to stop the nonlinear projection step,
with the parameters $N_{\rm iter}^{\rm max} = 2$ and ${\tt tol} = 10^{-11}$
(which again corresponds to the error on the energy preservation of the fine
scheme).

\begin{figure}[htbp]
\begin{center}
\includegraphics[width=6cm,angle=0]{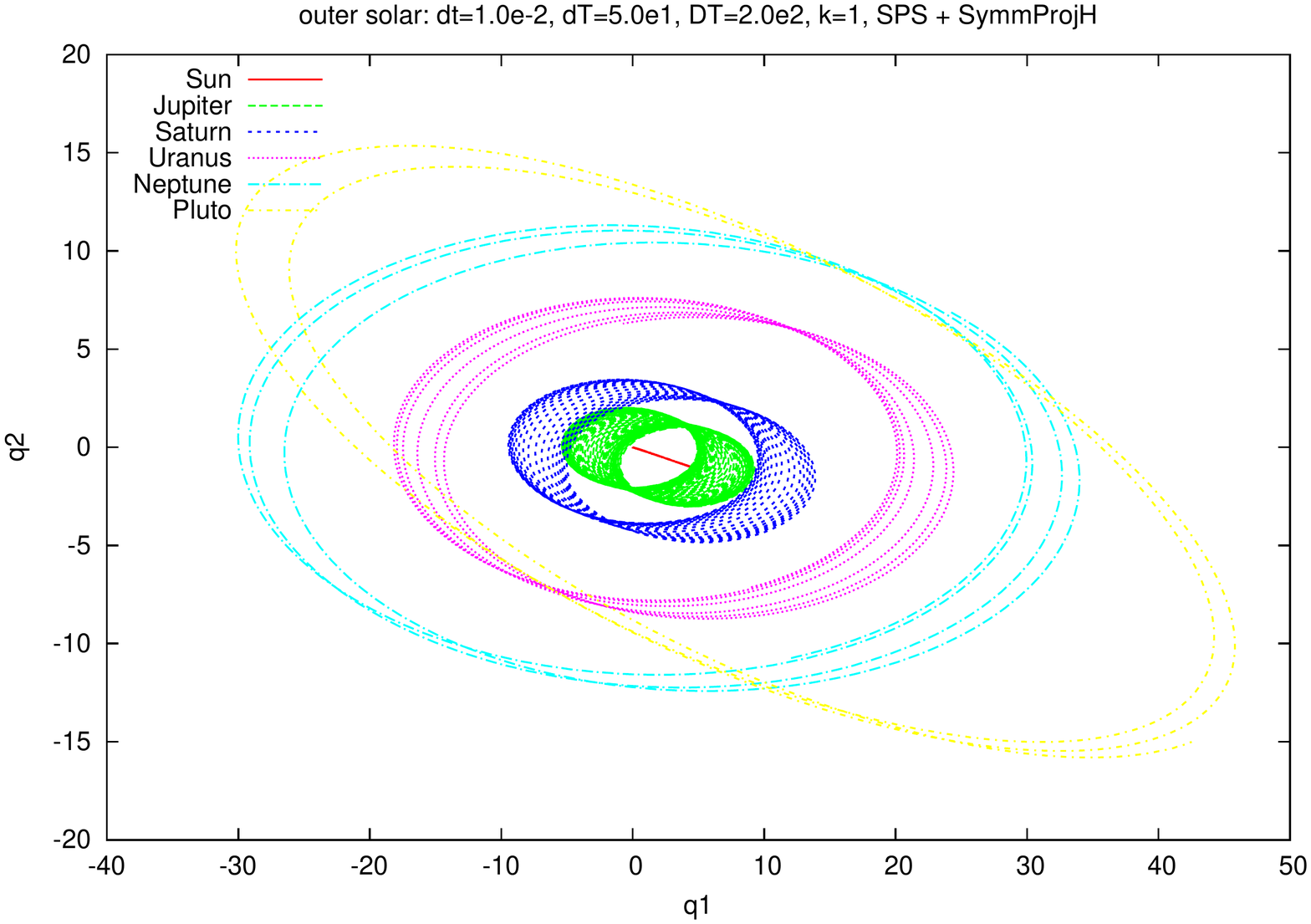}
\includegraphics[width=6cm,angle=0]{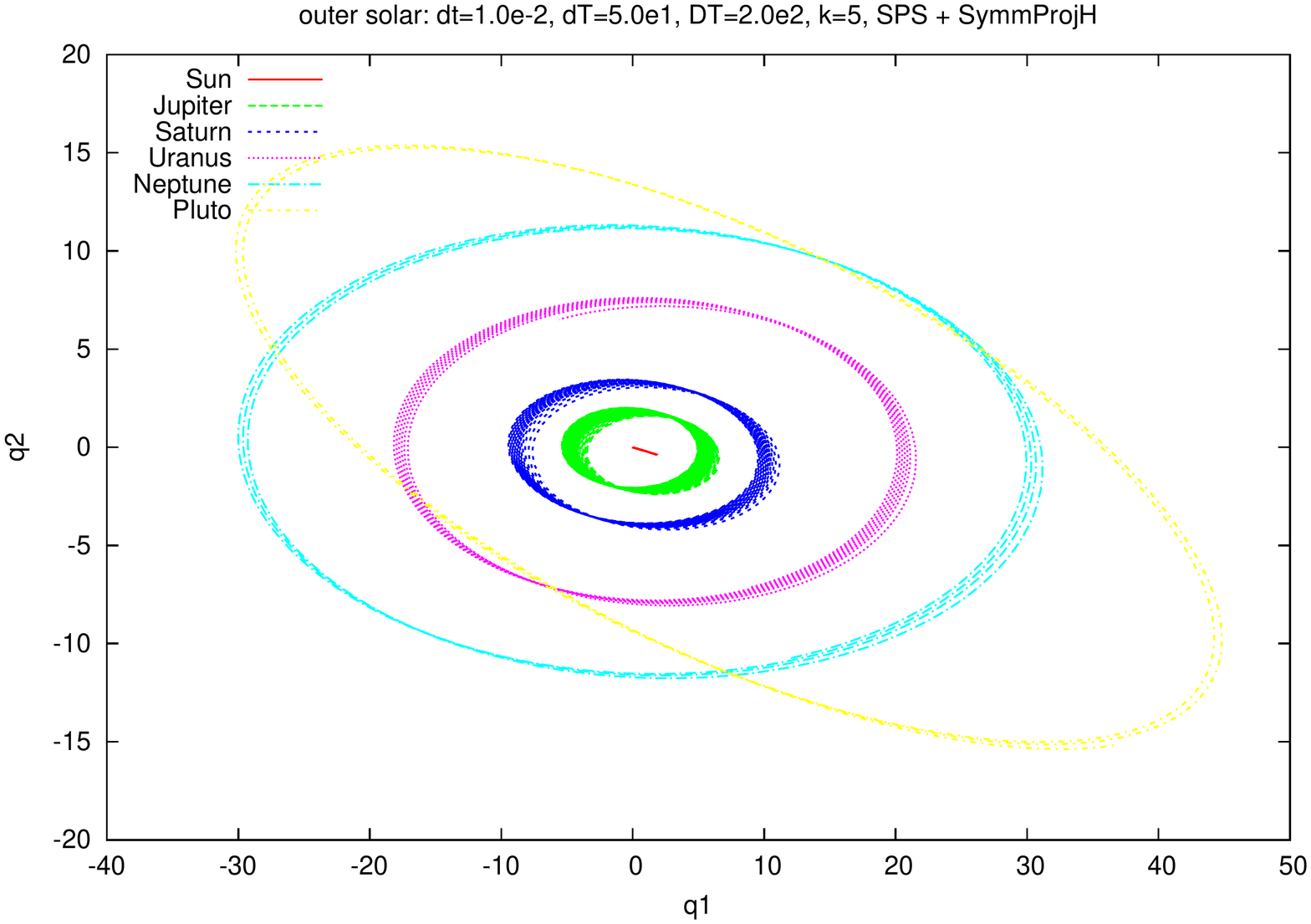}
\includegraphics[width=6cm,angle=0]{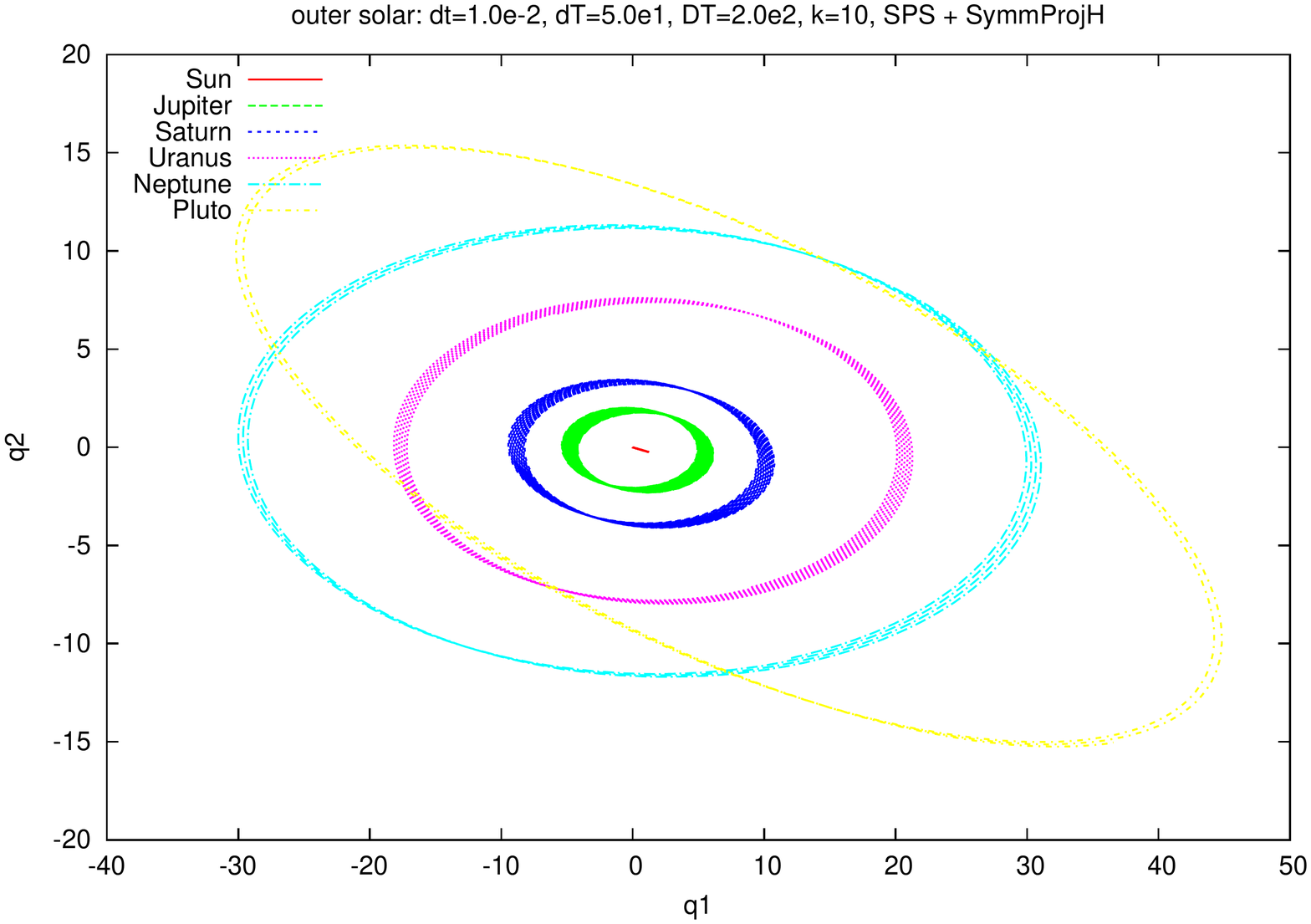}
\includegraphics[width=6cm,angle=0]{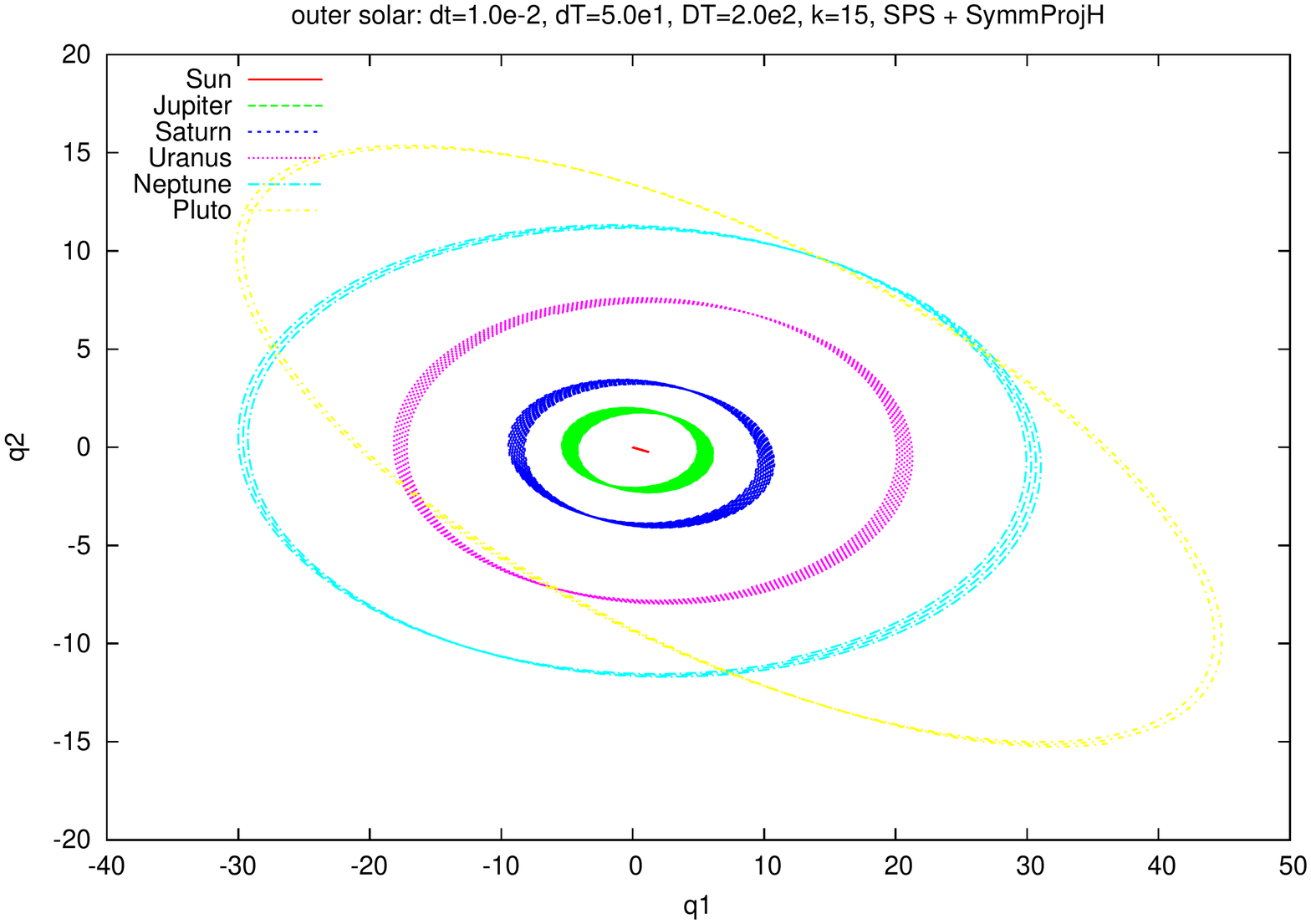}
\caption{ \label{fig-solar-GCS-case2-SPS-SymmProjH-q1-q2-5-xz}
Trajectory of the outer solar system obtained using
Algorithm~\ref{algorithm-sps-symm-proj} ($\delta t=10^{-2}$, $dT=50$,
$\Delta T = 200$), at iteration $k=1$ (top left), $k=5$ (top right),
$k=10$ (bottom left), $k=15$ (bottom right).} 
 \end{center}
\end{figure}

For $k \geq 8$, the energy is preserved up to the tolerance prescribed
by criteria C1 on the whole time range. The first convergence criteria that is
fulfilled is thus C1 (hence convergence in the
projection step is achieved with less than $N_{\rm iter}^{\rm max} = 2$
iterations). Note that a larger value of $N_{\rm iter}^{\rm max}$ would
be required to reach convergence for $k \leq 7$. We also obtain a 
good preservation of the angular momentum, when $k \geq 5$ (the error is
then always smaller than 1~\%). At $k=15$, the trajectory
error is comparable to the error made by the fine
propagator. Note that the error on the trajectory is already small
when $k \geq 9$ (the absolute error~\eqref{eq:err_traj} is smaller than
0.01, and the relative error is even smaller, as $\| q \|$ is of order 1
or larger, see Fig.~\ref{fig-solar-GCS-case2-SPS-SymmProjH-q1-q2-5-xz}).

\begin{remark}
We have checked that, if we set $N_{\rm iter}^{\rm max} = 5$, we obtain a
better energy preservation (the error on the energy is then smaller than
the tolerance prescribed by criteria C1 over the complete time range, at any
parareal iteration $k$), 
and comparable results for the angular momentum and the
trajectory. However, the speed-up is then a little smaller. 
\end{remark}

Algorithm~\ref{algorithm-sps-symm-proj} also provides very good
qualitative results in terms of 
trajectory. On Fig.~\ref{fig-solar-GCS-case2-SPS-SymmProjH-q1-q2-5-xz},
we plot the actual trajectories of the different planets obtained at
parareal iterations $k=1$, 5, 10 and 15. Except for $k=1$, trajectories
are qualitatively similar to the ones obtained with a symplectic or
symmetric integration of the solar system. We also note that, even
though the trajectory is quantitatively wrong for $k=5$, it is
qualitatively correct. This is a standard observation in geometric
integration. In turn, for $k \geq 10$, the trajectory is quantitatively
correct, as we observed from
Fig.~\ref{fig-solar-GCS-case2-SPS-SymmProjH-P-0-5}. 

\bigskip

We conclude this section by discussing the complexity of the algorithm. 
Note that the coarse solver integrates a simpler system than the fine
solver. Its complexity, for an equal time step, is hence smaller. We
are in the case of assumption (\ref{assumption}).
In addition, the average number of iterations in the nonlinear
projection problem
is 
$$
m_{\rm sym-proj} = 1.12. 
$$

With these choices, we need
$$
K_{\rm SP/sym-proj} = 15 \text{ parareal iterations}
$$
to reach convergence, namely to obtain an error on the trajectory as small as
that obtained with the fine solver used sequentially over the complete
time range (see
Fig.~\ref{fig-solar-GCS-case2-SPS-SymmProjH-P-0-5}). 

Using our algorithm, and assuming that we have 1000 processors at our
disposal, we thus obtain a speed-up of the order of $1000/15 \approx
66$: for a similar accuracy, the computational time is 66 times smaller.

\section{Conclusion}
\label{sec:conclusion}

The parareal algorithm allows to take benefit from a parallel
architecture to speed up the numerical integration of systems
of ODEs or time-dependent PDEs. It consists in the iterative use of a
coarse propagator and a fine one. 

In case the system is Hamiltonian, 
discretization schemes that preserve some geometrical features of the
underlying Hamiltonian dynamics are known to yield the best results, in
terms of preservation of energy, but also preservation of other first
invariants, trajectory accuracy, \ldots
These properties are not compatible with the plain parareal
algorithm, which has some issues when used for extremely long
simulations of Hamiltonian systems. 

We have identified in this article two ingredients that, when
combined, yield an algorithm that is better adapted to the Hamiltonian
context. The first ingredient consists in making 
the scheme {\em symmetric}. The idea is simple, but the adaptation to
the parareal context is not so simple, because the original parareal algorithm
does not naturally write as a one-step method. We have shown here how
to obtain a symmetric algorithm that retains the parallel features of the
original, plain, parareal algorithm. This is yet
not sufficient to provide an improvement on the geometric properties of
the resulting scheme and thus on the long time 
simulations. The second ingredient consists in {\em projecting} the solution on
the constant energy manifold. Here, this projection is not realized at
every time step but only at the coarse level. This may be the reason why
such a projection, which is generally not recommended as a viable
solution for a general solver, appears here to be the right complement
(at least in the form of the symmetric or quasi symmetric projector) to
the symmetrization of the parareal algorithm.  
Note also that the projection is inexpensive, as only one or two
iterations are needed to solve this nonlinear procedure.  

The conjunction of both ingredients above, in a case where the fine and
coarse solvers have good geometric properties, is shown here to be the
best choice. The symmetrized parareal algorithm has, as expected, good
geometric features, apart from the resonance problem. The
projection at the end of each propagation interval $[T_n, T_{n+1}]$
takes care of this resonance problem. After a limited
number of parareal iterations, we obtain a trajectory which is as
accurate as the one obtained using the sequential fine solver over the
whole time range $[0,T]$. In turn, this convergence stabilizes the
values of the invariants, besides the energy. 

This allows to achieve a substantial speed-up. Let us recall that the
full efficiency offered by the parallel architecture with the parareal
in time algorithm can very scarcely be obtained on systems of ODE's. In
order to get a full speed-up, the parareal algorithm (in its symmetric
version or not) has to be combined with other iterative procedures (such
as domain decomposition approaches for PDEs, see
e.g.~\cite{maday-turinici-04}).  

Finally, we emphasize that, in this article, in addition to energy
conservation, we have chosen to illustrate the 
quality of the integrators by the accuracy with which the
trajectories are computed. In many examples of Hamiltonian systems, for 
even longer simulations, the trajectories cannot be well approximated
and the interest is more in e.g. averages of some quantities along the
trajectories. The problem should then be considered from a different
viewpoint as the one adopted here.

\section*{Acknowledgements}

The work of Xiaoying Dai and Yvon Maday is supported in part by the
Agence Nationale de la Recherche, under the grant ANR-06-CIS6-007
(PITAC). The work of Claude Le Bris and Fr\'ed\'eric Legoll is supported
by INRIA under the grant ``Action de Recherche Collaborative'' HYBRID
and by the Agence Nationale de la Recherche, under the grant
ANR-09-BLAN-0216-01 (MEGAS).


\begin{thebibliography}{99}

\bibitem{rattle1} {\sc H.C. Andersen},
{\em Rattle: a ``velocity'' version of the Shake algorithm for molecular
  dynamics calculations}, J. Comput. Phys. {\bf 52}, 24--34, 1983.

\bibitem{baffico-bernard-maday-02} 
{\sc L. Baffico, S. Bernard, Y. Maday, G. Turinici, and G. Z\'erah}, 
{\em Parallel in time molecular dynamics simulations}, 
Phys. Rev. E {\bf 66}, 057701, 2002.

\bibitem{bal-03}
{\sc G. Bal}, {\em On the convergence and the stability of the
parareal algorithm to solve partial differential equations}, in
Domain decomposition methods in science and engineering,
Lect. Notes Comput. Sci. Eng. {\bf 40}, R. Kornhuber, R. Hoppe,
J. Périaux, O. Pironneau, O. Widlund and J. Xu eds., 
425--432, Springer Verlag, 2005.

\bibitem{bal-maday-02} {\sc G. Bal and Y. Maday},
{\em A parareal time discretization
for nonlinear PDE's with application to the pricing of an American
put}, in Recent developments in domain decomposition methods, 
Lect. Notes Comput. Sci. Eng. {\bf 23}, L.F. Pavarino and A. Toselli eds., 
189--202, Springer Verlag, 2002.

\bibitem{bal} {\sc G. Bal and Q. Wu}, {\em Symplectic parareal}, 
in Domain decomposition methods in science and engineering,
Lect. Notes Comput. Sci. Eng. {\bf 60}, 
U. Langer, M. Discacciati, D.E. Keyes, O.B. Widlund and W. Zulehner
eds., 401--408, Springer Verlag, 2008.

\bibitem{batchelor} {\sc D.B. Batchelor}, personnal communication.

\bibitem{BZ89} {\sc A. Bellen and M. Zennaro},
{\em Parallel algorithms for initial value problems for
nonlinear vector difference and differential equations},
J. Comput. Appl. Math. {\bf 25}, 341--350, 1989.

\bibitem{benettin} {\sc G. Benettin and A. Giorgilli}, {\em On the
Hamiltonian interpolation of near to the identity symplectic
mappings with application to symplectic integration algorithms},
J. Stat. Phys. {\bf 74}, 1117--1143, 1994.

\bibitem{Bur:book} {\sc K. Burrage},
{\em Parallel and sequential methods for ordinary differential equations},
Numerical Mathematics and Scientific Computation, Oxford Science
Publications, The Clarendon Press, Oxford University Press, New York,
1995. 

\bibitem{burrage_review} {\sc K. Burrage},
{\em Parallel methods for ODEs},
Advances in Computational Mathematics {\bf 7}(1-2), 1--3, 1997.

\bibitem{chartier93} {\sc P. Chartier and B. Philippe}, 
{\em A parallel shooting technique for solving dissipative ODE's},
Computing {\bf 51}(3-4), 209--236, 1993.

\bibitem{farhat06} {\sc C. Farhat, J. Cortial, C. Dastillung, and
H. Bavestrello}, {\em Time-parallel implicit integrators for the
near-real-time prediction of linear structural dynamic responses}, 
Int. J. Numer. Meth. Engng. {\bf 67}(5), 697--724, 2006.

\bibitem{ficher-hecht-maday-05}
{\sc P. Fischer, F. Hecht, and Y. Maday}, {\em A parareal in
time semi implicit approximation of the Navier Stokes equations},
in Domain decomposition methods in science and engineering,
Lect. Notes Comput. Sci. Eng. {\bf 40}, R. Kornhuber, R. Hoppe,
J. Périaux, O. Pironneau, O. Widlund and J. Xu eds., 
433--440, Springer Verlag, 2005. 

\bibitem{frenkel} {\sc D. Frenkel and B. Smit}, {\em
Understanding molecular simulation, from algorithms to
applications}, 2nd ed., Academic Press, 2002.

\bibitem{gander-vandewalle-05} {\sc M. Gander and S. Vandewalle},
{\em On the superlinear and linear convergence of the parareal
  algorithm}. Proceedings of the 16th International Conference on Domain
Decomposition Methods, 2005. 

\bibitem{gander-vandewalle-07} {\sc M. Gander and S. Vandewalle},
{\em Analysis of the parareal time-parallel time-integration
method}, SIAM J. Sci. Comput. {\bf 29}(2), 556--578, 2007.

\bibitem{garido} {\sc I. Garrido, B. Lee, G.E. Fladmark, and M.S. Espedal},
{\em Convergent iterative schemes for time parallelization}, Mathematics
of Computation {\bf 75}(255), 1403--1428, 2006. 

\bibitem{Hack} {\sc W. Hackbusch}, 
{\em Parabolic multigrid methods}, 
Computing methods in applied sciences and engineering VI (Versailles,
1983), 189--197, North-Holland, Amsterdam, 1984.

\bibitem{hairer97} {\sc E. Hairer and C. Lubich}, {\em The life span
of backward error analysis for numerical integrators},
Numer. Math. {\bf 76}, 441--462, 1997.

\bibitem{hlw} {\sc E. Hairer, C. Lubich, and G. Wanner}, 
{\em Geometric numerical integration: structure-preserving
algorithms for ordinary differential equations}, Springer Series in
Computational Mathematics {\bf 31}, 2002. 

\bibitem{joly2} {\sc P. Joly}, {\em Numerical methods for elastic wave
propagation}, in Waves in nonlinear pre-stressed materials,
M. Destrade and G. Saccomandi eds., 181--281, Springer-Verlag, 2007.

\bibitem{joly1} {\sc P. Joly}, {\em The mathematical model for elastic
wave propagation. Effective computational methods for wave
propagation}, in Numer. Insights, 247--266, Chapman \& Hall/CRC, 2008.

\bibitem{celestial1} {\sc J. Laskar}, {\em A numerical experiment on the
chaotic behavior of the Solar system}, Nature {\bf 338}, 237--238, 1989. 

\bibitem{celestial2} {\sc J. Laskar}, {\em Chaotic diffusion in the
Solar system}, Icarus {\bf 196}, 1--15, 2008.

\bibitem{leim_reich} {\sc B. Leimkuhler and S. Reich}, {\em Simulating
Hamiltonian dynamics}, Cambridge University Press, 2004.

\bibitem{rattle2} {\sc B.J. Leimkuhler and R.D. Skeel}, {\em Symplectic
numerical integrators in constrained Hamiltonian systems},
J. Comput. Phys. {\bf 112}, 117--125, 1994.

\bibitem{LRSV82} 
{\sc E. Lelarasmee, A.E. Ruehli, and A.L. Sangiovanni-Vincentelli},
{\em The waveform relaxation method for time-domain analysis of large
scale integrated circuits}, 
IEEE Trans. on CAD of IC and Syst. {\bf 1}, 131--145, 1982.

\bibitem{lions-maday-turinici-01}
{\sc J.-L. Lions, Y. Maday, and G. Turinici}, {\em A parareal in
time discretization of PDE's}, C. R. Acad. Sci. Paris, S\'erie I {\bf 332},
661--668, 2001.

\bibitem{maday-08} {\sc Y. Maday}, {\em The parareal in time
algorithm}, in Substructuring Techniques and Domain Decomposition
Methods, F. Magoul\`es ed., Chapter 2, pp.~19--44, Saxe-Coburg
Publications, Stirlingshire, UK, 2010. doi:10.4203/csets.24.2 

\bibitem{maday-turinici-02} {\sc Y. Maday and G. Turinici}, 
{\em A parareal in time procedure for the control of partial
differential equations},  
C. R. Acad. Sci. Paris, S\'erie I {\bf 335}, 387--392, 2002.

\bibitem{maday-turinici-03} {\sc Y. Maday and G. Turinici},
{\em A parallel in time approach for quantum control: the parareal
algorithm}, Int. J. Quant. Chem. {\bf 93}, 223--228, 2003.

\bibitem{maday-turinici-04} {\sc Y. Maday and G. Turinici}, 
{\em The parareal in time iterative solver: a further direction to
parallel implementation}, 
in Domain decomposition methods in science and engineering,
Lect. Notes Comput. Sci. Eng. {\bf 40}, R. Kornhuber, R. Hoppe,
J. Périaux, O. Pironneau, O. Widlund and J. Xu eds., 441--448, Springer
Verlag, 2005. 

\bibitem{Quarteroni} {\sc A. Quarteroni and A. Valli}, 
{\em Domain decomposition methods for partial differential equations}, 
Numerical Mathematics and Scientific Computation, Oxford Science
Publications, The Clarendon Press, Oxford University Press, New York,
1999.

\bibitem{reich99} {\sc S. Reich}, {\em Backward error analysis for
numerical integrators}, SIAM J. Numer. Anal. {\bf 36}, 1549--1570,
1999. 

\bibitem{shake} {\sc J.-P. Ryckaert, G. Ciccotti, and H.J.C. Berendsen},
{\em Numerical integration of the cartesian equations of motion of a
system with constraints: molecular dynamics of $n$-alkanes},
J. Comput. Phys. {\bf 23}, 327--341, 1977.

\bibitem{tremaine} {\sc P. Saha, J. Stadel and S. Tremaine},
{\em A parallel integration method for Solar system dynamics}, 
Astron. J. {\bf 114}(1), 409--414, 1997.

\bibitem{celestial3} {\sc P. Saha and S. Tremaine}, {\em Symplectic
integrators for solar system dynamics}, Astron. J. {\bf 104},
1633--1640, 1992. 

\bibitem{ss-calvo} {\sc J.M. Sanz-Serna and M.P. Calvo},
{\em Numerical Hamiltonian problems}, Chapman \& Hall, 1994.

\bibitem{staff-ronquist-03} {\sc G.A. Staff and E.M. R{\o}nquist},  
{\em Stability of the parareal algorithm}, 
in Domain decomposition methods in science and engineering,
Lect. Notes Comput. Sci. Eng. {\bf 40}, R. Kornhuber, R. Hoppe,
J. Périaux, O. Pironneau, O. Widlund and J. Xu eds., 
449--456, Springer Verlag, 2005.

\bibitem{Widlund} {\sc A. Toselli and O. Widlund},  
{\em Domain decomposition methods---algorithms and theory}, 
Springer Series in Computational Mathematics {\bf 34},
Springer Verlag, 2005. 

\end{thebibliography}
\end{document}